\documentclass[journal ]{new-aiaa}
\usepackage[utf8]{inputenc}
\usepackage{textcomp}

\usepackage{graphicx}
\usepackage{amsmath}
\usepackage[version=4]{mhchem}
\usepackage{siunitx}
\usepackage{longtable,tabularx}
\usepackage{mathtools}
\usepackage{mathrsfs}
\usepackage{xcolor}

\usepackage{tikz}
\graphicspath{{images/}}
\usepackage{verbatim}
\usetikzlibrary{decorations.pathmorphing}
\usetikzlibrary{shapes,arrows}
\usepackage[linesnumbered,ruled,vlined]{algorithm2e}

\usepackage{amsthm}
\usepackage{enumerate}
\usepackage{leftidx}

\usepackage{tabu}
\usepackage{arydshln}
\usepackage{booktabs,etoolbox}
\usepackage{multirow}
\usepackage{wasysym}
\usepackage{subfig}
\usepackage{caption}

\title{Passively-Safe and Robust Multi-Agent Optimal Control \\
with Application to Distributed Space Systems}

\author{Tommaso Guffanti \footnote{Postdoctoral Scholar, Aeronautics and Astronautics, 496 Lomita Mall; tommaso@stanford.edu.}\ and Simone D'Amico \footnote{Associate Professor, Aeronautics and Astronautics, 496 Lomita Mall; Associate Fellow AIAA.}\footnote{A preliminary version of this work was presented as Paper 1101 at the AIAA Scitech 2021 Forum, Virtual Event, January 11-15 \& 19-21.}} \affil{Stanford University, Stanford, California 94305}

\begin{document}

\maketitle

\begin{abstract}
This paper presents a novel methodology to enforce motion safety guarantees even in the event of a sudden loss of control capabilities by any agent within a multi-agent system. This passive safety methodology permits the replacement of point-evaluated constraints on the trajectory following a contingency with a single function of the integration constants of the equations of motion at contingency instants. The effects of uncertainties on the trajectory are compensated for through the method of variation of parameters. This permits a reduction in the number of required constraints by one polynomial degree in the number of discrete time samples, and allows computationally efficient enforcement of passive safety within a multi-agent optimal control problem, solvable using direct methods. The main application is distributed space systems employing miniaturized low size-weight-and-power and commercial-off-the-shelf technology, which reduce mission financial costs at the expense of reliability. Experimental results for the upcoming Virtual Super-resolution Optics with Reconfigurable Swarms (VISORS) mission, as well as complementary formation-flying test cases in eccentric orbits, demonstrate the advantages of the proposed methodology in achieving fault-tolerant safety guarantees, as well as computational efficiency.
\end{abstract}

\section{Introduction}
Distributed space systems (DSS) demand high levels of safety, robustness and fault-tolerance. This risk-aversion is mainly due to the potential consequences a collision in space may have in terms of debris generation and investment losses. In recent years, the trends of spacecraft bus miniaturization and use of low-size-weight-and-power (low-SWaP) as well as commercial-off-the-shelf (COTS) technologies have fostered the distribution of payload and tasks among multiple coordinated agents enabling functionalities that are otherwise not achievable by single monolithic systems \cite{damico_pavone_miniat}. On the other hand, these increased the likelihood of contingencies that directly or indirectly cause permanent or temporary inability to control a spacecraft \cite{jacklin2019small}. Fault-tolerant strategies are therefore fundamental to mitigate the consequences and keep the risks associated to flying these systems within tolerable levels \cite{dezfuli2011nasa}. From a collision safety perspective, different strategies are possible, including: reactive \cite{Scharf_11}, proactive \cite{starek_2017_fast}, and passive \cite{hubert1984stationkeeping, damico_proximity_2006, breger_08, barbee_11}. In general terms, reactive means that spacecraft have to re-compute trajectories on-line to avoid collision after a contingency. Proactive can be intended as a variant of reactive, in which the availability of a controlled escape trajectory/maneuver has to be guaranteed at any time even in the presence of a degradation, but not complete loss, of control capabilities. Finally, passive means that the controlled trajectories have been pre-designed to guarantee safe separation even in the event of sudden loss of control capabilities, being it temporary or permanent. To fly a real space mission, possibly all of these strategies are needed, and must be implemented within the guidance and control pipeline. Nevertheless, relying exclusively on a reactive or even proactive approach has limitations, such as not being robust to contingencies that do not permit trajectory re-planning and/or control on the short term of one or multiple spacecraft. These types of contingencies (e.g., sudden bus safe modes, sudden temporary or permanent failure of critical sub-systems) is increasingly common on miniaturized buses employing low-SWAP and COTS components \cite{jacklin2019small}, and demand for a passively-safe approach to fault-tolerant  motion safety. Looking at the literature, examples of reactive and proactive collision avoidance strategies are relatively recent \cite{Scharf_11, starek_2017_fast}. Whereas, the first  examples of passively-safe strategies date back to the 60's with the co-elliptic rendezvouses of the Apollo missions \cite{Woffinden_07}, and to the 80's with the colocation of geostationary satellites in shared longitude slots \cite{hubert1984stationkeeping}. Subsequently, in the last decades,  passive safety (PS) became an enabler of the first binary formation-flying missions in low-Earth orbit  such as GRACE, TANDEM-X, PRISMA and BIROS \cite{Errico_13}.  Where the relative orbits were designed  such that safe separation was always guaranteed in the radial/normal orbital directions, irrespective of failures of controlling and/or estimating the separation in the along-track direction \cite{damico_proximity_2006}. This work has been then generalized to $N$-spacecraft swarms in near-circular orbits \cite{koenig_2018_robust}. The limitations of this literature are  multiple. Firstly, except few cases in the close-proximity operations literature \cite{breger_08, starek_2017_fast}, fault-tolerant motion safety strategies are explored at orbit design level \cite{damico_proximity_2006, barbee_11, koenig_2018_robust} and not as a formal constraint within the optimal control problem (OCP). Secondly, PS guarantees are given under assumptions on the orbit type and dynamics, e.g., circularity \cite{breger_08, starek_2017_fast} or near-circularity \cite{damico_proximity_2006, koenig_2018_robust}, and linearized two-body dynamics \cite{breger_08, starek_2017_fast}, or including just specific perturbations \cite{koenig_2018_robust}. Thirdly, PS is enforced through possibly over-conservative constraints on the relative orbit design, e.g., energy matching conditions, sub-state designs \cite{damico_proximity_2006, koenig_2018_robust}. In general, these limitations cause the current state-of-the-art to fall short in addressing the needs required by future DSS missions, both in terms of applicability,  breadth and robustness.

 This paper explores the mathematical method of variation of parameters (VoP) \cite{VoP_book, vallado} to present a novel methodology for the modeling of dynamics-dependent constraints with application to fault-tolerant control, and to enforce them efficiently within an OCP solvable through direct methods. Specifically, VoP is a general method for the solution of both linear and nonlinear ordinary differential equations (ODE) \cite{VoP_book}, which models the effects of non-integrable dynamics on the integration constants (IC) of an integrable portion of the governing dynamics itself. VoP was firstly developed for the orbital dynamics \cite{vallado} under the form of the Gauss variational equations (GVE). GVE model the effect of non-integrable orbital perturbations and control accelerations on the Keplerian orbital elements (OE), which are IC of the homogeneous part of the orbital ODE. Similarly, the IC of the homogeneous part of the ODE governing the relative orbital dynamics linearized for small spacecraft separations are first-order equivalent to the relative orbital elements (ROE), combination of the OE of two spacecraft or a reference orbit and a spacecraft. This holds in both near-circular \cite{clohessy_terminal_1960, damico_ROEIC_2005} and eccentric orbits \cite{yamanaka_new_2002, sullivan_nonlinear_2017} by tuning appropriately the ROE definition.
Using VoP and isolating the effect of non-integrable dynamics on the IC state has various advantages. Firstly, it fosters the formulation of highly efficient analytical closed-form dynamics models \cite{gim_state_2003, koenig_2017_new, guffanti_long-term_2017, guffanti_linear_2019}. Secondly, it  permits direct compensation of the effects of non-integrable dynamics in guidance and control strategies \cite{chernick_phd_2021, koenig_2018_robust}. Thirdly, it  permits direct specification of the "long-term" behavior of the natural motion of the system, which is directly connected to the value of the IC themselves. This  is useful when constraints such as collision avoidance have to be enforced over a specified time horizon even in absence of applied control input. This is evident by looking at the DSS literature, where PS is easily obtained by properly selecting the OE or ROE of two resident space objects,  such as for co-elliptic rendezvous \cite{barbee_11}, colocation \cite{hubert1984stationkeeping},  and formation-flying \cite{damico_proximity_2006}. In addition, given the applicability of VoP to a wider range of dynamical systems \cite{VoP_book}, the advantages  brought by its use are possibly generalizable.

From an algorithmic perspective, this paper  explores the optimization and robotics literature in  five topics of motion planning and control.
The first topic is open-loop trajectory optimization under state and control constraints. In particular, this paper uses sequential convex programming (SCP), which is a local direct optimization method that seeks convergence to a local stationary point of a non-convex problem by successively convexifying the non-convex cost and constraints and solving a series of convexified problems \cite{scp_2000, liu_scp_2014, malyuta_scp_2022}.
The second topic is robust motion planning, to achieve probabilistic guarantees of constraints satisfaction  in the presence of system uncertainties. Usually the literature seeks either probabilistic guarantees of safety for a stochastic model of uncertainty or considers a bounded model of uncertainty \cite{singh_phd_2019}.  This paper models the uncertainty as linearized Gaussian \cite{carpenter2018navigation}, uses the Kalman filter equations to propagate the error covariance matrix through the dynamics linearized around the mean estimate, and computes safety margins by using the relevant diagonal components. Moreover, the effect of applied control on the uncertainty propagation is explicitly taken into account, being relevant to model the actuation uncertainty present, for example, on spacecraft thrusters.
The third topic is model predictive control (MPC) \cite{mpc_book}, to retrieve closed-loop performances on-line. In particular, in this paper the open-loop robust motion planner is used within a MPC pipeline to retrieve closed-loop performances on-line  in the presence of unmodeled or mismodeled effects. 
The fourth topic is multi-agent control, to deploy algorithms on multi-vehicle systems. In general, the computational complexity of checking separation between any  pair of agents within an SCP setting scales quadratically with their number \cite{morgan_2014}, this makes the solutions of the centralized multi-agent collision avoidance problem computationally intense for large number of agents. Possible mitigation strategies are  to limit the collision checks to a specified subset of neighboring agents \cite{morgan_2014}, or  to decouple the problem in a centralized open-loop planning and a decentralized closed-loop tracking  that keeps the agents within specified tracking error bounds \cite{singh_phd_2019}. In this paper the latter approach is followed.
Finally, fault-tolerance is a  topic largely addressed in automatic control literature, assuming diverse connotations in different application fields \cite{ZHANG_08}. Looking at the multi-agent control literature specifically, consensus schemes \cite{Kazerooni_10} are often used, and fault-tolerance is not included as a formal constraint within the optimization problem. Therefore, this paper addresses the relevant, and not exhaustively explored, problem of efficiently and explicitly enforcing fault-tolerant motion safety guarantees within a multi-agent OCP solvable through direct methods.

This paper builds over the previous work of the authors \cite{guffanti_2018_integration, guffanti_2021_multi, guffanti_2022_acc, guffanti_phd_2022} to provide five main contributions to the state-of-the-art. Firstly, in Section \ref{Sec_II},  it presents the novel methodology for the modeling of dynamics-dependent constraints with application to fault-tolerant control,  in the presence of nonlinear non-integrable dynamics and realistic system uncertainties from sensing, actuation, and unmodeled system dynamics.  This methodology permits replacing point-evaluated constraints on the trajectory following a contingency with a single function of the IC at contingency instants, while compensating for the effects of uncertainties through VoP. Secondly, in Section \ref{Sec_III}, it  demonstrates that the number of constraints required to guarantee multi-agent PS is polynomial (super-linear) in the number of time samples, and that the proposed methodology reduces them by one polynomial degree. Thirdly, in Section \ref{Sec_IV}, it  shows how this methodology can be used within a sequential optimization framework for the solution of the multi-agent passively-safe OCP. This includes a SCP-based centralized algorithm that enforces guarantees of multi-agent PS, and a decentralized tracking logic that reinforces them within an MPC-like closed-loop pipeline. Fourthly, in Section \ref{Sec_V}, it specializes the methodology to the DSS dynamics, and it derives novel closed-form solutions of PS in eccentric orbits, which can be used both for relative orbit design and within a constrained OCP. Finally, in Section \ref{Sec_VI}, it applies the  algorithmic framework to the Virtual Super-resolution Optics with Reconfigurable Swarms (VISORS) mission \cite{VISORS}, a first-of-a-kind nanosatellite segmented telescope  due to launch in 2024 for high-resolution imaging of the solar corona, and to complementary formation-flying test cases in eccentric orbits.  These test cases demonstrate the advantages of the proposed methodology in achieving fault-tolerant safety guarantees, as well as computational efficiency. Moreover, a comparison with closed-form impulsive control solutions with flight-heritage \cite{gaias_imp_2015} shows the potential of the proposed algorithmic framework in terms of gained control accuracy and fuel-efficiency, while enforcing the desired safety guarantees in closed-loop.

\section{The Integration Constants Approach}\label{Sec_II}
\subsection{Mathematical Preliminaries}\label{Sec_VoP}
 Consider a system characterized by an $s$-dimensional operational state $\boldsymbol{\chi}(t) \in \mathbb{R}^s$, where $t \in \mathbb{R}_{\geq 0}$ is the independent variable here identified with time. The  ordinary differential equation (ODE) that governs the operational state is assumed nonlinear control-affine as 
\begin{equation} \label{ODE_OP}
\dot{\boldsymbol{\chi}}(t) = \boldsymbol{f}(\boldsymbol{\chi}(t)) + \boldsymbol{B}(\boldsymbol{\chi}(t))\boldsymbol{u}(t)
\end{equation}
with $\boldsymbol{f} : \mathbb{R}^s \rightarrow \mathbb{R}^s$ the drift vector term, $\boldsymbol{u}(t) \in \mathbb{R}^r$ the $r$-dimensional applied control input, and $\boldsymbol{B} : \mathbb{R}^s \rightarrow \mathbb{R}^{s \times r}$ the control input matrix. Both $\boldsymbol{f}$ and $\boldsymbol{B}$ are assumed continuously differentiable. Moreover,  assume the drift vector term in Eq. \ref{ODE_OP} can be partitioned as
\begin{equation}\label{ODE_ASS}
\dot{\boldsymbol{\chi}}(t) = \boldsymbol{\iota}(\boldsymbol{\chi}(t)) + \boldsymbol{d}(\boldsymbol{\chi}(t)) + \boldsymbol{B}(\boldsymbol{\chi}(t))\boldsymbol{u}(t)
\end{equation}
where, the $\boldsymbol{\iota} : \mathbb{R}^s \rightarrow \mathbb{R}^s$, is a reduced part of the governing ODE integrable analytically in closed-form, and $\boldsymbol{d} : \mathbb{R}^s \rightarrow \mathbb{R}^s$ is the part of the governing ODE that makes it non-integrable. The reduced part $\dot{\boldsymbol{\chi}}(t) = \boldsymbol{\iota}(\boldsymbol{\chi}(t))$ can be integrated analytically starting at instant $t_i$ as
\begin{equation}\label{IC_map}
\boldsymbol{\chi}(t) = \boldsymbol{\psi}(t, \boldsymbol{c}(t_i))
\end{equation}
where, $\boldsymbol{c} \in \mathbb{R}^s$ is the integration constants (IC) state, and $\boldsymbol{\psi} : \mathbb{R}^{s+1} \rightarrow \mathbb{R}^s$ is the integrated dynamics map defined in closed-form upon integration, and is here assumed to be continuous, differentiable, invertible (with the inverse defined as $\boldsymbol{\psi}^{-1}(t, \boldsymbol{\chi}(t)) : \mathbb{R}^{s+1} \rightarrow \mathbb{R}^s$), and possible ambiguities for multivalued functions can be discriminated. The IC state is guaranteed to exist as long as Eq. \ref{ODE_ASS} holds, and its components are an independent linear or nonlinear combination of the operational state components at the integration instant $\boldsymbol{\chi}(t_i)$ (i.e., the initial conditions), and possibly depend on $t_i$ itself. The IC can have a clear physical meaning (e.g., conserved energy or angular momentuum), or, more simply, they can be the initial conditions of the system at integration instant. If the integrable dynamics is nonlinear, the map $\boldsymbol{\psi}$ is generally a nonlinear function of the IC. Therefore, its Jacobian $\boldsymbol{\Psi}(t, \boldsymbol{c}) = \partial\boldsymbol{\psi}(t, \boldsymbol{c})/\partial \boldsymbol{c} : \mathbb{R}^{s+1} \rightarrow \mathbb{R}^{s \times s}$ is a function of the IC, too. On the contrary, if the integrable dynamics is linear or linearized (i.e., $\boldsymbol{\iota}(\boldsymbol{\chi}(t)) \approx \boldsymbol{A}(t) \boldsymbol{\chi}(t)$, with $\boldsymbol{A} \in \mathbb{R}^{s \times s}$), the map $\boldsymbol{\psi}$ becomes linear with respect to the IC as $\boldsymbol{\psi}(t, \boldsymbol{c}) = \boldsymbol{\Psi}(t) \boldsymbol{c}$, where now its Jacobian coincides with the fundamental matrix $\boldsymbol{\Psi}(t) \in \mathbb{R}^{s \times s}$ solution of the linear ODE and just dependent on time. Note that the non-integrable part may include truncated higher-order terms if a linearization is performed to isolate a linear integrable part.

Variation of parameters (VoP) \cite{VoP_book} is a mathematical tool that  permits the modelling of the effect of non-integrable dynamics and control on the IC. In particular, Eq.  \ref{ODE_ASS} can be rewritten exactly as an ODE of IC variation as
\begin{equation}\label{VOP}
\dot{\boldsymbol{c}}(t) = \boldsymbol{d}_c(t, \boldsymbol{c}(t)) + \boldsymbol{B}_c(t,\boldsymbol{c}(t))\boldsymbol{u}(t)
\end{equation}
where, the drift vector is defined as $\boldsymbol{d}_c(t, \boldsymbol{c}(t)) = \boldsymbol{\Psi}^{-1}(t, \boldsymbol{c}(t))\boldsymbol{d}(\boldsymbol{c}(t)) : \mathbb{R}^{s+1} \rightarrow \mathbb{R}^s$, whereas control input matrix is defined as $\boldsymbol{B}_c(t, \boldsymbol{c}(t)) = \boldsymbol{\Psi}^{-1}(t, \boldsymbol{c}(t)) \boldsymbol{B}(\boldsymbol{c}(t)) : \mathbb{R}^{s+1} \rightarrow \mathbb{R}^{s \times r}$. Eq. \ref{VOP} shows how both non-integrable dynamics and control have a direct effect on the IC state. On the other hand, they do not alter the structure of the map $\boldsymbol{\psi}$. In particular, the motion trajectory of $\boldsymbol{\chi}(t)$ governed by Eq. \ref{ODE_ASS} is instantaneously contained and tangent to $\boldsymbol{\psi}(t, \boldsymbol{c}(t))$,  such that $\boldsymbol{\chi}(t) = \boldsymbol{\psi}(t, \boldsymbol{c}(t))$, with the evolution in time of $\boldsymbol{c}(t)$ governed by Eq. \ref{VOP}. From a topological perspective, $\boldsymbol{\psi}(t, \boldsymbol{c}(t))$ describes a submanifold of $\mathbb{R}^s$ on which $\boldsymbol{\chi}(t)$ lies instantaneously. This idea is known as the "osculation principle" \cite{vallado}.  As mentioned in the introduction and further detailed in Section \ref{Sec_V}, the assumption in Eq. \ref{ODE_ASS} and VoP apply both to the orbital dynamics and to the relative orbital dynamics linearized for small spacecraft separations.

\subsection{Efficient Modeling of Dynamics-Dependent Constraints with Application to Fault-Tolerant Control}\label{Sec_IC_Constr}
In this section, the advantage  brought by the IC representation regarding the modeling of dynamics-dependent constraints with application to fault-tolerant control is presented in general terms. Then, in Section \ref{Sec_III}, it is applied specifically to the modeling of multi-agent  passive safety (PS), and, in Section \ref{Sec_V}, it is specialized to the distributed space systems (DSS) dynamics. From a fault-tolerant perspective, there is interest in providing guarantees of constraints satisfaction  even in the event of a sudden loss of control capability. In particular,  consider the state-dependent scalar constraint $\mathcal{M}(\boldsymbol{\chi}(t)) \leq 0$, where $\mathcal{M} : \mathbb{R}^s \rightarrow \mathbb{R}$ is a generic non-convex scalar function of the state, assumed continuously differentiable. The constraint function is assumed  to be scalar for simplicity and because many higher dimensional constraints are decomposable into separated scalar ones. The objective is to enforce  this constraint at each instant in the time interval $[t_0, t_f]$, and in addition, at each instant in $(t_i, t_i+T]$, where $T \in \mathbb{R}_{\geq 0}$, considering possible loss of control capability $\forall t_i \in [t_0, t_f]$. Under  Eq. \ref{ODE_ASS}, this can be formulated in mathematical terms as
\begin{equation}\label{Eq_constr_fault}
\begin{aligned}
 \{ & \mathcal{M}(\boldsymbol{\chi}(t)) \leq 0, \hspace{0.5cm} \forall t \in [t_0, t_f],  \\
 & \textrm{such that:} \hspace{0.2cm} \dot{\boldsymbol{\chi}}(t) = \boldsymbol{\iota}(\boldsymbol{\chi}(t)) + \boldsymbol{d}(\boldsymbol{\chi}(t)) + \boldsymbol{B}(\boldsymbol{\chi}(t))\boldsymbol{u}(t) \}\\
 & \textrm{AND} \\
 \{ &  \mathcal{M}(\boldsymbol{\chi}(t)) \leq 0, \hspace{0.5cm} \forall t \in (t_i, t_i + T], \hspace{0.2cm} \forall t_i \in [t_0, t_f], \\ 
 & \textrm{such that:} \hspace{0.2cm} \dot{\boldsymbol{\chi}}(t) = \boldsymbol{\iota}(\boldsymbol{\chi}(t)) + \boldsymbol{d}(\boldsymbol{\chi}(t)) \}
 \end{aligned}
\end{equation}
In a discrete time setting, if the interval $[t_0, t_f]$ is discretized in $n$ samples, and the intervals $(t_i, t_i + T]$ in $m$ samples, and if contingencies are  accounted for at each of the $n$ samples, the total number of constraints to be enforced is $nm$. In the following, the IC representation and VoP are used to reduce  this required number of constraints. 

To do so,  assume an instant of worst case occurrence of $\mathcal{M}$ exists within $[t_i, t_i + T]$. In such a way, the satisfaction of the constraint at this instant is a sufficient condition for the satisfaction of the constraint at any other instant within the interval. Mathematically,  assuming $\mathcal{M}$ continuously differentiable,
\begin{equation}\label{Eq_t_star_ex}
\begin{aligned}
    \exists t^* \in [t_i, t_i + T]: \{ \mathcal{M}(\boldsymbol{\chi}(t^*)) \leq 0 \} \Longrightarrow \{ \mathcal{M}(\boldsymbol{\chi}(t)) \leq 0, \hspace{0.2cm} \forall t \in [t_i, t_i + T] \}
\end{aligned}
\end{equation}
The instant of worst case occurrence is defined as
\begin{equation}\label{Eq_t_star0}
  t^* = \arg{\max_{t \in [t_i, t_i + T]}{\mathcal{M}(\boldsymbol{\chi}(t))}}
\end{equation}
Moreover,  assume preliminarily that the non-integrable part of the dynamics is set to zero after a contingency together with the control.  Under this assumption, the dynamics after a contingency is fully integrable and, therefore it is expressible in closed-form using Eq. \ref{IC_map}, $\forall t \in [t_i, t_i + T]$. Under this assumption, Eq. \ref{Eq_t_star0} can be rewritten as
\begin{equation}\label{Eq_t_star1}
  t^* = \arg\max_{t \in [t_i, t_i + T]}{\mathcal{M}\big(\boldsymbol{\psi}(t, \boldsymbol{c}(t_i))\big)}
\end{equation}
where,  by definition, the IC state $\boldsymbol{c}(t_i)$ is a constant linear or nonlinear combination of $\boldsymbol{\chi}(t_i)$ and possibly $t_i$ itself. Note that in Eq. \ref{Eq_t_star1}, $\mathcal{M}$ is a closed-form expression of time, which is the optimization variable, whereas the IC state acts as a parameter vector in the optimization.  $\mathcal{M}$ is a generic non-convex function of the state, and in general, $\boldsymbol{\psi}$ is a non-convex function of time. Therefore, Eq. \ref{Eq_t_star1} is a non-convex parametric optimization problem, in the parameter vector $\boldsymbol{c}(t_i)$. By definition, the solution of a parametric optimization problem is a function of the parameters themselves. Therefore, by definition, $t^*$ is a function of the IC. This function is assumed identifiable as $t^* = t^*(\boldsymbol{c}(t_i)) : \mathbb{R}^{s} \rightarrow \mathbb{R}_{\geq 0}$. Under the assumption of integrability of the dynamics  after a contingency, this  permits rewriting the sufficient condition in Eq. \ref{Eq_t_star_ex} as
\begin{equation}\label{Eq_constr1}
 \mathcal{M}(\boldsymbol{\chi}(t^*)) = \mathcal{M}\Big(\boldsymbol{\psi}\big(t^*\left(\boldsymbol{c}(t_i)\right), \boldsymbol{c}(t_i) \big)\Big) = \mathcal{M}^{*}(\boldsymbol{c}(t_i)) \leq 0
\end{equation}
which makes the constraint a function of the IC state evaluated at the contingency instant independent of the time along the uncontrolled trajectory: $\mathcal{M}^{*}(\boldsymbol{c}(t_i)) : \mathbb{R}^s \rightarrow \mathbb{R}$.  By removing the dependency on the time along the uncontrolled trajectories, fault-tolerant constraint satisfaction is achieved for at least $T$ after a contingency by enforcing as few constraints as the number of possible contingency instants. In particular, enforcing Eq. \ref{Eq_constr1} in place of Eq. \ref{Eq_constr_fault} implies a reduction of the number of constraints by one polynomial degree in the number of discrete time samples, from $nm$ to $n$. In Eq. \ref{Eq_constr1}, the formalization of the function $\mathcal{M}^{*}$ assumes the capability of identifying the function $t^*(\boldsymbol{c}(t_i))$. In this paper, a general algorithmic procedure is proposed in Section \ref{Sec_IV}, moreover, closed-form solutions applied to DSS PS are presented in Section \ref{Sec_V}. Furthermore,  Eq. \ref{Eq_constr1} assumes the integrability of the dynamics  after a contingency. Nevertheless, VoP can be used to make  it robust to non-integrable dynamics effects. This is done by bounding the integrated effect of non-integrable dynamics on $\mathcal{M}^{*}$ over $[t_i , t_i + T]$ as
\begin{equation}\label{Eq_bound_d}
 \left| \frac{\partial \mathcal{M}^{*}(\boldsymbol{c})}{\partial \boldsymbol{c}}\bigg\lvert_{\boldsymbol{c}(t_i)} \int_{t_i}^{t_i + T} \boldsymbol{d}_c(\tau, \boldsymbol{c}(\tau)) d\tau \right| \leq \beta(t_i, T)
\end{equation}
In Eq. \ref{Eq_bound_d}, $\partial \mathcal{M}^{*} / \partial \boldsymbol{c} \in \mathbb{R}^s$ is a Jacobian vector mapping the integrated effect of non-integrable dynamics to a net variation of $\mathcal{M}^{*}$, whereas, $\beta \in \mathbb{R}_{\geq 0}$ is a bounding constant. This constant is added to Eq. \ref{Eq_constr1} as 
\begin{equation}\label{Eq_constr2}
    \mathcal{M}^{*}(\boldsymbol{c}(t_i)) + \beta(t_i, T) \leq 0
\end{equation}
to compensate for non-integrable dynamics effects on constraint satisfaction.  Eq. \ref{Eq_bound_d} represents a first-order sensitivity model, which, depending on the nature of the function $\mathcal{M}^{*}(\boldsymbol{c})$, may lose accuracy the further away from the neighborhood of $\boldsymbol{c}(t_i)$, i.e., the larger the integrated effect of non-integrable dynamics on the IC state gets over $[t_i, t_i + T]$. If needed, this first-order model can be upgraded to a higher-order one, without impacting the proposed constraint reformulation advantage. As analyzed in Section \ref{Sec_VI}, the proposed first-order model is sufficiently accurate for the DSS PS application. In the following, the deterministic analysis developed in this section is generalized to account for bounded effects of uncertainties on the constraint satisfaction.

\subsection{Inclusion of Uncertainties}
 Consider now a closed-loop control setting where the system has to find a control plan $\tilde{\boldsymbol{u}}(t) \in \mathbb{R}^{r}$, for $t \in [t_0, t_f]$,  to move from an initial configuration at $t_0$ to a target one at $t_f$, acquiring at instants $t_e \in [t_0, t_f]$ an updated estimate of its IC state $\tilde{\boldsymbol{c}} (t_e) \in \mathbb{R}^s$. This estimate is the mean of a Gaussian distribution $\mathcal{N}(\tilde{\boldsymbol{c}} (t_e), \boldsymbol{C} (t_e))$, provided by on-board sensing and navigation, where $\boldsymbol{C} (t_e) \in \mathbb{R}^{s \times s}$ is the error covariance matrix. Note that, depending on  the estimation filter implementation, the Gaussian distribution may be provided on the operational state $\mathcal{N}(\tilde{\boldsymbol{\chi}} (t_e), \boldsymbol{X} (t_e))$, with $\tilde{\boldsymbol{\chi}} (t_e) \in \mathbb{R}^s$ and $\boldsymbol{X} (t_e) \in \mathbb{R}^{s \times s}$, and then mapped to the IC state inverting the integrated dynamics closed-form map as $\tilde{\boldsymbol{c}} (t_e) = \boldsymbol{\psi}^{-1}(t_e, \tilde{\boldsymbol{\chi}} (t_e))$ and $\boldsymbol{C} (t_e) = \boldsymbol{\Psi}^{-1}(t_e, \tilde{\boldsymbol{c}} (t_e)) \boldsymbol{X} (t_e) \boldsymbol{\Psi}^{-T}(t_e, \tilde{\boldsymbol{c}} (t_e))$, where $(.)^T$ means transpose. Moreover,  assume the system has available on-board a nonlinear dynamics model, $\tilde{\boldsymbol{d}}_c : \mathbb{R}^{s+1} \rightarrow \mathbb{R}^s$, approximating the real effect of non-integrable dynamics on the IC. Using this model and VoP, the system propagates on-board the estimated mean IC state as
\begin{equation}\label{ODE_GOV_mean_IC}
\dot{\tilde{\boldsymbol{c}}}(t) = \tilde{\boldsymbol{d}}_c(t, \tilde{\boldsymbol{c}}(t)) + \boldsymbol{B}_c(t,\tilde{\boldsymbol{c}}(t))\tilde{\boldsymbol{u}}(t)
\end{equation}
Given this on-board dynamics,  assume the system uses a linearized approach of propagation of uncertainties. The evolution in time of the Gaussian approximation of the uncertainty distribution is modeled through a linear ODE governing the state error covariance matrix \cite{carpenter2018navigation}

\begin{equation}\label{ODE_GOV_cov_IC}
\begin{aligned}
\dot{\boldsymbol{C}}(t) = & \tilde{\boldsymbol{D}}_c(t, \tilde{\boldsymbol{c}}(t), \tilde{\boldsymbol{u}}(t)) \boldsymbol{C}(t) + \boldsymbol{C}(t) \tilde{\boldsymbol{D}}_c^T(t,\tilde{\boldsymbol{c}}(t), \tilde{\boldsymbol{u}}(t)) + \boldsymbol{B}_c(t,\tilde{\boldsymbol{c}}(t)) \boldsymbol{U}(t) \boldsymbol{B}_c^T(t,\tilde{\boldsymbol{c}}(t)) + \boldsymbol{Q}(t) + \boldsymbol{E}(t) \end{aligned}
\end{equation}
which holds accurately in the neighborhood of the mean state estimate. In particular, $\tilde{\boldsymbol{D}}_c = \partial \tilde{\boldsymbol{d}}_c/\partial \boldsymbol{c}\big\lvert_{\tilde{\boldsymbol{c}}} + \partial \boldsymbol{B}_c/\partial \boldsymbol{c}\big\lvert_{\tilde{\boldsymbol{c}}} \tilde{\boldsymbol{u}} : \mathbb{R}^{s+r+1} \rightarrow \mathbb{R}^{s \times s}$ is a Jacobian matrix evaluated at the mean state estimate. $\boldsymbol{Q}(t) \in \mathbb{R}^{s \times s}$ is the power spectral density of a zero-mean white Gaussian process-noise modeling the discrepancy between the real non-integrable dynamics and the available approximation: $\tilde{\boldsymbol{d}}_c(t, \tilde{\boldsymbol{c}}(t)) - \boldsymbol{d}_c(t, \boldsymbol{c}(t))$. The process noise power spectral density and covariance may be time-dependent, but no time correlation is assumed to exist. $\boldsymbol{U}(t) \in \mathbb{R}^{r \times r}$ is the power spectral density of a zero-mean white Gaussian noise modeling the discrepancy introduced by the actuation system between the planned control input and the real applied control input: $\tilde{\boldsymbol{u}}(t) - \boldsymbol{u}(t)$. Also in this case no time correlation is assumed to exist. On the other hand, the relation between planned control input and associated execution uncertainty is modeled as: $\boldsymbol{U}(t) = \boldsymbol{\mathscr{U}}(\tilde{\boldsymbol{u}}(t)) + \bar{\boldsymbol{U}}(t)$ where $\boldsymbol{\mathscr{U}} : \mathbb{R}^r \rightarrow \mathbb{R}^{r \times r}$ and $\bar{\boldsymbol{U}}(t) \in \mathbb{R}^{r \times r}$ is a remaining aleatory component independent of the planned control. $\boldsymbol{E}(t) = \delta(t - t_{\tilde{e}}) \boldsymbol{\mathcal{E}}(\boldsymbol{C} (t_{\tilde{e}}))$ is an impulsive function modeling the effect of future expected measurements on the state error covariance, where the Dirac-delta equals one at the future predicted measurements instants $t_{\tilde{e}}$, and $\boldsymbol{\mathcal{E}} : \mathbb{R}^{s \times s} \rightarrow \mathbb{R}^{s \times s}$ is a function of the predicted future estimated state covariance $\boldsymbol{C} (t_{\tilde{e}})$. 
Modeling the effect to future measurements on the covariance is important, since it informs the controller of the expected closed-loop behavior of the covariance even when it has to solve for an open-loop plan. This prevents possible unfeasibility of the open-loop control solution, due to the growth of the uncertainty propagated in open-loop.

\subsubsection{Robustification of Fault-Tolerant Constraint Formulation}
It is of interest to make Eq. \ref{Eq_constr2} applicable in this real closed-loop control setting, robustifying it against the modeled uncertainties evaluated at a specified $q$-$\sigma$ confidence level (with $q \in \mathbb{N}_{>0}$ and $\sigma$ being a statistical standard deviation). To do so, the integrated effect of non-integrable dynamics and uncertainty on $\mathcal{M}^{*}$ over $[t_i , t_i + T]$ are bounded using VoP as
\begin{equation}\label{Eq_dist_eff}
\begin{aligned}
&\left| \frac{\partial \mathcal{M}^{*}(\boldsymbol{c})}{\partial \boldsymbol{c}}\bigg\lvert_{\tilde{\boldsymbol{c}}(t_i)} \int_{t_i}^{t_i + T} \tilde{\boldsymbol{d}}_{c}(\tau, \tilde{\boldsymbol{c}}(\tau)) d\tau \right| \\
& + \left| \frac{\partial \mathcal{M}^{*}(\boldsymbol{c})}{\partial \boldsymbol{c}}\bigg\lvert_{\tilde{\boldsymbol{c}}(t_i)} q \sqrt{ \left| Diag\left( \int_{t_i}^{t_i + T} \left( \tilde{\boldsymbol{D}}_c(\tau, \tilde{\boldsymbol{c}}(\tau)) \boldsymbol{C}(\tau) + \boldsymbol{C}(\tau) \tilde{\boldsymbol{D}}_c^T(\tau, \tilde{\boldsymbol{c}}(\tau)) + \boldsymbol{Q}(\tau) \right) d\tau + \boldsymbol{C}(t_i) \right) \right|} \right| \leq \beta(t_i, T, q)
\end{aligned}
\end{equation}
where $\partial \mathcal{M}^{*} / \partial \boldsymbol{c} \in \mathbb{R}^s$ is a Jacobian vector, the $Diag(.)$ operator extracts the diagonal vector from a square matrix, the absolute value $|.|$ is applied element-wise for vector quantities, $q \in \mathbb{N}_{>0}$ models the $q$-$\sigma$ uncertainty level considered, and $\tilde{\boldsymbol{c}}(t_i)$ and $\boldsymbol{C}(t_i)$ are obtained by propagating Eq. \ref{ODE_GOV_mean_IC}-\ref{ODE_GOV_cov_IC} up to the contingency instant $t_i$. In Eq. \ref{Eq_dist_eff}, $\beta \in \mathbb{R}_{\geq 0}$ is an enlarged bounding constant with respect to the one defined in Eq. \ref{Eq_bound_d}, which  allows the compensation of both non-integrable dynamics and uncertainty effects on the fault-tolerant constraint satisfaction, as
\begin{equation}\label{Eq_constr4}
  \mathcal{M}^{*}(\tilde{\boldsymbol{c}}(t_i)) + \beta(t_i, T, q) \leq 0
\end{equation}
where now the function $\mathcal{M}^{*}$ is evaluated at the mean IC state estimate.  Eq. \ref{Eq_dist_eff} represents a first-order sensitivity model, which, depending on the nature of the function $\mathcal{M}^{*}(\boldsymbol{c})$, may lose accuracy the further away from the neighborhood of $\boldsymbol{c}(t_i)$, i.e., the larger the integrated effect of non-integrable dynamics and uncertainty on the IC state gets over $[t_i, t_i + T]$. If needed, this first-order model can be upgraded to a higher-order one, without impacting the proposed constraint reformulation advantage. As analyzed in Section \ref{Sec_VI}, the proposed first-order model is sufficiently accurate for the DSS PS application.

\subsection{Efficient Formulation of the Fault-Tolerant Optimal Control Problem} \label{Sec_NumOpt}
 The optimal control problem (OCP) the system must solve in closed-loop, possibly at every instant $t_e \in [t_0, t_f]$ an updated state estimate is available, until a specified target set $\boldsymbol{\mathcal{C}}_{f}$ is reached at $t_f$, is
\begin{equation}\label{OCP_IC}
\begin{aligned}
& \underset{\tilde{\boldsymbol{u}}(t)}{\text{minimize}} & & \int_{t_e}^{t_f} \mathcal{J}(\tilde{\boldsymbol{c}}(t), \boldsymbol{C}(t), \tilde{\boldsymbol{u}}(t), t) dt &\\
& \text{subject to} & & \dot{\tilde{\boldsymbol{c}}}(t) = \tilde{\boldsymbol{d}}_c(t, \tilde{\boldsymbol{c}}(t)) + \boldsymbol{B}_c(t,\tilde{\boldsymbol{c}}(t))\tilde{\boldsymbol{u}}(t) & \forall t \in [t_e, t_f] \\
&&& \dot{\boldsymbol{C}}(t) = \tilde{\boldsymbol{D}}_c(t, \tilde{\boldsymbol{c}}(t), \tilde{\boldsymbol{u}}(t)) \boldsymbol{C}(t) + \boldsymbol{C}(t) \tilde{\boldsymbol{D}}_c^T(t,\tilde{\boldsymbol{c}}(t), \tilde{\boldsymbol{u}}(t)) + & \forall t \in [t_e, t_f] \\
&&& \hspace{1.2cm} + \boldsymbol{B}_c(t,\tilde{\boldsymbol{c}}(t)) \boldsymbol{U}(t) \boldsymbol{B}_c^T(t,\tilde{\boldsymbol{c}}(t)) + \boldsymbol{Q}(t) + \boldsymbol{E}(t) & \\
&&& \mathcal{M}^{*}(\tilde{\boldsymbol{c}}(t)) + \beta(t, T, q) \leq 0 & \forall t \in [t_e, t_f] \\
&&& \tilde{\boldsymbol{u}}(t) \in \boldsymbol{\mathcal{U}} & \forall t \in [t_e, t_f] \\
&&&  \tilde{\boldsymbol{c}}(t_f) \in \boldsymbol{\mathcal{C}}_{f}  & \\
&&&   \mathcal{N}(\tilde{\boldsymbol{c}}(t_e), \boldsymbol{C}(t_e)) \equiv \textrm{state estimate} & \\
\end{aligned}
\end{equation}
Where, $\mathcal{J} \in \mathbb{R}_{> 0}$ is a generic cost functional modeling  for example control input and/or time minimization, $\boldsymbol{\mathcal{U}}$ is the set of admissible control inputs, and the bound $\beta(t_i, T, q)$ is defined as in Eq. \ref{Eq_dist_eff}. The control horizon progressively contracts from $[t_0, t_f]$ to zero, while $t_e \rightarrow t_f$. This represents a generic MPC-like implementation that  permits the retrieval of closed-loop properties by recomputing the control plan when an updated state estimate is available.  In the last line of Eq. \ref{OCP_IC}, the state estimate available at instant $t_e$ is the initial condition of the optimization problem solved at that instant. The boundary condition at $t_f$ can be enforced either deterministically on the mean state, or probabilistically on the distribution $\mathcal{N}(\tilde{\boldsymbol{c}}(t_f), \boldsymbol{C}(t_f))$. In the second to last line in Eq. \ref{OCP_IC}, the former approach is followed. The OCP in Eq. \ref{OCP_IC} enforces fault-tolerant constraint satisfaction for at least $T$  after a contingency, at $q$-$\sigma$ confidence, at the cost of checking as many constraints as the number of possible contingency instants in $[t_e, t_f]$. This implies a reduction of the number of constraints by one polynomial degree in the number of discrete time samples with respect to the formulation in Eq. \ref{Eq_constr_fault} (i.e., from $nm$ to $n$).  Eq. \ref{OCP_IC} is a non-convex optimization problem  that can be solved, for example, using direct methods, as discussed in Section \ref{Sec_IV}.

\section{Multi-Agent Passively-Safe Optimal Control}\label{Sec_III}
\begin{figure}[ht!]
    \centering
    \subfloat[Collision avoidance \label{Fig_CAvsPS_1}]{\includegraphics[width=0.34\linewidth]{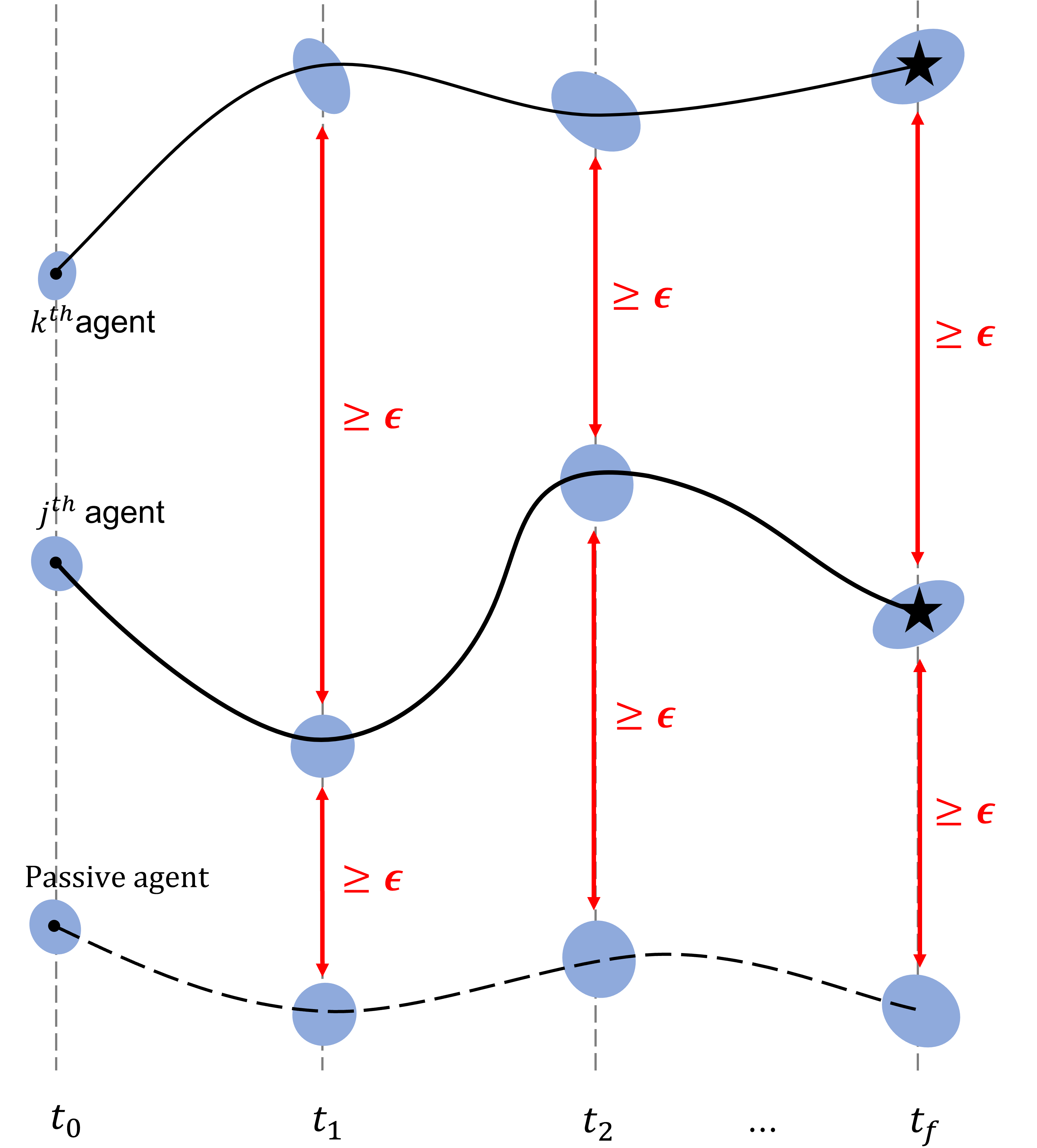}}\hspace{1.5cm}
    \subfloat[ Passive safety as in Definition 1 \label{Fig_CAvsPS_2}]{\includegraphics[width=0.47\linewidth]{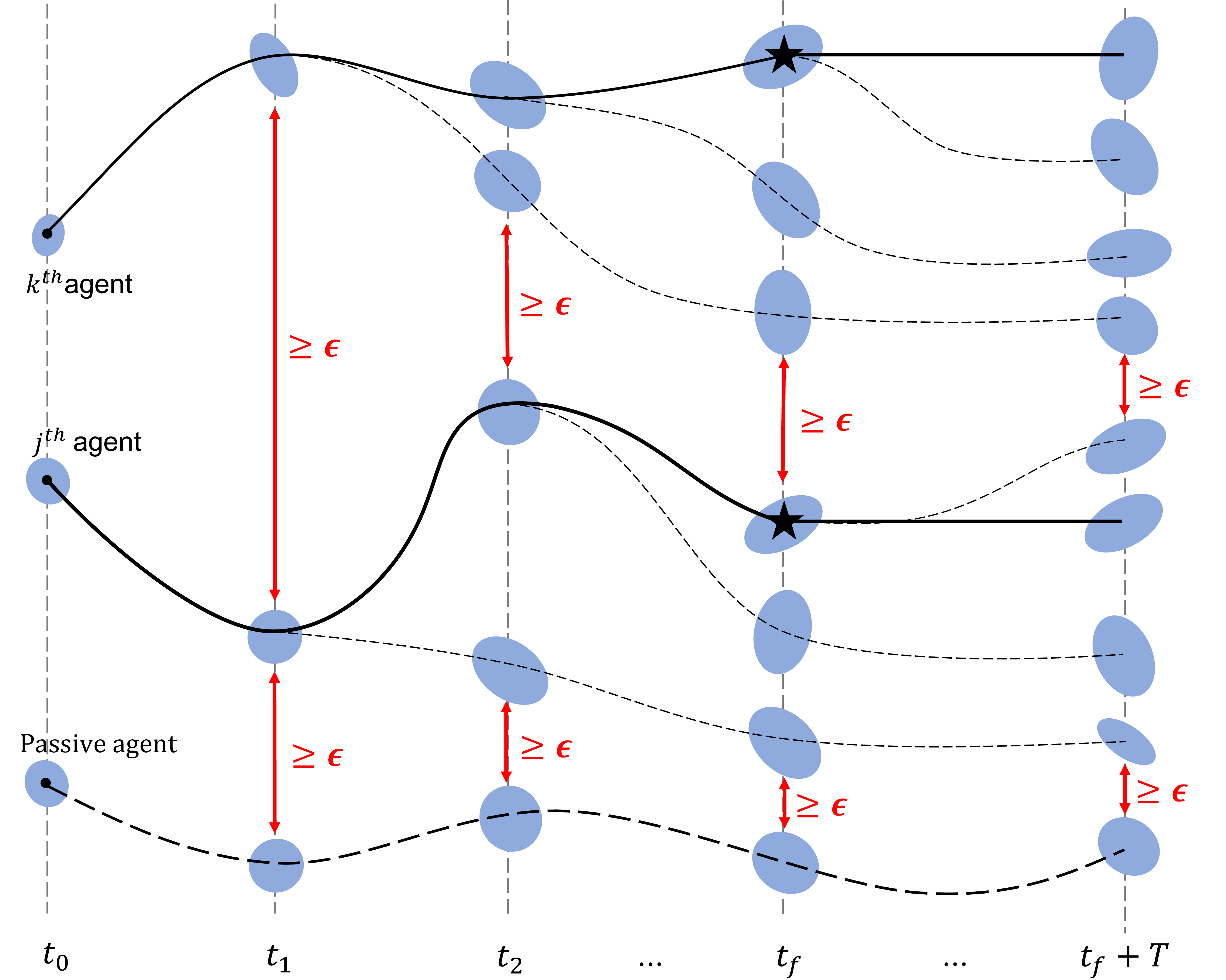}}
     \caption{Comparison between standard collision avoidance and passive safety}\label{Fig_CAvsPS}
\end{figure}

In this section, the  developed methodology is applied to the constraint of multi-agent  passive safety (PS), which is defined in this paper as follows. \\
 \textit{Definition 1} -- (Passive safety on $[t_0, t_f]$, for at least $T$, of at least $\epsilon$, at $q$-$\sigma$ confidence): Along the trajectory towards the target condition, a multi-agent system is passively-safe if all the agents maintain a separation of at least $\epsilon$, accounting for $q$-$\sigma$ uncertainties,  even in the event that any agent loses control at any instant $t \in [t_0, t_f]$ without possibility of re-planning until $t_f+T$. \\
Note that, in principle, it is possible to distinguish between different PS classes depending on the type of failure causing the loss of control capabilities. For example, it is possible to distinguish between hard failures, which prevent re-planning capabilities on the short term, or temporary failures, which allow re-planning capabilities on the short term. Furthermore, in the multi-agent setting, it is possible to just  account for the loss of control of a specified subset of agents. PS as in  Definition 1, is a rather conservative class, in which any subset or combination of agents may lose control at any instant in $[t_0, t_f]$, without being able to recover it and re-plan to avoid collision until $t_f + T$. The intention is to set the bar for milder PS classes  deducible by down-grading the one considered in this paper. Note that  Definition 1 directly accounts for the presence of possible "permanently passive" agents, whose trajectory is uncontrolled since $t_0$.  In a distributed space system, these may be a subset of spacecraft orbiting passively.  Moreover, the presented definition relies on the capability to retrieve an estimate of active and permanently passive agents states through sensing, navigation and prediction. A graphical comparison between this PS class and the standard collision avoidance is presented in Figure \ref{Fig_CAvsPS}. The controlled trajectories of the active agents are in solid lines, the uncontrolled ones in dashed lines, the shaded blue ellipses represent the uncertainty distributions of the Cartesian positions evaluated at $q$-$\sigma$ confidence interval.

In order to describe PS mathematically,  consider a multi-agent system in the Cartesian space,  composed of $N$ agents, of which $N_a$ are active and controlled, whereas $N_p$ are permanently passive from $t_0$. By definition, $N_a + N_p = N$. The operational state $\boldsymbol{\chi}_k(t)$ of each $k^{th}$ agent includes at the top the Cartesian position of the agent's center of mass, and is governed by Eq. \ref{ODE_OP}. For the permanently passive agents, $\boldsymbol{u}(t) = \boldsymbol{0}$, $\forall t \geq t_0$. Whereas, for active agents, control goes to zero  after a contingency.  In particular, define with $\boldsymbol{\chi}_k(t| t_i)$ the operational state of $k^{th}$ agent on the uncontrolled trajectory generated at instant $t_i$. If $k^{th}$ agent is active, $t_i$ is the contingency instant; if it is permanently passive $t_i = t_0$ by definition. Furthermore,  define as $\mathcal{S}_{kjil}(t) = ||\boldsymbol{S} (\boldsymbol{\chi}_k(t| t_i) - \boldsymbol{\chi}_j(t| t_l))||_2$ the separation reached at time $t$ by $k^{th}$ and $j^{th}$ agents on the uncontrolled trajectories generated at instants $t_i$ and $t_l$ respectively. The matrix $\boldsymbol{S} \in \mathbb{N}^{3 \times s}$ selects the components of the Cartesian position within the state, and $||.||_2$ is the $L_2$-norm. 
Reasoning deterministically, PS can be formulated by using Eq. \ref{Eq_constr_fault}, with the $\mathcal{M}$ function defined as
\begin{equation}\label{PS_M}
\mathcal{M}\big(\boldsymbol{\chi}_k(t), \boldsymbol{\chi}_j(t)\big) = \epsilon - \mathcal{S}_{kjil}(t) = \epsilon - ||\boldsymbol{S} (\boldsymbol{\chi}_k(t| t_i) - \boldsymbol{\chi}_j(t| t_l))||_2
\end{equation}
$\forall t_i \in [t_0, t_f]$, $\forall t_l \in [t_0, t_f]$, $\forall t \in [\max\{t_i, t_l\}, t_f + T]$, $k = 1, \dots N - 1$, $j = k + 1, \dots N$. Where, for example, if $j^{th}$ agent is permanently passive then $t_l = t_0$. Using the $\mathcal{M}$ formulation defined in Eq. \ref{PS_M}, the computational requirements in terms of number of constraints to be enforced are reported in Table \ref{tab_PS} on the top, and formal proof is developed in Appendix \ref{App:PS}. The interval $[t_0, t_f]$ is discretized in $n$ samples, whereas the interval $(t_f, t_f+T]$ in $m$ samples. A distinction between "continuous" and "impulsive" control types is made.  This distinction is motivated by the fact that failure to apply control can occur just when control application is actually planned. In particular, in the continuous type, control  is assumed to be planned and possibly lost at each of the $n$ sample instants in $[t_0, t_f]$. In the impulsive type, control/maneuvers  are assumed to be planned and possibly failed just in a subset $n_{M} \leq n$ of the $n$ instants,  where for simplicity the ratio $n/n_{M} \in \mathbb{N}_{>0}$. In both types, safe separation must be enforced at each of the $n$ and $m$ instants. This guarantees that even the uncontrolled arcs in between impulsive maneuvers are collision free. The number of constraints required to enforce PS between active agents are reported in Table \ref{tab_PS} on the left, whereas between active and permanently passive agents on the right. To summarize, using the $\mathcal{M}$ formulation, for both control types, the number of constraints to be enforced is cubic in the number of discrete time samples for PS between active agents, whereas  it is quadratic in the number of discrete time samples for PS between active and permanently passive agents. For $n_M \rightarrow n$, the number of constraints required in impulsive control tends to the ones required in continuous control, without affecting the polynomial degree behavior. The number of constraints scales quadratically with the number of agents. 
\begin{table}[ht!]
\centering
\caption{Passive safety constraint requirements, numerical examples are shown in Figure \ref{Fig_constr_num}}
\resizebox{0.98\textwidth}{!}{
\begin{tabular}{c c c c c c c}
 \hline 
 \hline
  \multicolumn{1}{c}{Formulation} & \multicolumn{1}{c}{Control type} & \multicolumn{2}{c}{Constraints Between Active Agents} & & \multicolumn{2}{c}{Constraints Between Active and Permanently Passive Agents} \\
  \cline{3-4} \cline{6-7}
  \multicolumn{1}{c}{} & \multicolumn{1}{c}{} & $\#$ & $\mathcal{O}(.)$ &  & $\#$ & $\mathcal{O}(.)$ \\ \hline
   $\mathcal{M}$ & Continuous & $\Big\{\frac{n(n+1)(2n+1)}{6}+(n+1)^2m\Big\} \frac{N_a(N_a-1)}{2}$ & $(n+m) n^2 N_a^2$ &  & $\Big\{\frac{n(n+1)}{2}+(n+1)m\Big\} N_a N_p$ & $(n+m) n N_a N_p$ \\
   & Impulsive & $\{\frac{n}{n_{M}}(\frac{n_{M}(n_{M}+1)(2n_{M}+1)}{6}) + (n_{M}+1)^2 m\} \frac{N_a(N_a-1)}{2}$ & $(n+m) n_{M}^2 N_a^2$ & & $\{\frac{n}{n_{M}}(\frac{n_{M}(n_{M}+1)}{2})+(n_{M}+1)m\} N_a N_p$ & $(n+m) n_{M} N_a N_p$ \\
  \cline{2-7}
  $\mathcal{M}^*$ &  Continuous & $n^2 \frac{N_a(N_a-1)}{2}$ & $n^2 N_a^2$ &  & $n N_a N_p$ & $n N_a N_p$ \\
  & Impulsive & $n_M^2 \frac{N_a(N_a-1)}{2}$ & $n_{M}^2 N_a^2$ & & $n_M N_a N_p$ & $n_{M} N_a N_p$ \\
  \hline
  \hline
\end{tabular}}\label{tab_PS}
\end{table}

\begin{figure}[ht!]
    \centering
    \subfloat[ Passive safety between active agents, continuous control \label{Fig_constr_num_1}]{\includegraphics[width=0.32\linewidth,trim=2 2 2 4,clip]{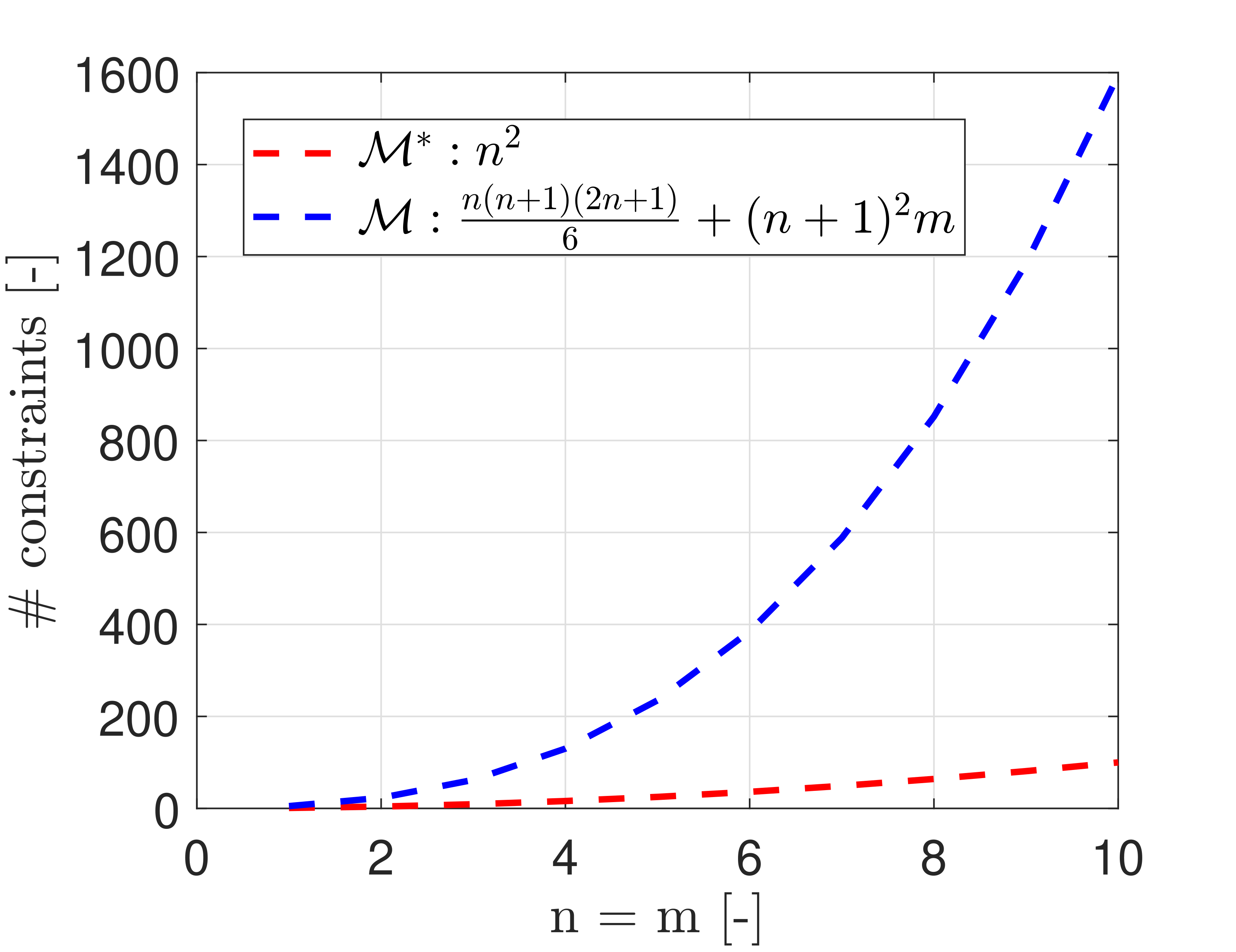}}
    \subfloat[ Passive safety between active agents, impulsive control \label{Fig_constr_num_2}]{\includegraphics[width=0.32\linewidth,trim=2 2 2 4,clip]{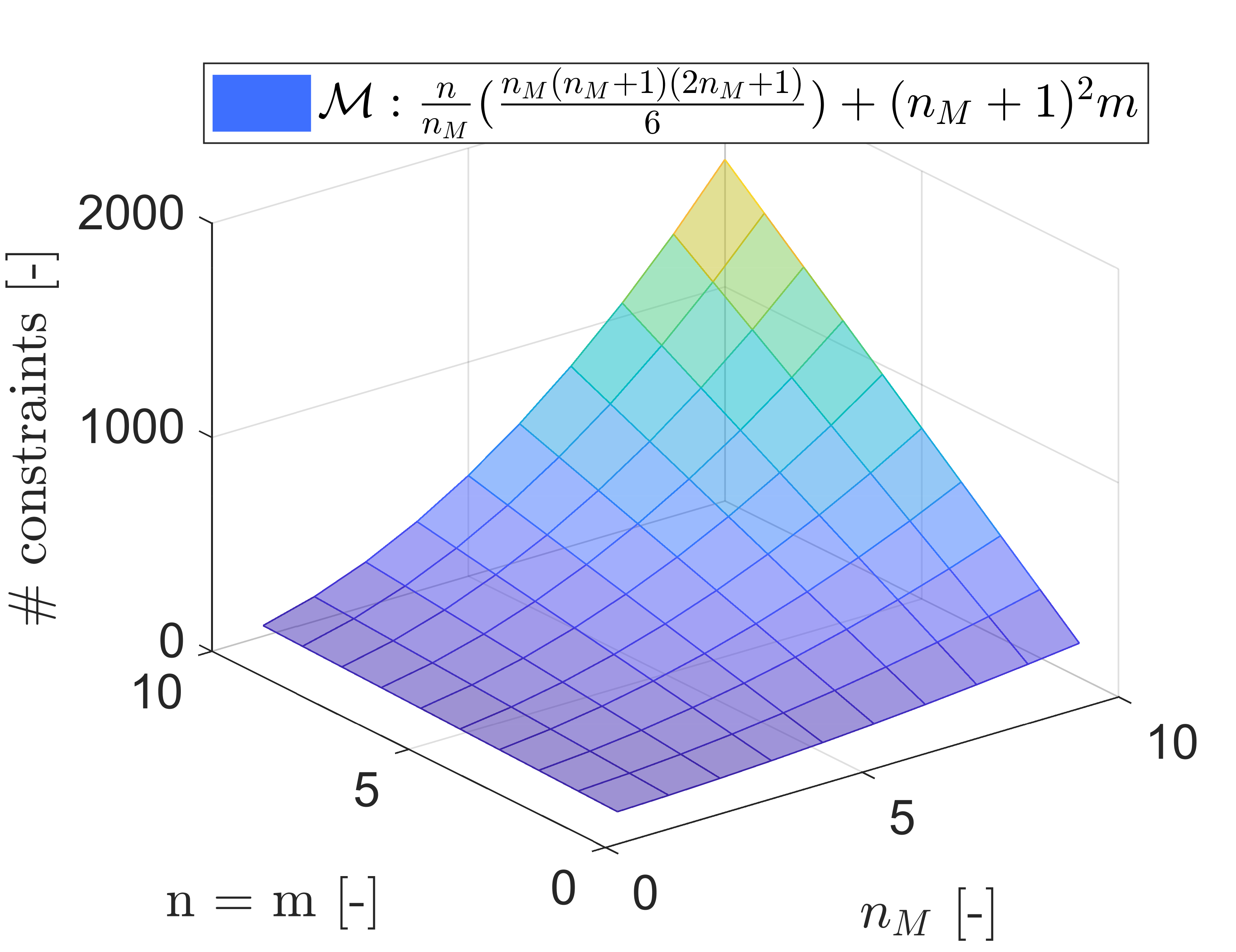} \includegraphics[width=0.32\linewidth,trim=2 2 2 4,clip]{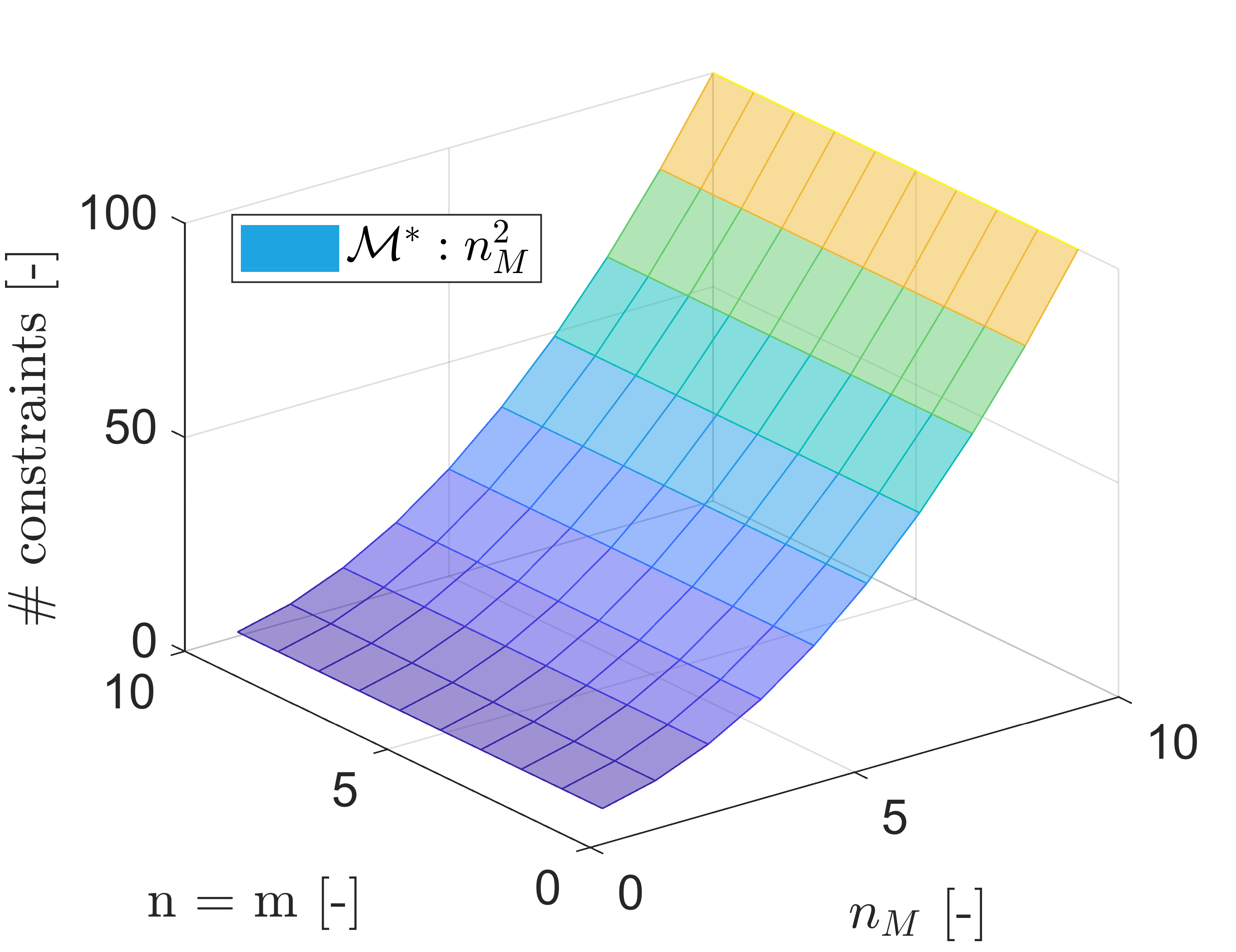}} \\ 
    \subfloat[ Passive safety between active and permanently passive agents, continuous control \label{Fig_constr_num_3}]{ \includegraphics[width=0.32\linewidth,trim=2 2 2 4,clip]{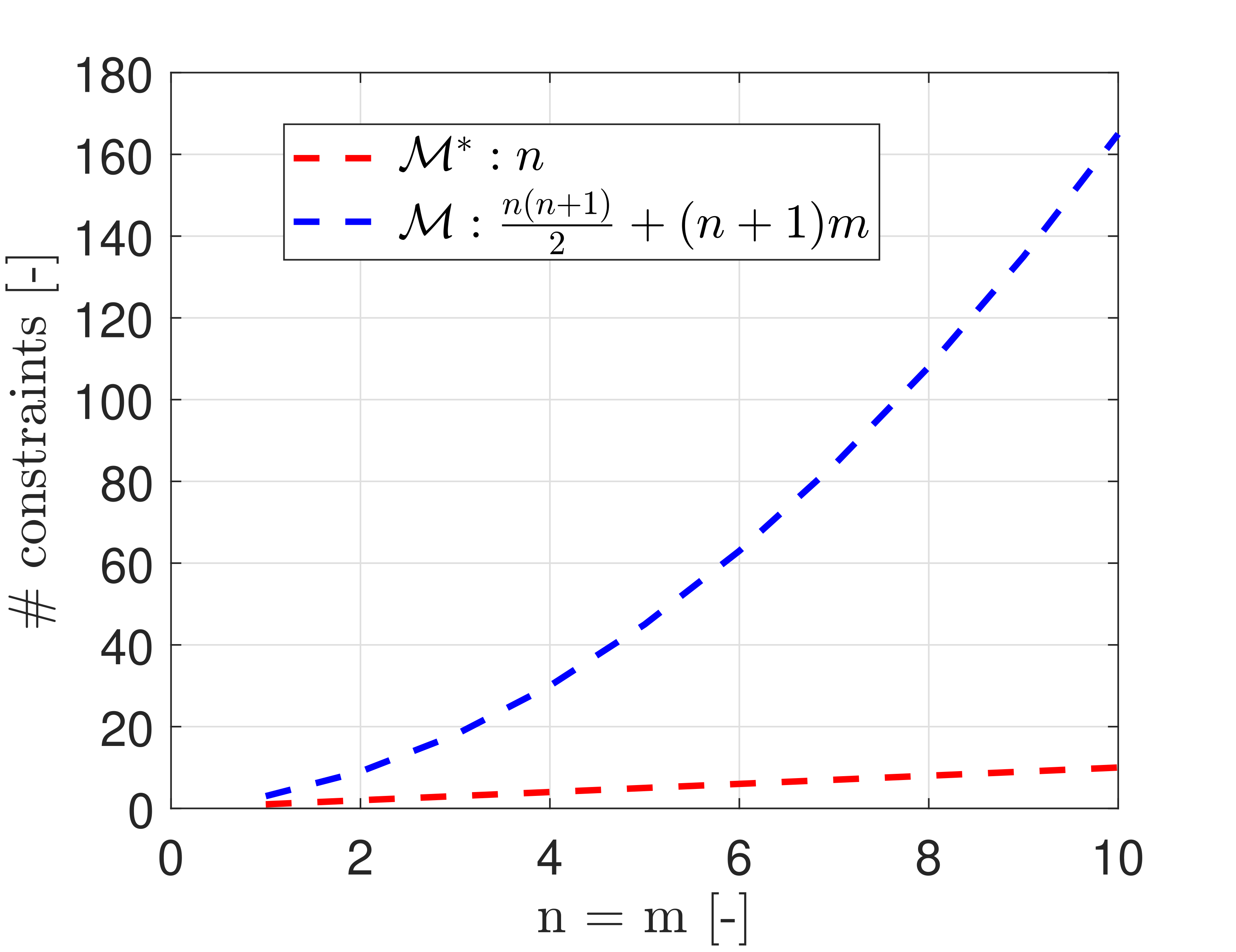}}
    \subfloat[ Passive safety between active and permanently passive agents, impulsive control \label{Fig_constr_num_4}]{ \includegraphics[width=0.32\linewidth,trim=2 2 2 4,clip]{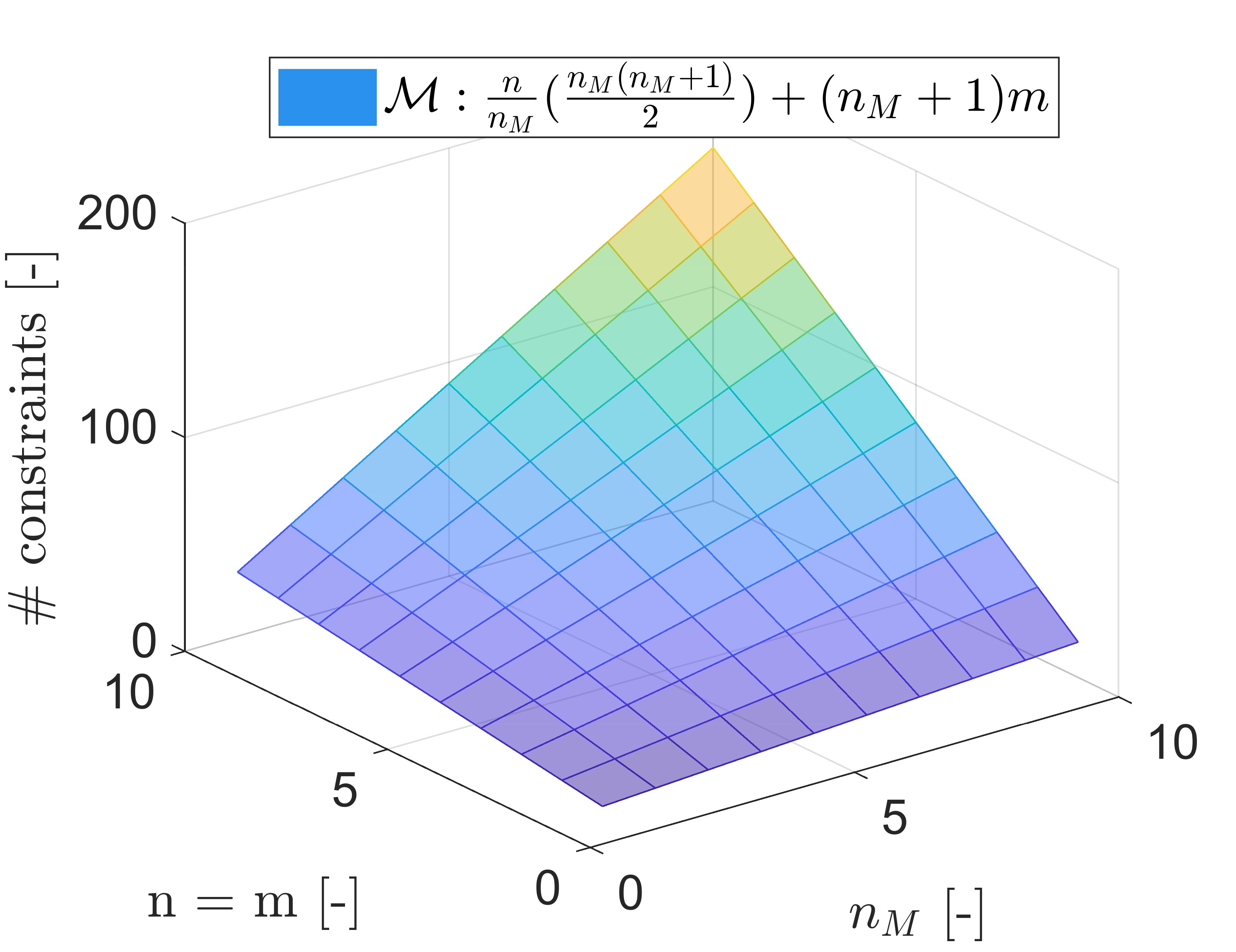}
    \includegraphics[width=0.32\linewidth,trim=2 2 2 4,clip]{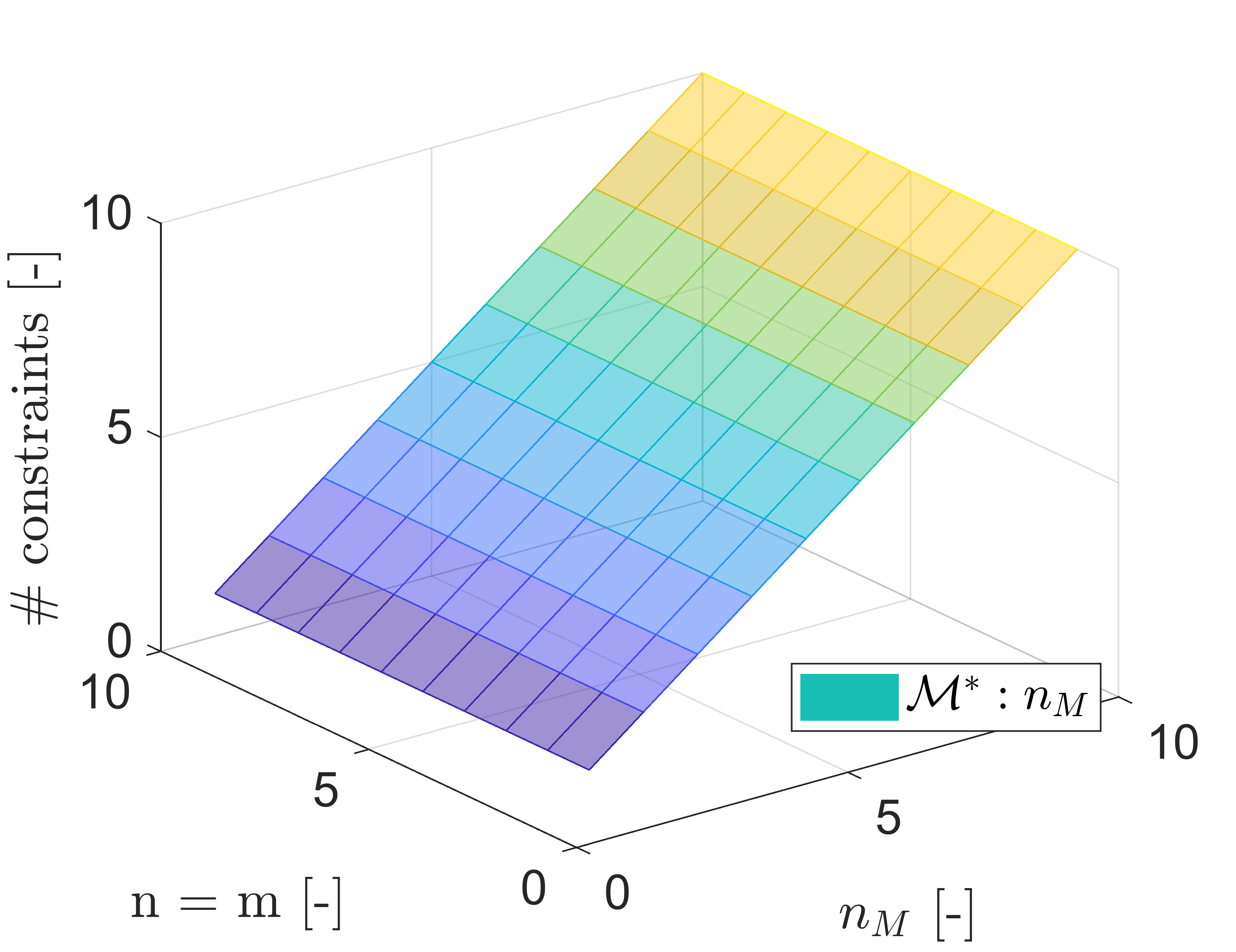}}
     \caption{Number of constraints as a function of the number of discrete time samples}\label{Fig_constr_num}
\end{figure}

By introducing  the assumption in Eq. \ref{ODE_ASS}, this high number of constraints is reduced. In particular, the procedure presented in Section \ref{Sec_II} is applied and, under the assumption of integrability of the dynamics  after a contingency, the $\mathcal{M}^*$ function is defined as
\begin{equation}\label{PS_M_star}
\resizebox{0.92\linewidth}{!}{$
\mathcal{M}^*\big(\boldsymbol{c}_{k}(t_i), \boldsymbol{c}_{j}(t_l)\big) =  \epsilon - \mathcal{S}^*_{kjil}\big(\boldsymbol{c}_{k}(t_i), \boldsymbol{c}_{j}(t_l)\big) = \epsilon - \left|\left|\boldsymbol{S} \left( \boldsymbol{\psi}\big(t^*_{kjil}\big(\boldsymbol{c}_{k}(t_i), \boldsymbol{c}_{j}(t_l)\big), \boldsymbol{c}_{k}(t_i)\big) - \boldsymbol{\psi}\big(t^*_{kjil}\big(\boldsymbol{c}_{k}(t_i), \boldsymbol{c}_{j}(t_l)\big), \boldsymbol{c}_{j}(t_l)\big) \right) \right|\right|_2 $}
\end{equation}
$\forall t_i \in [t_0, t_f]$, $\forall t_l \in [t_0, t_f]$, $k = 1,\dots N-1$, $j = k+1,\dots N$. In Eq. \ref{PS_M_star}, $\mathcal{S}^*_{kjil}$ is the minimum separation (in terms of $L_2$-norm distance) reached by agents $k^{th}$ and $j^{th}$ on the uncontrolled trajectories generated at instants $t_i$ and $t_l$ respectively. Using Eq. \ref{Eq_t_star1}, the instant $t^*_{kjil}$ of minimum separation on the uncontrolled trajectories is computed as
\begin{equation}\label{Eq_tstar}
t^*_{kjil}(\boldsymbol{c}_{k}(t_i), \boldsymbol{c}_{j}(t_l)) = \arg \min_{t \in [\max{\{t_i,t_l\}}, t_f + T]} \mathcal{S}_{kjil}(t) = \arg \min_{t \in [\max{\{t_i,t_l\}}, t_f + T]} \left|\left| \boldsymbol{S} \left( \boldsymbol{\psi}(t, \boldsymbol{c}_{k}(t_i)) - \boldsymbol{\psi}(t, \boldsymbol{c}_{j}(t_l)) \right) \right|\right|_2
\end{equation}
As derived in Eq. \ref{Eq_t_star1} and \ref{Eq_constr1}, $t^*_{kjil}$ is a function of the  integration constants (IC) states at contingency instants. This makes $\mathcal{S}^*_{kjil}$ a function of these IC states, and independent of the time along the uncontrolled trajectories.  This implies that safe separation all-along the uncontrolled trajectories is guaranteed by enforcing as many constraints as the number of combinations of contingency instants, i.e, once per each $kjil$-combination. As presented in Table \ref{tab_PS} on the bottom and proved in Appendix \ref{App:PS}, this reduces the number of constraints required to enforce PS by one polynomial degree of the number of discrete time samples. Specifically, from cubic to quadratic for each  pair of active agents, and from quadratic to linear for each  pair of active-permanently passive agents. Moreover, using the $\mathcal{M}^*$ formulation, the number of constraints required to enforce PS  becomes independent of the length of the PS horizon $T$ itself.
Figure \ref{Fig_constr_num} presents numerical examples of the trends in Table \ref{tab_PS}. To isolate the benefit of the $\mathcal{M}^*$ formulation, the number of agents is set to $N_a = 2$ in Figures \ref{Fig_constr_num_1}-\ref{Fig_constr_num_2}, and to $N_a = N_p = 1$ in Figures \ref{Fig_constr_num_3}-\ref{Fig_constr_num_4}, moreover for simplicity, $m$ is set equal to $n$.

\begin{figure}[ht!]
    \centering
    \includegraphics[width=0.62\linewidth]{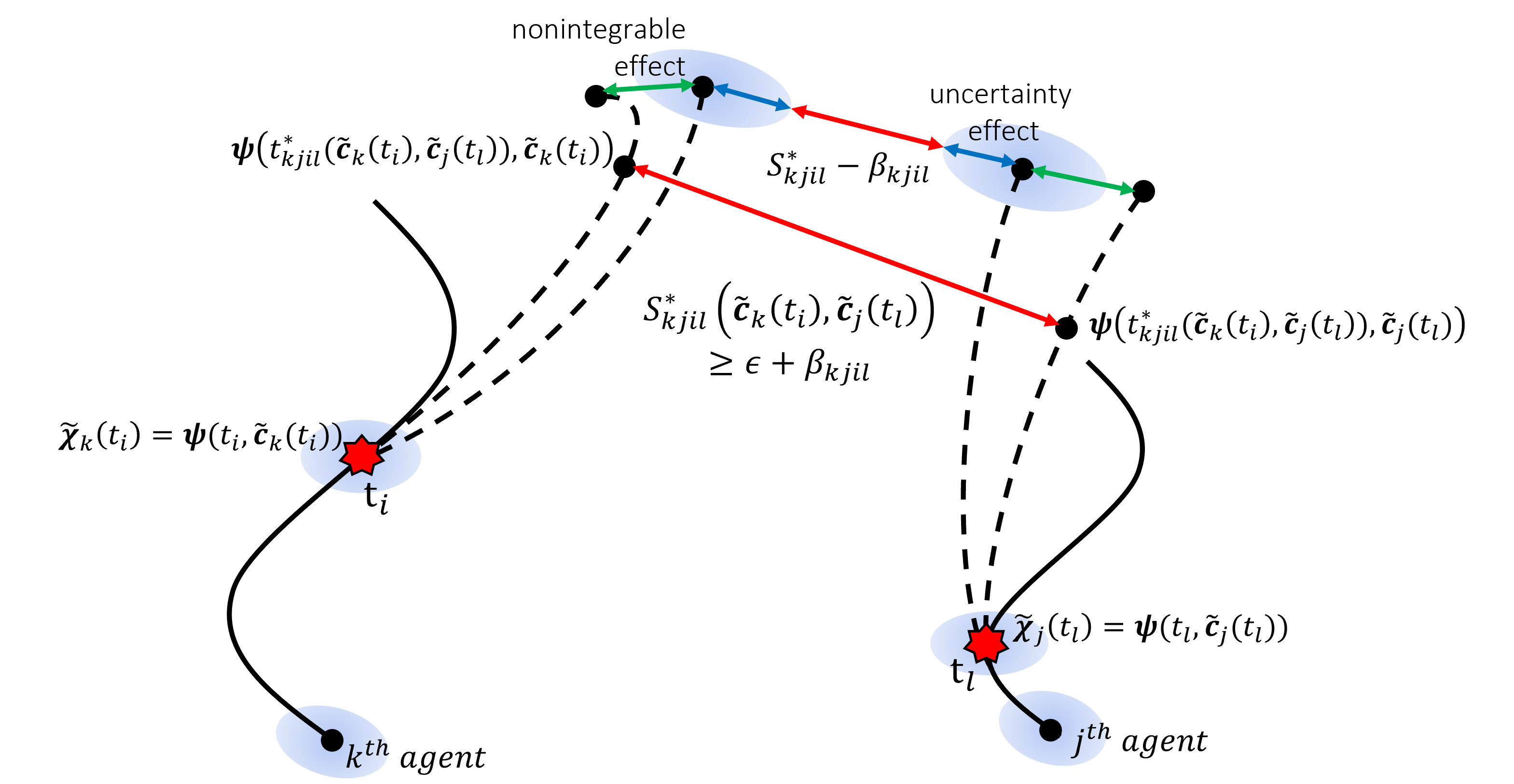}
     \caption{Compensation of non-integrable dynamics and uncertainty within the $\mathcal{M}^*$ formulation of  passive safety} \label{Fig_OSC_PS}
\end{figure}

The $\mathcal{M}^*$ formulation assumes integrability of the dynamics  after a contingency, therefore, as presented in Section \ref{Sec_II}, it requires the compensation of the effects of non-integrable dynamics and uncertainties within a bound defined as in Eq. \ref{Eq_dist_eff}. For multi-agent PS, assume these effects can be superimposed as
\begin{equation}\label{Eq_bound}
\resizebox{0.92\linewidth}{!}{$
\begin{aligned}
&\left| \frac{\partial \mathcal{M}^{*}(\boldsymbol{c}_k, \boldsymbol{c}_j)}{\partial \boldsymbol{c}_k}\bigg\lvert_{\tilde{\boldsymbol{c}}_k(t_i)} \int_{t_i}^{t_f + T} \tilde{\boldsymbol{d}}_{c}(\tau, \tilde{\boldsymbol{c}}_k(\tau)) d\tau \right| + \left| \frac{\partial \mathcal{M}^{*}(\boldsymbol{c}_k, \boldsymbol{c}_j)}{\partial \boldsymbol{c}_j}\bigg\lvert_{\tilde{\boldsymbol{c}}_j(t_l)} \int_{t_l}^{t_f + T} \tilde{\boldsymbol{d}}_{c}(\tau, \tilde{\boldsymbol{c}}_j(\tau)) d\tau \right| \\
& + \left| \frac{\partial \mathcal{M}^{*}(\boldsymbol{c}_k, \boldsymbol{c}_j)}{\partial \boldsymbol{c}_k}\bigg\lvert_{\tilde{\boldsymbol{c}}_k(t_i)} q \sqrt{ \left| Diag\left( \int_{t_i}^{t_f + T} \left( \tilde{\boldsymbol{D}}_c(\tau, \tilde{\boldsymbol{c}}_k(\tau)) \boldsymbol{C}_k(\tau) + \boldsymbol{C}_k(\tau) \tilde{\boldsymbol{D}}_c^T(\tau, \tilde{\boldsymbol{c}}_k(\tau)) + \boldsymbol{Q}_k(\tau) \right) d\tau + \boldsymbol{C}_k(t_i) \right) \right|} \right| \\
& + \left| \frac{\partial \mathcal{M}^{*}(\boldsymbol{c}_k, \boldsymbol{c}_j)}{\partial \boldsymbol{c}_j}\bigg\lvert_{\tilde{\boldsymbol{c}}_j(t_l)} q \sqrt{ \left| Diag\left( \int_{t_l}^{t_f + T} \left( \tilde{\boldsymbol{D}}_c(\tau, \tilde{\boldsymbol{c}}_j(\tau)) \boldsymbol{C}_j(\tau) + \boldsymbol{C}_j(\tau) \tilde{\boldsymbol{D}}_c^T(\tau, \tilde{\boldsymbol{c}}_j(\tau)) + \boldsymbol{Q}_j(\tau) \right) d\tau + \boldsymbol{C}_j(t_l) \right) \right|} \right| \leq \beta_{kjil}(t_i, t_l, T, q)
\end{aligned}$}
\end{equation}
where $\beta_{kjil} \in \mathbb{R}_{\geq 0}$. Note that Eq. \ref{Eq_bound} assumes the capability of retrieving an estimate of both active and permanently passive agents' states through navigation and prediction. Using Eq. \ref{Eq_bound} and Eq. \ref{PS_M_star} evaluated at the mean IC states, Eq. \ref{Eq_constr4} is specialized for multi-agent PS as
\begin{equation}\label{Eq_PS_SuffC}
\begin{aligned}
  \mathcal{S}^*_{kjil}(\tilde{\boldsymbol{c}}_{k}(t_i), \tilde{\boldsymbol{c}}_{j}(t_l)) & \geq \epsilon + \beta_{kjil}(t_i, t_l, T, q)
\end{aligned}
\end{equation}
A graphical representation is provided in Figure \ref{Fig_OSC_PS}. Assuming integrability of the dynamics  after a contingency, $\mathcal{S}_{kjil}^*$ is the minimum separation between agents $k^{th}$ and $j^{th}$ who experienced control losses at instants $t_i$ and $t_l$ respectively. The bound $\beta_{kjil}$ accounts for non-integrable dynamics and uncertainties effects on  this minimum separation.

Given these developments, the multi-agent passively-safe  optimal control problem (OCP) can be formulated as
\begin{equation}\label{OCP_IC_multi}
\resizebox{0.9\linewidth}{!}{$
\begin{aligned}
& \underset{\tilde{\boldsymbol{u}}_k(t)}{\text{minimize}} & & \sum_{k = 1}^{N_a} \int_{t_e}^{t_f} \mathcal{J}(\tilde{\boldsymbol{c}}_k(t), \boldsymbol{C}_k(t), \tilde{\boldsymbol{u}}_k(t), t) dt &\\
& \text{subject to} & & \dot{\tilde{\boldsymbol{c}}}_k(t) = \tilde{\boldsymbol{d}}_c(t, \tilde{\boldsymbol{c}}_k(t)) + \boldsymbol{B}_c(t,\tilde{\boldsymbol{c}}_k(t))\tilde{\boldsymbol{u}}_k(t) & \forall t \in [t_e, t_f], \hspace{0.2cm} k = 1,\dots N_a \\
&&& \dot{\tilde{\boldsymbol{c}}}_k(t) = \tilde{\boldsymbol{d}}_c(t, \tilde{\boldsymbol{c}}_k(t)) & \forall t \in [t_e, t_f], \hspace{0.2cm} k = 1,\dots N_p\\
&&& \dot{\boldsymbol{C}}_k(t) = \tilde{\boldsymbol{D}}_c(t, \tilde{\boldsymbol{c}}_k(t), \tilde{\boldsymbol{u}}_k(t)) \boldsymbol{C}_k(t) + \boldsymbol{C}_k(t) \tilde{\boldsymbol{D}}_c^T(t,\tilde{\boldsymbol{c}}_k(t), \tilde{\boldsymbol{u}}_k(t)) + & \forall t \in [t_e, t_f], \hspace{0.2cm} k = 1,\dots N_a \\
&&& \hspace{1.2cm} + \boldsymbol{B}_c(t,\tilde{\boldsymbol{c}}_k(t)) \boldsymbol{U}_k(t) \boldsymbol{B}_c^T(t,\tilde{\boldsymbol{c}}_k(t)) + \boldsymbol{Q}_k(t) + \boldsymbol{E}_k(t) & \\
&&& \dot{\boldsymbol{C}}_k(t) = \tilde{\boldsymbol{D}}_c(t, \tilde{\boldsymbol{c}}_k(t)) \boldsymbol{C}_k(t) + \boldsymbol{C}_k(t) \tilde{\boldsymbol{D}}_c^T(t,\tilde{\boldsymbol{c}}_k(t)) + \boldsymbol{Q}_k(t) + \boldsymbol{E}_k(t) & \forall t \in [t_e, t_f], \hspace{0.2cm} k = 1,\dots N_p \\
&&& \mathcal{S}^*_{kjil}(\tilde{\boldsymbol{c}}_{k}(t_i), \tilde{\boldsymbol{c}}_{j}(t_l)) \geq \epsilon + \beta_{kjil} (t_i, t_l, T, q) & \forall t_i \in [t_e, t_f], \hspace{0.1cm} \forall t_l \in [t_e, t_f], \\
&&&& k = 1,\dots N-1, \hspace{0.1cm} j = k+1,\dots N \\
&&& \tilde{\boldsymbol{u}}_k(t) \in \boldsymbol{\mathcal{U}}_k & \forall t \in [t_e, t_f], \hspace{0.2cm} k = 1,\dots N_a \\
&&&  \tilde{\boldsymbol{c}}_k(t_f) \in \boldsymbol{\mathcal{C}}_{k,f} & k = 1,\dots N_a \\
&&&  \mathcal{N}(\tilde{\boldsymbol{c}}_k(t_e), \boldsymbol{C}_k(t_e)) \equiv \textrm{$k^{th}$ agent's state estimate} & k = 1,\dots N\\
\end{aligned}$}
\end{equation}
 The OCP in Eq. \ref{OCP_IC_multi} enforces PS, as in  Definition 1, over $[t_e, t_f]$, $\forall t_e \in [t_0, t_f]$, at the cost of checking the reduced number of constraints presented in Table \ref{tab_PS} on the bottom. From a practical implementation perspective, to achieve the envisioned reformulation of the PS constraint, the minimum separation $\mathcal{S}^*_{kjil}$ has to be computed as a function of the IC states at contingency instants. This requires the preliminary solution of the optimization problem in Eq. \ref{Eq_tstar}. In the following section, an algorithmic solution framework is presented. Finally, note that, if PS is enforced explicitly within the optimization problem, the trend of the number of required constraints is directly connected to the similar trend of the computational effort required to solve the optimization problem itself. Therefore the proposed efficient formulation is expected to provide benefits in terms of solution run-time. This is verified and analyzed in Section \ref{Sec_VI}, applied to the solution of the passively-safe reconfiguration of VISORS \cite{VISORS}.

\section{Algorithmic Solution Framework}\label{Sec_IV}
\subsection{Solution of the Centralized Problem} \label{Sec_Alg_Centr}
 Eq. \ref{OCP_IC_multi} is a non-convex centralized optimization problem. To solve it in a centralized fashion, various optimization algorithms can be used. Here a direct approach based on  sequential convex programming (SCP) is followed \cite{scp_2000, liu_scp_2014, malyuta_scp_2022}. Specifically, this Section does not aim to derive a novel SCP algorithm, but to show how the proposed constraint reformulation can be used within a sequential optimization routine, whereby the desired fault-tolerant motion safety guarantees can be achieved. To apply SCP, non-convex constraints and cost function have to be linearized with respect to the previous solution in the series, defined as "SCP reference", and noted with $\bar{(.)}$: $\{ \bar{\tilde{\boldsymbol{c}}}_k(t), \bar{\tilde{\boldsymbol{u}}}_k(t) \}$. Moreover, the optimization problem has to be time discretized. To not replicate analyses already available in literature \cite{scp_2000, liu_scp_2014, malyuta_scp_2022}, the focus here is on the newly formulated  passive safety (PS) constraint. The turning point from an algorithmic perspective lies in the availability or not of an analytical expression of $\mathcal{S}^*_{kjil}$  as a function of the  integration constants (IC). This is presented in Algorithm \ref{ALG_C}, and analyzed in the following.
\begin{algorithm}[ht]\small

\caption{Enforces passive safety on $[t_0, t_f]$ for at least $T$, of at least $\epsilon$, at $q$-$\sigma$ confidence} \label{ALG_C}
   \DontPrintSemicolon
   
   \KwIn{Initial state estimate: $\mathcal{N}(\tilde{\boldsymbol{c}} _{k}(t_0), \boldsymbol{C}_{k}(t_0))$ , $k = 1, \dots N$. Initial guess on: $\{ \tilde{\boldsymbol{c}}_k(t_i), \tilde{\boldsymbol{u}}_k(t_i) \}$, $t_i \in [t_0, t_f]$, $k = 1, \dots N_a$.}
   
   \KwOut{ $\tilde{\boldsymbol{u}}_k(t_i)$, $t_i \in [t_0, t_f]$, $k = 1, \dots N_a$.}
   
   \KwData{$t_0, t_f, T, \epsilon, q, \boldsymbol{\mathcal{C}}_{k,f}, MAX\_ITER$, SCP parameters.}
   
   \Begin{ 
   
   \For{$t_e = t_0$ {\bf to} $t_f$}{
   
   Get $\mathcal{N}(\tilde{\boldsymbol{c}}_{k}(t_e), \boldsymbol{C}_{k}(t_e))$, $k = 1, \dots N$. \;
   
   $\{ \bar{\tilde{\boldsymbol{c}}}_k(t_i), \bar{\tilde{\boldsymbol{u}}}_k(t_i) \} \leftarrow \{ \tilde{\boldsymbol{c}}_k(t_i), \tilde{\boldsymbol{u}}_k(t_i) \}$, $t_i \in [t_e, t_f]$, $k = 1, \dots N_a$ \;
   
   \If{Closed-form expression of $\mathcal{S}^*_{kjil}$ is available}{ 
   
   Standard SCP, noting that $\bar{\beta}_{kjil}$ in Eq. \ref{PS_SCP_constr} can be precomputed at the reference.
   
   }
   \Else{
   
   Precompute $\{ \bar{t}^{*}_{kjil}, \bar{\beta}_{kjil} \}$ at the reference, $\forall kjil$.
   
   $PS \leftarrow 0$\;
   
   \While{$\mathrm{SCP}$ $\mathrm{not}$ $\mathrm{converged}$ {\bf or} $PS = 0$}{
   
   $\overline{PS} \leftarrow 0$, $ITER \leftarrow 0$\;
   
   \While{$\overline{PS} = 0$ {\bf and} $ITER < MAX\_ITER$}{
   
   Solve  Eq. \ref{OCP_IC_multi}, enforcing Eq. \ref{PS_SCP_constr2}. Get $\{ \tilde{\boldsymbol{c}}_k(t_i), \tilde{\boldsymbol{u}}_k(t_i) \}$, $t_i \in [t_e, t_f]$, $k = 1, \dots N_a$. 
   
   Recompute $\{ t^{*}_{kjil}, \beta_{kjil} \}$ at the new solution, $\forall kjil$.
   
   \If {$\mathcal{S}_{kjil}(t^*_{kjil}, \tilde{\boldsymbol{c}}_{k}(t_i), \tilde{\boldsymbol{c}}_{j}(t_l)) \geq \epsilon + \beta_{kjil}$, $\forall kjil$}{
   
   $\overline{PS} \leftarrow 1$\;
   
   }
   
   $\{ \bar{t}^{*}_{kjil}, \bar{\beta}_{kjil} \} \leftarrow \{ t^{*}_{kjil}, \beta_{kjil} \}$, $\forall kjil$.\;
   
   $ITER \leftarrow ITER + 1$\;
   
   }
   
   \If {$\overline{PS} = 0$}{
   
   Select among the $MAX\_ITER$ solutions computed the one that got closer to satisfy $\mathcal{S}_{kjil}(t^*_{kjil}, \tilde{\boldsymbol{c}}_{k}(t_i), \tilde{\boldsymbol{c}}_{j}(t_l)) \geq \epsilon + \beta_{kjil}$, $\forall kjil$. \;
   
   }
   
   $\{ \bar{\tilde{\boldsymbol{c}}}_k(t_i), \bar{\tilde{\boldsymbol{u}}}_k(t_i) \} \leftarrow \{ \tilde{\boldsymbol{c}}_k(t_i), \tilde{\boldsymbol{u}}_k(t_i) \}$, $t_i \in [t_e, t_f]$, $k = 1, \dots N_a$.\;
   
   Recompute $\{ \bar{t}^{*}_{kjil}, \bar{\beta}_{kjil} \}$ at the new reference, $\forall kjil$.
   
   $PS \leftarrow \overline{PS}$\;
   
   }
   
   }
   
   Get $\tilde{\boldsymbol{u}}_k(t_e)$, $k = 1, \dots N_a$.\;
   
   }
   
   \Return $\tilde{\boldsymbol{u}}_k(t_i)$, $t_i \in [t_0, t_f]$, $k = 1, \dots N_a$. \;
   
   }
   
\end{algorithm}

 If an analytical expression of $\mathcal{S}^*_{kjil}$ is available and is, in general, a non-convex function of the IC state, it can be sequentially linearized around the SCP reference as
\begin{equation}\label{PS_SCP_constr}
\resizebox{0.93\textwidth}{!}{$
\begin{aligned}
    & \frac{\partial \mathcal{S}^*_{kjil}(\tilde{\boldsymbol{c}}_{k}, \tilde{\boldsymbol{c}}_{j})}{\partial \tilde{\boldsymbol{c}}_{k}} \bigg\lvert_{\bar{\tilde{\boldsymbol{c}}}_{k}(t_i), \bar{\tilde{\boldsymbol{c}}}_{j}(t_l)} \Big( \tilde{\boldsymbol{c}}_k(t_i) - \bar{\tilde{\boldsymbol{c}}}_k(t_i) \Big) + \frac{\partial \mathcal{S}^*_{kjil}(\tilde{\boldsymbol{c}}_{k}, \tilde{\boldsymbol{c}}_{j})}{\partial \tilde{\boldsymbol{c}}_{j}}\bigg\lvert_{\bar{\tilde{\boldsymbol{c}}}_{k}(t_i), \bar{\tilde{\boldsymbol{c}}}_{j}(t_l)} \Big( \tilde{\boldsymbol{c}}_j(t_l) - \bar{\tilde{\boldsymbol{c}}}_j(t_l) \Big) + \mathcal{S}^*_{kjil}(\bar{\tilde{\boldsymbol{c}}}_{k}(t_i), \bar{\tilde{\boldsymbol{c}}}_{j}(t_l)) \geq \epsilon + \bar{\beta}_{kjil}
\end{aligned}$}
\end{equation}
where $\bar{\beta}_{kjil}$ (defined in Eq. \ref{Eq_bound}) is, by choice, directly evaluated at the reference. Since, from a PS enforcement perspective, the IC state covariances are used just to quantify this bound, this  permits: 1) pre-propagating them around the reference mean IC states and control inputs, 2) quantifying $\bar{\beta}_{kjil}$ evaluating the terms on the left-hand side of Eq. \ref{Eq_bound} at the reference, 3) enforcing Eq. \ref{PS_SCP_constr} within the optimization problem, without having to enforce explicitly the convariance dynamics constraint. This is an approximate but simpler solution than linearizing the terms on the left-hand side of Eq. \ref{Eq_bound}. The accuracy of this approximation can be regulated through additional trust region constraints \cite{liu_scp_2014, malyuta_scp_2022}.  Note that, by definition, the integrals in Eq. \ref{Eq_bound}, cannot be solvable explicitly and in closed-form. Therefore, they are pre-computed at the SCP reference using, for example, numerical integration, or even semi-analytical techniques if models of the non-integrable portion of the dynamics are available, as for the distributed space system application \cite{gim_state_2003, koenig_2017_new, guffanti_long-term_2017, guffanti_linear_2019}. The resulting bound $\beta$ is fed as a parameter in the optimization problem solved at each SCP iteration. Moreover, note that the IC states of permanently passive agents at $t_0$ can be treated as input parameters in the optimization, and not as optimization variables. This simplifies the left-hand side of Eq. \ref{PS_SCP_constr} for active-permanently passive agents  pairs. In particular, the linearization with respect to the permanently passive agents' IC states is not performed.

\begin{figure}[ht!]
    \centering
    \subfloat[Enforcement of PS within each SCP interation 
 \label{Fig_tstar_alg_1}]{\includegraphics[width=0.58\linewidth]{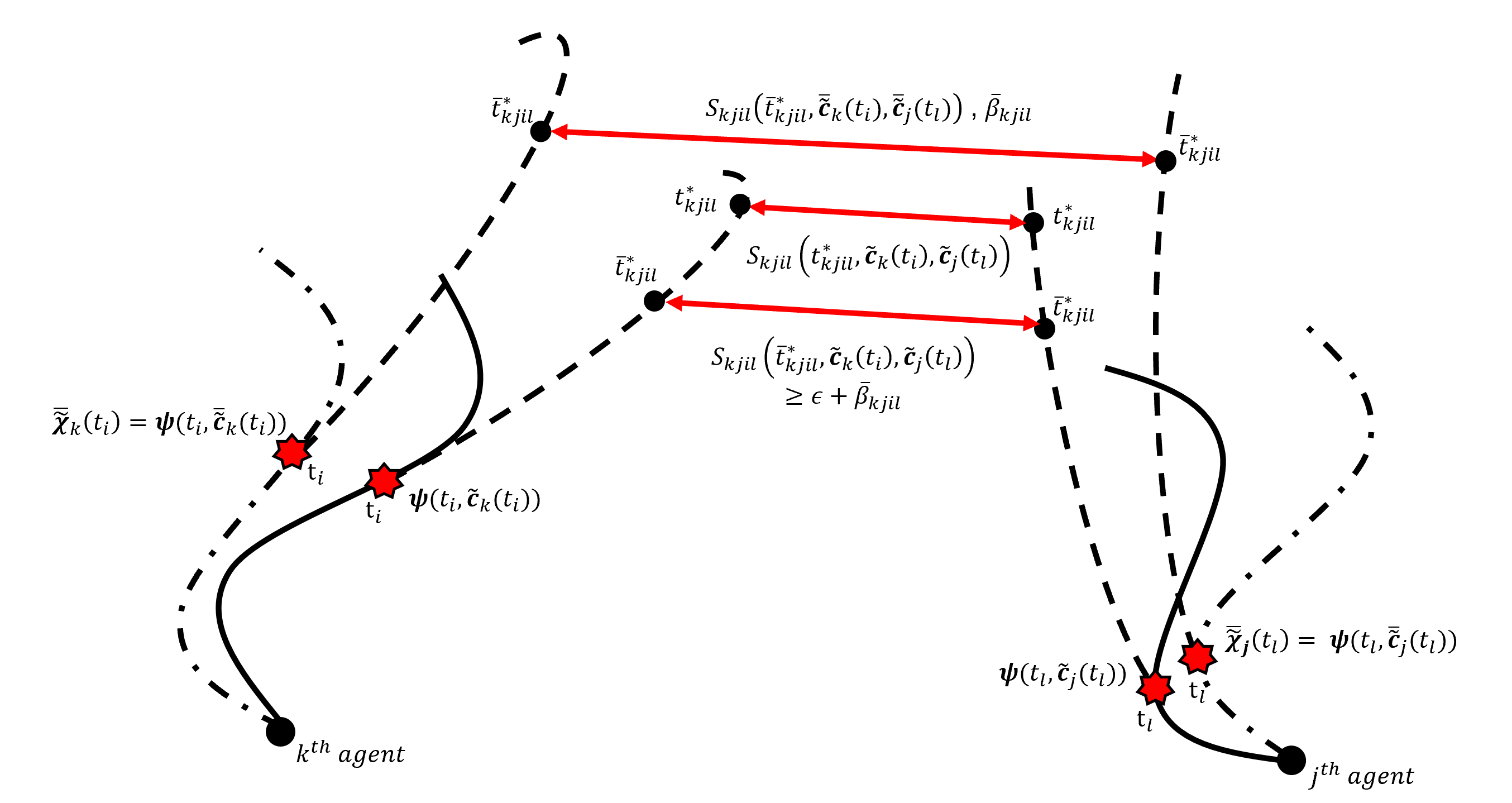}}
    \subfloat[Trend of convergence moving through the SCP iterations 
 \label{Fig_tstar_alg_2}]{\includegraphics[width=0.41\linewidth,trim=2 2 2 2,clip]{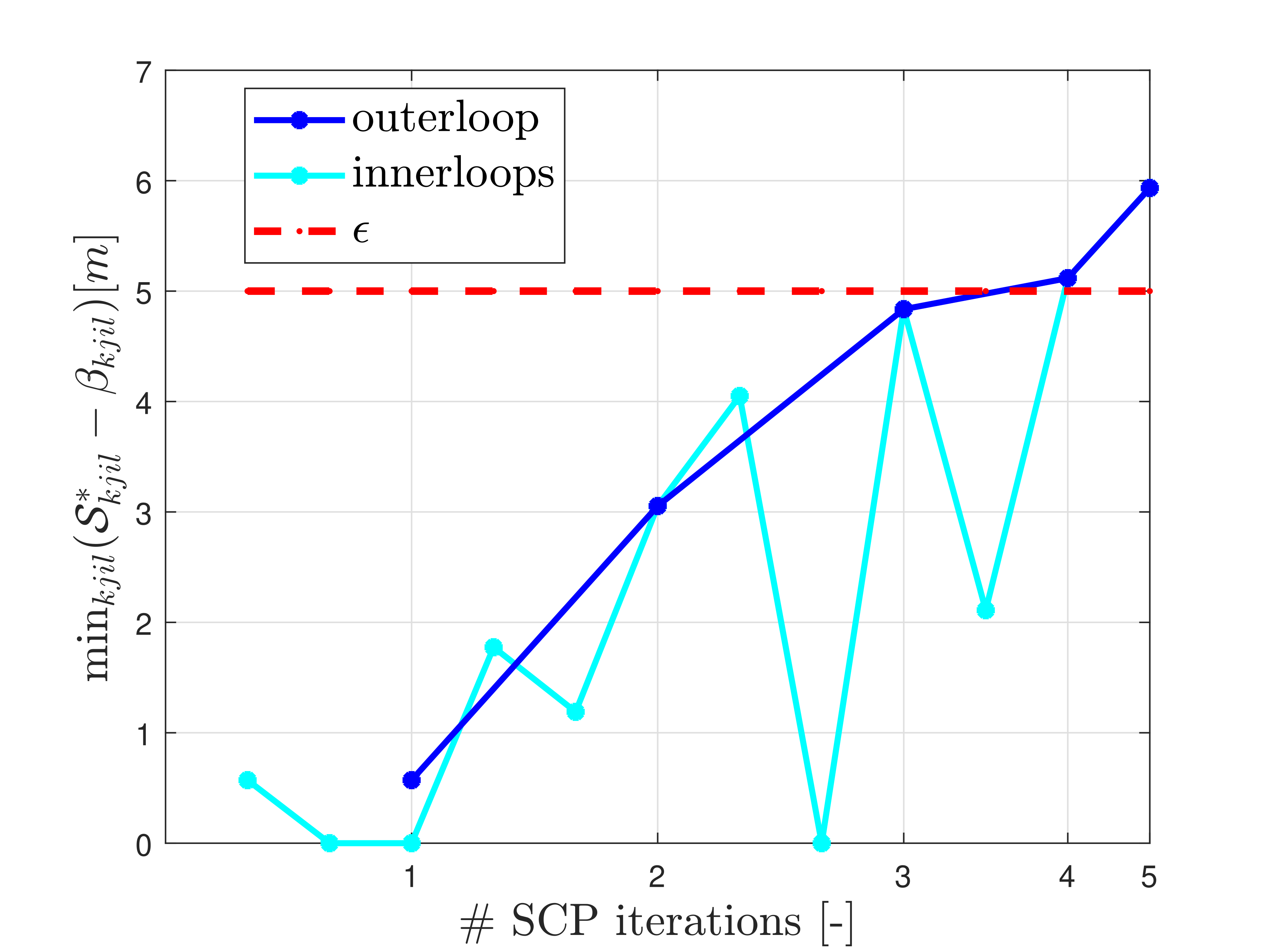}}
     \caption{Enforcement of  passive safety (PS) within the  sequential convex programming (SCP) pipeline}\label{Fig_tstar_alg}
\end{figure}

If an analytical expression for $\mathcal{S}^*_{kjil}$ is not available, due to the fact that Eq. \ref{Eq_tstar} is too complicated to be solved in closed-form, then an algorithmic development is needed. This is presented in Algorithm \ref{ALG_C} (lines 9-24) and in the following. In particular, at each iteration of the SCP loop, $\bar{t}^*_{kjil}$ is precomputed at the reference through a simple search algorithm on the discretized interval $[\max{\{t_i,t_l\}}, t_f + T]$. Using Eq. \ref{Eq_tstar}, this search algorithm selects the discrete instant/index that minimizes $\mathcal{S}_{kjil}(t, \bar{\tilde{\boldsymbol{c}}}_{k}(t_i), \bar{\tilde{\boldsymbol{c}}}_{j}(t_l))$, available as a closed-form expression of the reference mean IC states under assumption of integrability of the dynamics  after a contingency. Using the precomputed $\bar{t}^*_{kjil}$, the PS constraint is enforced as
\begin{equation}\label{PS_SCP_constr2}
\resizebox{0.93\textwidth}{!}{$
    \frac{\partial \mathcal{S}_{kjil}(\bar{t}^*_{kjil},\tilde{\boldsymbol{c}}_{k}, \tilde{\boldsymbol{c}}_{j})}{\partial \tilde{\boldsymbol{c}}_{k}}\bigg\lvert_{\bar{\tilde{\boldsymbol{c}}}_{k}(t_i), \bar{\tilde{\boldsymbol{c}}}_{j}(t_l)} \Big( \tilde{\boldsymbol{c}}_k(t_i) - \bar{\tilde{\boldsymbol{c}}}_k(t_i) \Big) + \frac{\partial \mathcal{S}_{kjil}(\bar{t}^*_{kjil},\tilde{\boldsymbol{c}}_{k}, \tilde{\boldsymbol{c}}_{j})}{\partial \tilde{\boldsymbol{c}}_{j}} \bigg\lvert_{\bar{\tilde{\boldsymbol{c}}}_{k}(t_i), \bar{\tilde{\boldsymbol{c}}}_{j}(t_l)} \Big( \tilde{\boldsymbol{c}}_j(t_l) - \bar{\tilde{\boldsymbol{c}}}_j(t_l) \Big) + \mathcal{S}_{kjil}(\bar{t}^*_{kjil},\bar{\tilde{\boldsymbol{c}}}_{k}(t_i), \bar{\tilde{\boldsymbol{c}}}_{j}(t_l)) \geq \epsilon + \bar{\beta}_{kjil}$}
\end{equation}
With arguments similar to the one used for the standard collision avoidance in \cite{morgan_2014}, it is possible to prove that Eq. \ref{PS_SCP_constr2} represents a sufficient condition for the satisfaction of $\mathcal{S}_{kjil}(\bar{t}^*_{kjil},\tilde{\boldsymbol{c}}_{k}(t_i), \tilde{\boldsymbol{c}}_{j}(t_l)) \geq \epsilon + \bar{\beta}_{kjil}$. The formal proof is  presented by the authors in \cite{guffanti_2021_multi, guffanti_phd_2022}. Nevertheless, Eq. \ref{PS_SCP_constr2} does not constitutes a sufficient condition of having the actual minimum separation at $t^*_{kjil}$ being greater than $\epsilon + \bar{\beta}_{kjil}$. This is represented graphically in Figure \ref{Fig_tstar_alg_1}. At the SCP reference, minimum separation on the uncontrolled trajectories is reached at instant $\bar{t}^{*}_{kjil}$. At the current solution, at instant $t^{*}_{kjil}$. The separation reached at instant $\bar{t}^{*}_{kjil}$ by the current solution may be larger than the minimum one since $t^{*}_{kjil} \neq \bar{t}^{*}_{kjil}$.
Therefore, $\bar{t}^*_{kjil}$ has to be sequentially updated within the SCP pipeline to guarantee that, at convergence, the actual worst case minimum separation is constrained. In particular, 1) the optimization problem is solved by enforcing Eq. \ref{PS_SCP_constr2} which uses $\bar{t}^*_{kjil}$ and $\bar{\beta}_{kjil}$ evaluated at the reference, then 2) $t^*_{kjil}$ and $\beta_{kjil}$ are computed at the new solution, then 3) the actual minimum separation $\mathcal{S}_{kjil}(t^*_{kjil},\tilde{\boldsymbol{c}}_{k}(t_i), \tilde{\boldsymbol{c}}_{j}(t_l))$ is checked a-posteriori to be $\geq$ than $\epsilon + \beta_{kjil}$. If this condition is true, the algorithm proceeds to the next SCP-loop iteration. If it is not true, the optimization problem is solved again keeping the same linearization reference but using the newly computed $t^{*}_{kjil}$ and $\beta_{kjil}$.  This inner-loop breaks if PS is satisfied a-posteriori, or if a maximum number of iterations ($MAX\_ITER \in \mathbb{N}_{>0}$) is reached. If the maximum number of iterations is reached, the solution in the inner-loop  closest to satisfying PS is taken as reference and the algorithm proceeds to the subsequent SCP iteration. Introducing a maximum number of iterations in the inner-loop prevents from getting stuck into intermediate solutions, in which the PS constraint does not get satisfied in few iterations. Figure \ref{Fig_tstar_alg_2} presents a sample trend of  PS enforcement moving through the inner-loops of Algorithm \ref{ALG_C} (lines 12-18, light-blue line, with $MAX\_ITER = 3$, and $\epsilon = 5$m) and the SCP outer-loop (blue line).  The algorithm is considered converged to a feasible solution, only if the final solution out of the SCP loop is passively-safe. Otherwise, the problem is deemed unfeasible.

 The feasibility of finding a passively-safe solution trajectory depends on various factors, including: the boundary conditions, the level of uncertainty present in the system (which influences the size of the PS bounds), and optimization parameters such as the initial warm-starting trajectory, and the defined time discretization. Certain optimization parameters can be properly selected and tuned to increase the chances of finding a feasible solution. Other factors such as the boundary conditions or the level of uncertainty are defined given the system and control scenario. Therefore, given a certain system, control scenario and optimization set-up, the possibility of not finding a feasible passively-safe solution is present. In risk sensitive applications, if a feasible solution is repeatedly not found, it is necessary to have a logic that falls back on an alternative action. For distributed space systems, this alternative action can be, for example, an escape maneuver performed by those agents that have not completely lost control.

Both if $\mathcal{S}^*_{kjil}$ is available in closed-form and not, Eq. \ref{PS_SCP_constr} and \ref{PS_SCP_constr2} enforce the reduced number of constraints presented in Table \ref{tab_PS} on the bottom. This possibly cuts the computational expense, and increases the chances of converging in useful time to a feasible solution using off-the-shelf solvers \cite{boyd_book, ecos_13}.
As final mention, as in any SCP algorithm, to implement Algorithm \ref{ALG_C} an initial "warm-start" guess on the profile: $\{ \tilde{\boldsymbol{c}}_{k}(t_i), \tilde{\boldsymbol{u}}_k(t_i) \}$, $\forall ki$, is required as input.  This initial guess does not have to be feasible, i.e., it does not have to satisfy the PS constraint. An option is to initialize the algorithm with a profile that simply solves for the two-point-boundary-value-problem from initial to final conditions, or, better, that satisfies the standard collision avoidance constraint but is not passively-safe. Obviously, the computational effort is reduced if a sub-optimal passively-safe  solution is already available, as analyzed in Section \ref{Sec_VI}.

\subsection{Decentralization}\label{Sec_Alg_Dec}
In a multi-agent setting,  using Algorithm \ref{ALG_C} at each $t_e \in [t_0, t_f]$ requires: 1) the communication of the state estimates of each agent to the one assigned with the role of computing the centralized plan, and 2) the communication back to each agent of the newly computed plan. From an operational stand-point,  it may be beneficial to do the computation of the centralized plan at the initial instant $t_0$ and then at subsequent instants $t_e$ only when necessary. This subtends the need of complementing the centralized problem with an additional decentralized optimization problem, defined conventionally as the "tracking problem". This enables each agent to independently compute a "tracking control plan" that tracks a previously computed centralized plan, any time an updated state estimate is available. The tracking problem has to compensate for the presence of initial tracking errors, as well as for the acting dynamics, while reinforcing guarantees of constraint satisfaction. In \cite{guffanti_phd_2022}, the authors present a possible formulation of the tracking problem, as well as the logic regulating the interplay between centralized and decentralized problems solution. Here, for  conciseness, the key aspects are summarized.  Introduce the superscript $(.)^{\mathcal{G}}$ to indicate variables associated to the centralized "guidance" and the superscript $(.)^{\mathcal{T}}$ to indicate variables associated to the decentralized "tracking". Algorithm \ref{ALG_C} is used at $t_0$ to produce the centralized plan $\{ \tilde{\boldsymbol{u}}^{\mathcal{G}}_k(t), \tilde{\boldsymbol{c}}^{\mathcal{G}}_k(t), \boldsymbol{C}^{\mathcal{G}}_k(t)\}$, $\forall t \in [t_0, t_f]$, for each active agent. If feasible, the centralized plan is certified as passively-safe. At each instant $t_e$, each agent acquires a state estimate $\mathcal{N}(\tilde{\boldsymbol{c}} _{k}(t_e),\boldsymbol{C} _{k}(t_e))$. Given this estimate and the available centralized plan, a tracking error is possibly present and defined as $\tilde{\boldsymbol{c}}^{\mathcal{T}}_k(t_e) = \tilde{\boldsymbol{c}} _{k}(t_e) - \tilde{\boldsymbol{c}}^{\mathcal{G}}_k(t_e)$. This tracking error propagates through the dynamics of the system. In addition, the newly estimated covariance $\boldsymbol{C} _{k}(t_e)$ is not necessarily equal to the one predicted during the solution of the centralized problem, and used to evaluate the bounds in Eq. \ref{PS_SCP_constr} and \ref{PS_SCP_constr2}. To make the guidance plan robust to these tracking effects, the idea is to enlarge the bound used at guidance level by an additional margin. In particular, for each involved agent, the following term can be added to the left-hand side of Eq. \ref{Eq_bound}: $\left|(\partial \mathcal{M}^{*}/\partial \boldsymbol{c}_k)|_{\tilde{\boldsymbol{c}}_{k}(t_i)} \cdot \boldsymbol{\nu}_{ki}^{\mathcal{T}} \right|$.  Here, $\boldsymbol{\nu}_{ki}^{\mathcal{T}} \in \mathbb{R}^s_{\geq 0}$ is a constant total worst-case tracking error which is then reinforced at tracking level on the mean tracking state as
\begin{equation}\label{Eq_track_bound0}
\big| \tilde{\boldsymbol{c}}^{\mathcal{T}}_k(t_i) \big| + \boldsymbol{\beta}_{ki}^{\mathcal{T}} \leq \boldsymbol{\nu}_{ki}^{\mathcal{T}}
\end{equation}
$\forall t_i \in [t_0, t_f]$, for $k = 1, \dots N_a$. In Eq. \ref{Eq_track_bound0}, $\boldsymbol{\beta}_{ki}^{\mathcal{T}} \in \mathbb{R}^s_{\geq 0}$ is quantified such that
\begin{equation}\label{Eq_track_bound1}
\begin{aligned}
    \boldsymbol{\beta}_{ki}^{\mathcal{T}} \geq & \left| \int_{t_i}^{t_f + T} \tilde{\boldsymbol{d}}_{c}(\tau, \tilde{\boldsymbol{c}}_k(\tau)) d\tau - \int_{t_i}^{t_f + T} \tilde{\boldsymbol{d}}_{c}(\tau, \tilde{\boldsymbol{c}}^{\mathcal{G}}_k(\tau)) d\tau \right| \\
    & + \left| q \sqrt{\left| Diag\left( \int_{t_i}^{t_f + T} \left( \tilde{\boldsymbol{D}}_c(\tau, \tilde{\boldsymbol{c}}_k(\tau)) \boldsymbol{C}_k(\tau) + \boldsymbol{C}_k(\tau) \tilde{\boldsymbol{D}}_c^T(\tau, \tilde{\boldsymbol{c}}_k(\tau)) + \boldsymbol{Q}_k(\tau) \right) d\tau + \boldsymbol{C}_k(t_i) \right) \right|} \right. \\
    & \left. - q \sqrt{\left| Diag\left( \int_{t_i}^{t_f + T} \left( \tilde{\boldsymbol{D}}_c(\tau, \tilde{\boldsymbol{c}}^{\mathcal{G}}_k(\tau)) \boldsymbol{C}^{\mathcal{G}}_k(\tau) + \boldsymbol{C}^{\mathcal{G}}_k(\tau) \tilde{\boldsymbol{D}}_c^T(\tau, \tilde{\boldsymbol{c}}^{\mathcal{G}}_k(\tau)) + \boldsymbol{Q}^{\mathcal{G}}_{k}(\tau) \right) d\tau + \boldsymbol{C}^{\mathcal{G}}_k(t_i) \right) \right|} \right|
\end{aligned}
\end{equation}
to capture the net effect of tracking on the non-integrable and uncertainty effects quantified at guidance level. In particular, in Eq. \ref{Eq_track_bound1}, the non-superscripted state quantities refer to the ones computed at tracking level, which depend on the tracking error itself.

Finally,  a high-level logic regulating the interplay between guidance and tracking is implemented. Specifically, at every instant $t_e$ an updated state estimate is available, the system has two options: a) tracking the centralized plan or, b) recomputing the centralized plan. The high-level logic has to decide autonomously on-board between a) or b). The triggers in this decision are mainly two: 1) the feasibility of the tracking, and 2) the prediction that the currently planned control solution leads to constraints satisfaction, given the latest state estimate and the on-board knowledge of the system dynamics. Therefore, at each instant $t_e$, each agent acquires an updated state estimate and tracks. If for some agents, the tracking is not feasible the agent(s) share the information and the "planner agent" recomputes the centralized plan using the latest state estimates. If the tracking is feasible for everyone, the agent(s) share their updated state estimates and the corrected control plans to the planner agent, which  performs a sanity-check. In particular, it propagates in open-loop both mean and covariance trajectories of each agent starting at the newly acquired state estimates using the corrected control plans that has been communicated. This open-loop propagation  permits checking if the expected worst case minimum separations reached on the uncontrolled trajectories are actually greater than $\epsilon$. If any of them is not greater than $\epsilon$, PS is not guaranteed by the current control plan given the newly acquired state estimate. Therefore, the centralized plan must be recomputed. Otherwise, if all of them are greater than $\epsilon$, the agents receive the go ahead to actuate the corrected control plan. This sanity-check performed by the planner agent introduces delay and computational overhead in the control pipeline, requiring the communication back and forth of data between agents as well as the open-loop propagation needed for checking the minimum separations. Depending on the risk sensitivity of the control scenario it is possible to skip this sanity-check or perform it just at key instants of the multi-agent reconfiguration (e.g., after the actuation of a largely uncertain maneuver). If the sanity-check is not performed, the only decision trigger on recomputing the centralized plan or not remains the possible unfeasibility of the tracking.

\section{Application to Distributed Space Systems}\label{Sec_V}
The orbital motion of a $N$-spacecraft  distributed space system can be described either by $N$ orbits in inertial space around the main attractor or by one reference orbit, either virtual or occupied by a chief spacecraft, and by $N$ (or $N-1$) relative trajectories with respect to this reference. These two ways of modeling the spacecraft motion influence the structure of the governing ODE, and therefore give rise to two different ways of applying  variation of parameters and the methodology proposed in Sections \ref{Sec_II}-\ref{Sec_IV}. As described in the introduction, one approach uses the  orbital elements (OE) as  integration constants (IC) state, while the other approach uses the  relative orbital elements (ROE). From an operational stand-point, for spacecraft in close-proximity, modeling directly the relative motion is advantageous both from a dynamics modeling standpoint \cite{gim_state_2003, koenig_2017_new, guffanti_long-term_2017, guffanti_linear_2019} and navigation accuracy standpoint \cite{giralo_phd_2021}. In particular, spacecraft absolute motion navigation solutions are  usually at least one order of magnitude less accurate than the relative motion ones. Navigation accuracy plays a fundamental role in the developed framework, since it determines the level of uncertainty present, and in face of which  passive safety (PS) guarantees have to be provided (see Eq. \ref{Eq_bound}). This is analyzed in the test cases in Section \ref{Sec_VI}.

Following the relative motion approach, the reference orbit is uniquely defined by a set of OE: $\textbf{\oe} \in \mathbb{R}^6$.  Use the quasi-nonsingular OE definition \cite{vallado}: $\textbf{\oe} = \left\{ a, \lambda_\nu, e_{x}, e_{y}, i, \Omega \right \}$, with $a$ the semi-major axis, $\lambda_\nu = \nu + \omega$ the true argument of latitude, $\nu$ the true anomaly, $\omega$ the argument of perigee, $\{e_x, e_y \} = \{ e \cos(\omega), e \sin(\omega) \}$ the eccentricity vector, $e$ the eccentricity, $i$ the inclination, and $\Omega$ the right ascension of the ascending node. The $N$ (or $N-1$) relative trajectories are defined by the relative Cartesian state expressed in the radial/tangential/normal (RTN) \cite{vallado} frame centered on the reference, as $\boldsymbol{\chi}_{k} = \{ \delta \boldsymbol{r}_{k}, \delta \dot{\boldsymbol{r}}_{k} \} \in \mathbb{R}^6$, for $k= 1,\dots N$, where $\delta \boldsymbol{r}_{k} \in \mathbb{R}^3$ and $\delta \dot{\boldsymbol{r}}_{k} \in \mathbb{R}^3$ are correspondingly the relative position and velocity of the $k^{th}$ spacecraft. The dynamics of $\boldsymbol{\chi}_{k}$ is governed by a nonlinear ODE \cite{casotto_equations_2016}
\begin{equation}\label{ODE_nl}
\dot{\boldsymbol{\chi}}_{k}(t) = \boldsymbol{z}(\boldsymbol{\chi}_{k}(t)) + \boldsymbol{p}(\boldsymbol{\chi}_{k}(t)) + \boldsymbol{B} \boldsymbol{u}_k(t)
\end{equation}
where, $\boldsymbol{z} \in \mathbb{R}^6$ models the relative effect of two-body spherical gravity as well as Euler, centrifugal, and Coriolis accelerations  stemming from the non-inertial (and potentially perturbed) nature of the rotating RTN frame, $\boldsymbol{p} \in \mathbb{R}^6$ models the differential effect of orbital perturbations, and $\boldsymbol{u}_k \in \mathbb{R}^3$ is the $k^{th}$ spacecraft control acceleration expressed in RTN frame, with $\boldsymbol{B} = [\boldsymbol{0}_{3 \times 3}; \boldsymbol{I}_{3 \times 3}] \in \mathbb{R}^{6 \times 3}$ the control input matrix. Although the relative dynamics in Eq. \ref{ODE_nl} depends on the reference OE, this dependence is left implicit in the notation for compactness. In literature, analytical solutions for the homogeneous nonlinear part of Eq. \ref{ODE_nl} are not available. Nevertheless, linearizing the differential spherical gravity terms in $\boldsymbol{z}$ for small spacecraft separations, Eq. \ref{ODE_nl} is rewritten exactly as \begin{equation}\label{ODE_lin}
\dot{\boldsymbol{\chi}}_{k}(t) = \boldsymbol{A}(t) \boldsymbol{\chi}_{k}(t) + \boldsymbol{d}(\boldsymbol{\chi}_{k}(t)) + \boldsymbol{B} \boldsymbol{u}_k(t)
\end{equation}
where  $\boldsymbol{A} \in \mathbb{R}^{6 \times 6}$ is the linear(ized) integrable part, and $\boldsymbol{d} \in \mathbb{R}^6$ is the nonlinear non-integrable part, which contains both $\boldsymbol{p}$ and the truncated higher-order terms of $\boldsymbol{z}$. Since the beginning of spaceflight, solutions for the linear part of Eq. \ref{ODE_lin} have been developed by different authors. For example,  Clohessy and Wiltshire (CW) developed a famous solution for spacecraft in near-circular orbits \cite{clohessy_terminal_1960}, whereas Yamanaka and Ankersen (YA) for spacecraft in eccentric orbits \cite{yamanaka_new_2002}. All these solutions share common structure: $\boldsymbol{\chi}_{k}(t) \approx \boldsymbol{\Psi}(t) \boldsymbol{c}_{k}$, where, $\boldsymbol{c}_k = \{ c_{1,k},c_{2,k},c_{3,k},c_{4,k},c_{5,k},c_{6,k} \} \in \mathbb{R}^{6}$ is the $k^{th}$ spacecraft IC state, and $\boldsymbol{\Psi}$ is the fundamental matrix solution of the linear ODE. The fundamental matrix depends on time, usually through the reference orbit's true anomaly $\nu = \nu(t)$, but does not depend on the IC since the integrated ODE is linear. In this paper, the chosen IC state is a tuned version of the YA's IC,  as represented by the map in Eq. \ref{ICYA} in Appendix \ref{App:Dyn}.  The associated $\boldsymbol{\Psi}$ can be expressed as a function of the reference OE as
\begin{equation}\label{CARTtoIC}
\resizebox{0.89\linewidth}{!}{$
\boldsymbol{\Psi}(t) = a \left[ \begin{array}{c;{2pt/2pt}c;{2pt/2pt}c;{2pt/2pt}c;{2pt/2pt}c;{2pt/2pt}c}
\frac{1}{\rho} - \frac{3}{2} \frac{e}{\eta^3} \sin(\nu) n t & 0 & -\cos(\lambda_\nu) & -\sin(\lambda_\nu) & 0 & 0 \\ \hdashline[2pt/2pt]
- \frac{3}{2} \frac{\rho}{\eta^3} n t & \frac{1}{\rho} & (\frac{1}{\rho}+1)\sin(\lambda_\nu) & -(\frac{1}{\rho}+1)\cos(\lambda_\nu) & 0 & 0 \\ \hdashline[2pt/2pt]
0 & 0 & 0 & 0 & \frac{1}{\rho} \sin(\lambda_\nu) & -\frac{1}{\rho} \cos(\lambda_\nu) \\ \hdashline[2pt/2pt]
\substack{ -\frac{1}{\rho^2} \dot{\rho} - \frac{3}{2} \frac{e}{\eta^3} \cos(\nu) \dot{\nu} n t \\ - \frac{3}{2} \frac{e}{\eta^3} \sin(\nu) n } & 0 & \sin(\lambda_\nu) \dot{\nu} & -\cos(\lambda_\nu) \dot{\nu} & 0 & 0 \\ \hdashline[2pt/2pt]
- \frac{3}{2} \frac{\dot{\rho}}{\eta^3} n t - \frac{3}{2} \frac{\rho}{\eta^3} n  & -\frac{1}{\rho^2}\dot{\rho} & \substack{-\frac{1}{\rho^2}\dot{\rho}\sin(\lambda_\nu) \\ + (\frac{1}{\rho}+1)\cos(\lambda_\nu) \dot{\nu}} & \substack{\frac{1}{\rho^2}\dot{\rho}\cos(\lambda_\nu) \\ + (\frac{1}{\rho}+1)\sin(\lambda_\nu) \dot{\nu}} & 0 & 0 \\ \hdashline[2pt/2pt]
0 & 0 & 0 & 0 & \substack{ -\frac{1}{\rho^2}\dot{\rho}\sin(\lambda_\nu) \\ + \frac{1}{\rho} \cos(\lambda_\nu) \dot{\nu} } & \substack{ \frac{1}{\rho^2}\dot{\rho}\cos(\lambda_\nu) \\ + \frac{1}{\rho} \sin(\lambda_\nu) \dot{\nu} } \\ 
\end{array} \right]$}
\end{equation}
where $\nu = \nu(t)$, $\lambda_\nu = \omega + \nu(t)$, $\rho = 1 + e \cos(\nu(t))$, $\dot{\rho} = - e \sin(\nu) \dot{\nu}$, $\eta = \sqrt{1-e^2}$, $n = \sqrt{\mu/a^3}$, and $\mu$ is the main gravitation parameter. $\boldsymbol{c}_k$ is non-dimensional, which makes the first three rows of $\boldsymbol{\Psi}$ have dimension of length, and the last three rows of length over time. $\boldsymbol{\Psi}$ models both bounded periodic motion and unbounded (unstable) motion. If $c_{1,k} \neq 0$, the motion is unstable in both along-track and radial direction if the reference orbit is eccentric, or just in along-track direction if the orbit is near-circular. If $c_{1,k} = 0$ the motion is bounded and periodic. As mentioned, the relative motion IC can be mapped at the first-order to a set of ROE, linear or nonlinear combination of the OE of two spacecraft or of a reference and a spacecraft. Eq. \ref{MapICtoROE} in Appendix \ref{App:Dyn} provides a linear map between the used IC and the quasi-nonsigular ROE.

\subsection{Closed-form solutions of Passive Safety in Eccentric Orbits}\label{Sec_Cl}
\begin{figure}[ht!]
    \centering
    \includegraphics[width=0.75\linewidth]{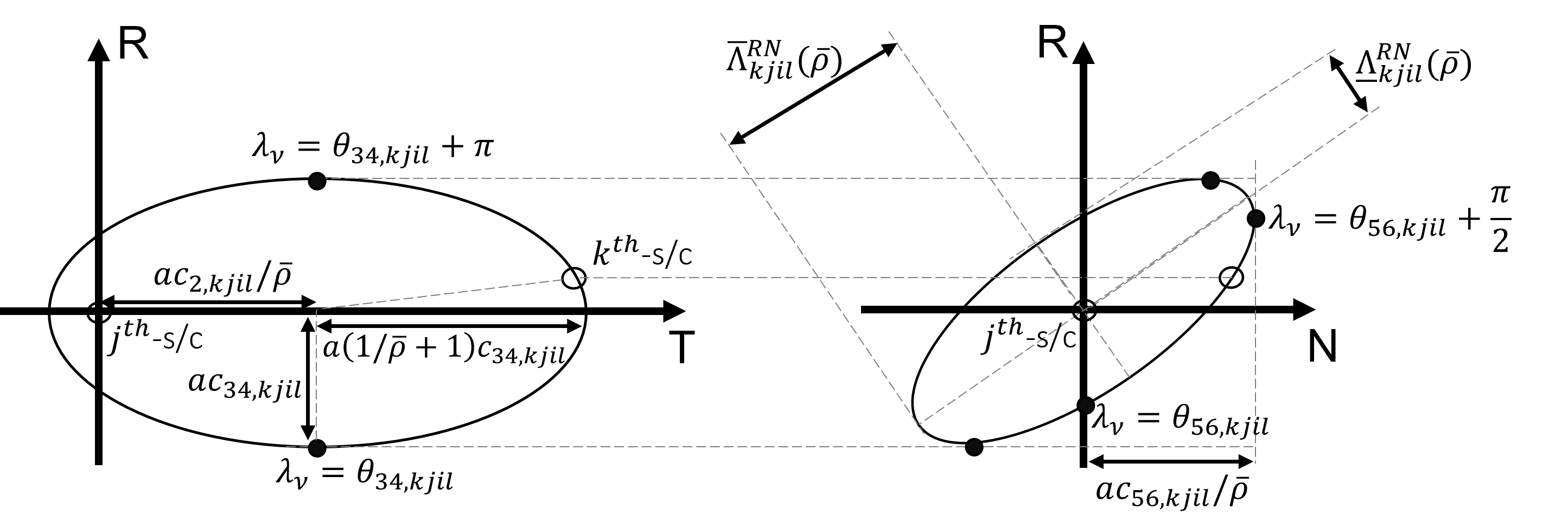}
     \caption{Geometrical representation of Eq. \ref{CARTtoIC2} fixing $\rho(\nu(t)) = \bar{\rho}$}\label{Fig_ROE_Cart}
\end{figure}

Closed-form analysis aims at computing an analytical expression of $\mathcal{S}^*_{kjil}$, or of a lowerbound $\underline{\mathcal{S}}^*_{kjil} \leq \mathcal{S}^*_{kjil}$,  as a function of the IC. To that end, Eq. \ref{Eq_tstar} has to be solved analytically using $\boldsymbol{\Psi}$ in Eq. \ref{CARTtoIC}. The full $3$-D solution for relative motion in eccentric orbits, accounting for both bounded and unbounded motion, is rather complicated. Therefore, literature \cite{damico_proximity_2006, koenig_2018_robust} has been developing reduced solutions by: 1) assuming $c_{1, k} \sim 0$, $\forall k$, and so enforcing the relative motion  to be bounded and periodic for all spacecraft, 2) decomposing the full $3$-D problem in the two $2$-D  radial/normal (RN) and radial/tangential (RT) sub-problems, 3) assuming near-circular orbits. In this section, the third assumption is removed, and two analytical lowerbounds, $\underline{\mathcal{S}}^{* RN}_{kjil}$ and $\underline{\mathcal{S}}^{* RT}_{kjil}$, are computed for relative motion in eccentric orbits for the first time. Here, the superscripts RN/RT refer to the sub-planes in which the minimum separation is computed. To that end, for compactness,  use the notations $\boldsymbol{c}_{kjil} = \boldsymbol{c}_{k}(t_i) - \boldsymbol{c}_{j}(t_l)$ and $\epsilon_{kjil} = \epsilon + \beta_{kjil}$. Assuming $c_{1, k} \sim 0$ and using Eq. \ref{CARTtoIC}, the first three components of $\boldsymbol{\Psi}(t) \boldsymbol{c}_{kjil}$ are
\begin{equation}\label{CARTtoIC2}
\begin{aligned}
& \delta r_{R,kjil}(t) \approx - a c_{34, kjil} \cos(\lambda_\nu(t) - \theta_{34, kjil})\\
& \delta r_{T,kjil}(t) \approx \frac{1}{\rho(\nu(t))} a c_{2, kjil} + \Big(\frac{1}{\rho(\nu(t))}+1\Big) a c_{34, kjil} \sin(\lambda_\nu(t) - \theta_{34, kjil})\\
& \delta r_{N,kjil}(t) \approx \frac{1}{\rho(\nu(t))} a c_{56, kjil} \sin(\lambda_\nu(t) - \theta_{56, kjil})
\end{aligned}
\end{equation}
which represents the relative position of $k^{th}$ and $j^{th}$ spacecraft along the uncontrolled trajectories generated by control losses at instant $t_i$ and $t_l$ respectively, assuming the integrability of the dynamics  after a contingency. In Eq. \ref{CARTtoIC2}, the IC state has been expressed in polar form $\boldsymbol{c}_{ab, kjil} = \{ c_{a,kjil}, c_{b,kjil}\} = c_{ab,kjil}\{\cos(\theta_{ab,kjil}), \sin(\theta_{ab,kjil})\}$. Looking at Eq. \ref{CARTtoIC2}, the terms function of $\rho(\nu(t))$ hinder a clear geometrical interpretation of the relative motion. In particular, $\rho(\nu) = 1 + e \cos(\nu) \in [1-e, 1+e]$, where the two extreme values are achieved at the apogee and at the perigee of the eccentric orbit respectively. If $\rho(\nu)$ is fixed to $\bar{\rho} \in [1-e, 1+e]$, Eq. \ref{CARTtoIC2} has geometrical representation as in Fig. \ref{Fig_ROE_Cart}.
 In near-circular orbits, $\rho(\nu(t)) \approx 1$ $\forall t$, therefore Fig. \ref{Fig_ROE_Cart} becomes the true representation of the relative motion. In eccentric orbits, instead, at each $\nu(t) = \bar{\nu}$, the real relative motion trajectory is tangent to the ellipses in Fig. \ref{Fig_ROE_Cart} at the point $\lambda_\nu = \omega + \bar{\nu}$, with $\bar{\rho} = 1 + e \cos(\bar{\nu})$. This implies that, a lowerbound of the real separation between spacecraft on the uncontrolled trajectories can be computed in either RN or RT planes by evaluating the minimum separations of the two ellipses from the origin. As presented in the following, these minimum separations assume  extreme values at the perigee and apogee of the eccentric orbit. The minima of these  extremes are by definition lowerbounds of $\mathcal{S}^*_{kjil}$, and can be used in Eq. \ref{Eq_PS_SuffC} to achieve PS.  Since $c_{1, k} \sim 0$, the uncontrolled relative motion is bounded and periodic, therefore PS is achieved infinitely ($T$=$\infty$). Note that, if assuming $c_{1, k} \sim 0$ is too limiting for a specific application, it can be avoided by following the numerical path to $t^*_{kjil}$ computation presented in Algorithm \ref{ALG_C} (lines 9-24).

\subsubsection{Passive Safety in RN-plane}
The minimum separation of the RN ellipse from the origin is computed as the minimum eigenvalue $\underline{\Lambda}^{RN}_{kjil}(\bar{\rho})$ of a singular value decomposition (Fig. \ref{Fig_ROE_Cart}, right) as 
\begin{equation}\label{drRNoptecc3}
\underline{\mathcal{S}}^{* RN}_{kjil}(\bar{\rho}) = \frac{\sqrt{2} a \left|\boldsymbol{c}_{34, kjil} \cdot \frac{\boldsymbol{c}_{56, kjil}}{\bar{\rho}}\right|}{\left[(c_{34, kjil})^2 + \left(\frac{c_{56, kjil}}{\bar{\rho}}\right)^2 + \left|\boldsymbol{c}_{34, kjil} + \frac{\boldsymbol{c}_{56, kjil}}{\bar{\rho}}\right| \cdot \left|\boldsymbol{c}_{34, kjil} - \frac{\boldsymbol{c}_{56, kjil}}{\bar{\rho}}\right|\right]^{1/2}}
\end{equation}
Consistently with  Eq. \ref{Eq_constr1}, $\underline{\mathcal{S}}^{* RN}_{kjil}(\bar{\rho})$ is, for a fixed value of $\bar{\rho}$, independent  of time and a function of the IC states evaluated at the contingency instants. Taking the partials with respect to $\bar{\nu}$, its smallest and largest values are achieved correspondingly at the perigee and apogee of the eccentric orbit. Therefore, a lowerbound over the full orbit is obtained substituting $\bar{\rho}$ with $1+e$. Using Eq. \ref{drRNoptecc3}, the sufficient condition of PS in Eq. \ref{Eq_PS_SuffC} can be rewritten as
\begin{equation}\label{dediSEPecc}
\left|\hat{\boldsymbol{c}}_{34, kjil} \cdot \hat{\boldsymbol{c}}_{56, kjil}\right| = \left|\cos(\theta_{34, kjil} - \theta_{56, kjil})\right| \geq \cos(\xi_{kjil}(\bar{\rho}))
\end{equation}
with 
\begin{equation}\label{dediSEPecc2}
\cos(\xi_{kjil}(\bar{\rho})) = \frac{\epsilon_{kjil} \bar{\rho}}{a c_{34, kjil} c_{56, kjil}}\sqrt{(c_{34, kjil})^2 + \left(\frac{c_{56, kjil}}{\bar{\rho}}\right)^2 - \frac{(\epsilon_{kjil})^2}{a^2}}
\end{equation}
$\forall t_i \in [t_0, t_f]$, $\forall t_l \in [t_0, t_f]$, for $k = 1,\dots N-1$, $j = k+1,\dots N$. Here, $\boldsymbol{c}_{ab, kjil} = c_{ab, kjil} \hat{\boldsymbol{c}}_{ab, kjil}$. Eq. \ref{dediSEPecc}-\ref{dediSEPecc2} have an interesting geometrical interpretation in IC space, as presented in Fig. \ref{Fig_PS_RN}. In particular, to maintain PS in RN-plane the phase angle between the vectors $\boldsymbol{c}_{34, kjil}$ and $\boldsymbol{c}_{56, kjil}$, i.e. $\theta_{34, kjil} - \theta_{56, kjil}$, must be either contained in the interval [$-\xi_{kjil}(\bar{\rho})$, $\xi_{kjil}(\bar{\rho})$] (parallel configuration) or in the interval [$\pi - \xi_{kjil}(\bar{\rho})$, $\pi + \xi_{kjil}(\bar{\rho})$] (anti-parallel configuration), where the most and least tight conditions are achieved correspondingly at the perigee and apogee of the eccentric orbit. Eq. \ref{dediSEPecc}-\ref{dediSEPecc2} represent the $\boldsymbol{c}_{34}$/$\boldsymbol{c}_{56}$-separation concept, which in near-circular orbits, using ROE, becomes the relative eccentricity/inclination-separation concept \cite{damico_proximity_2006}, used for the relative orbit design in the GRACE, TANDEM-X, PRISMA and BIROS missions flown in  low-Earth orbit \cite{Errico_13}.

\subsubsection{Passive Safety in RT-plane}
\begin{figure}[ht!]
    \centering 
    \subfloat[PS in RN plane, Cartesian space (left),  integration constants space (right) \label{Fig_PS_RN}]{\includegraphics[width=0.46\textwidth]{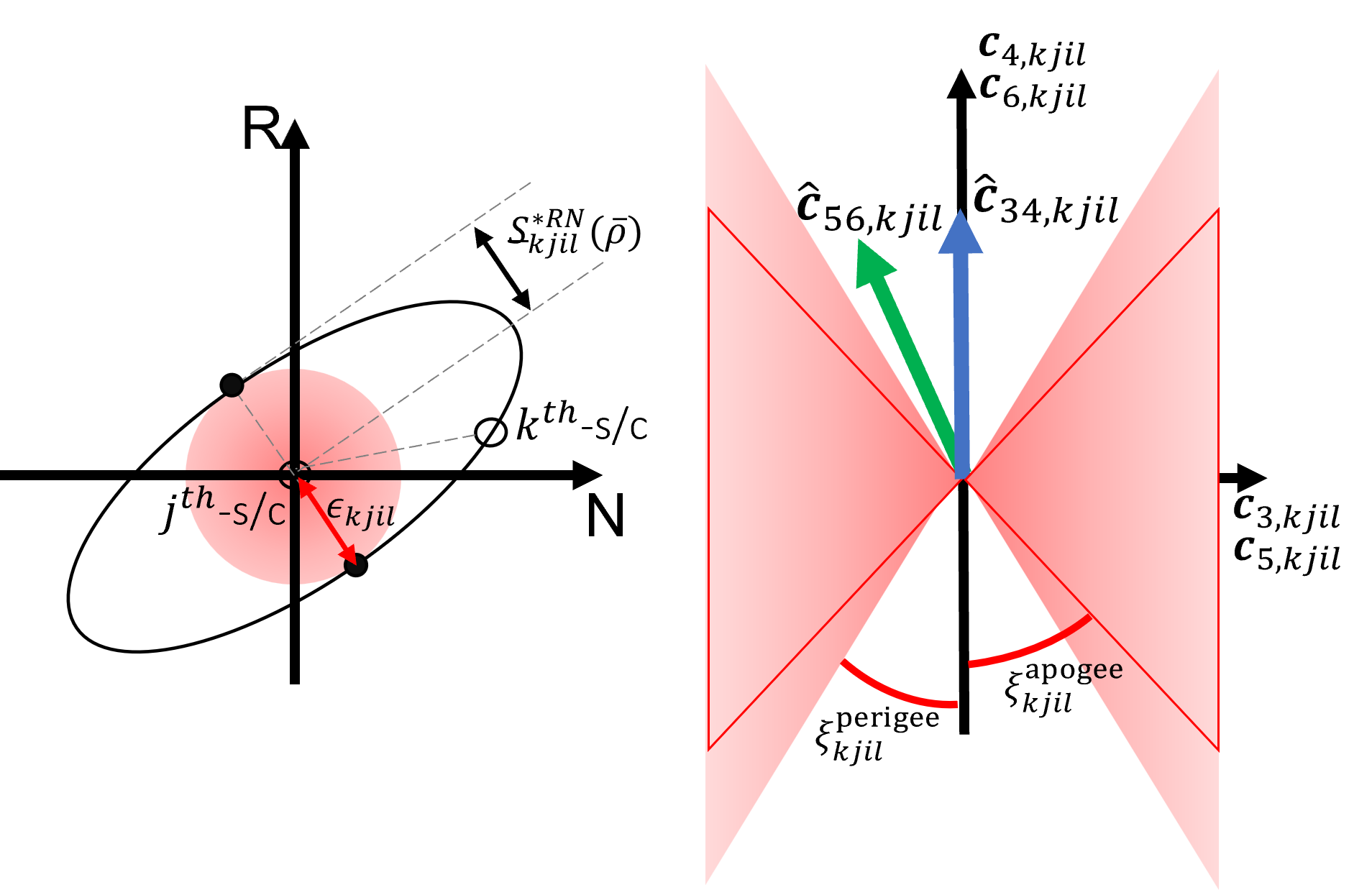}}\hfill
    \subfloat[PS in RT plane, Cartesian space (left),  integration constants space (right) \label{Fig_PS_RT}] {\includegraphics[width=0.54\textwidth]{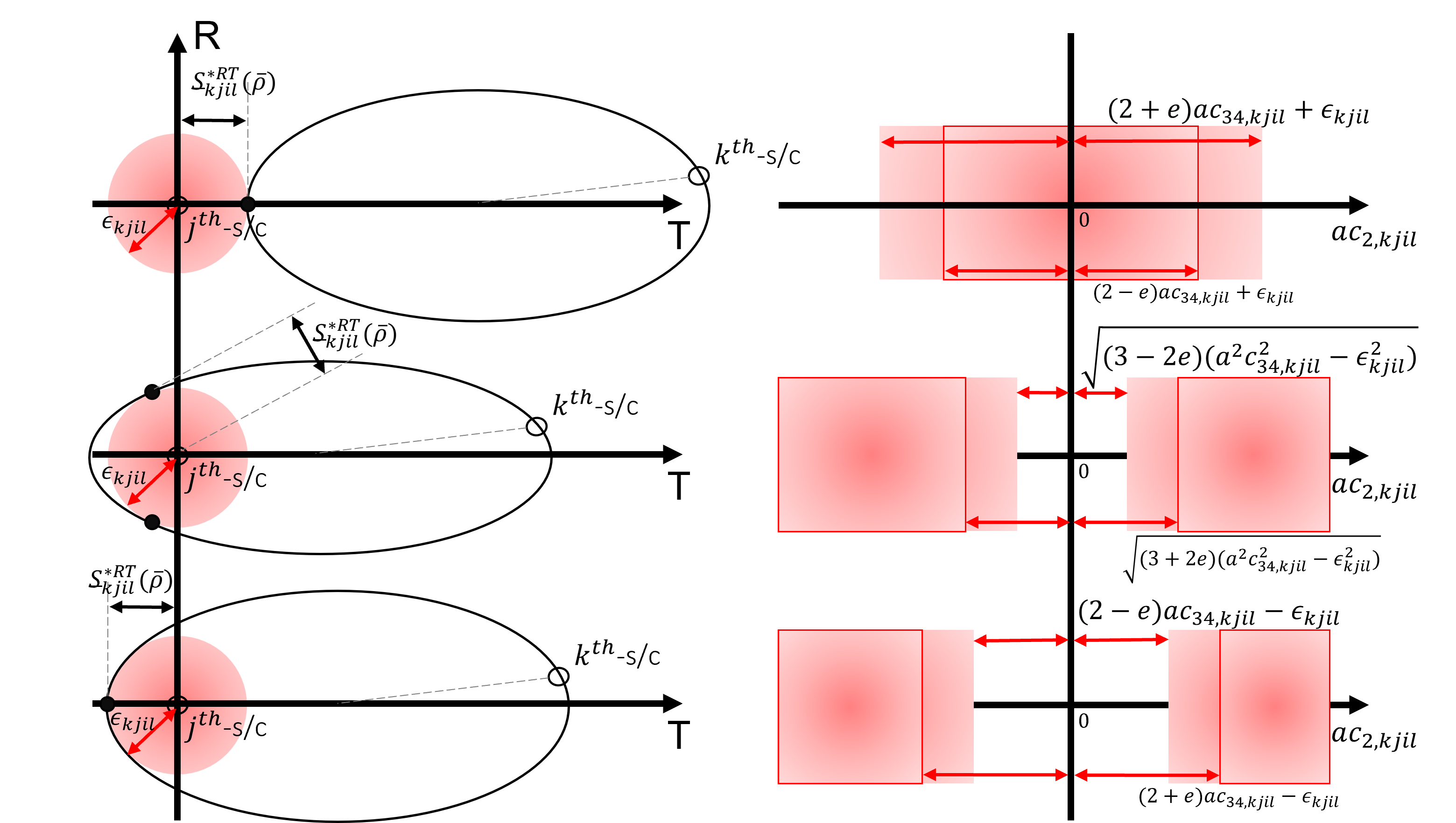}}\hfill
    \caption{Graphical representation of the sufficient conditions of  passive safety (PS) in eccentric orbits} \label{Fig_PS_RN_RT}
\end{figure}

The minimum separation of the RT ellipse from the origin is computed  by solving analytically Eq. \ref{Eq_tstar} formulated as a function of $\lambda_\nu (t)$, and by using Eq. \ref{CARTtoIC2} and fixing $\bar{\rho}$. The resulting minimum separation is
\begin{equation}\label{drRToptecc3}
\begin{aligned}
\underline{\mathcal{S}}^{* RT}_{kjil}(\bar{\rho}) \left\{ \begin{array}{ll}  = a \frac{1}{\bar{\rho}} c_{2, kjil} - a \Big( \frac{1}{\bar{\rho}} + 1 \Big) c_{34, kjil}  & \textrm{not-encircling} \\
 = \left\{ \begin{array}{ll} \sqrt{a^2 c_{34, kjil}^2 - \frac{1}{1+2\bar{\rho}} a^2 c_{2, kjil}^2} & \frac{\epsilon_{kjil}}{a} \leq c_{34, kjil} < \frac{1 + \bar{\rho}}{\bar{\rho}} \frac{\epsilon_{kjil}}{a} \\ a \Big( \frac{1}{\bar{\rho}} + 1 \Big) c_{34, kjil} - a \frac{1}{\bar{\rho}} c_{2, kjil} & c_{34, kjil} \geq \frac{1+ \bar{\rho}}{\bar{\rho}} \frac{\epsilon_{kjil}}{a}
\end{array} \right. & \textrm{encircling}
 \end{array} \right.
\end{aligned}
\end{equation}
where "not-encircling" refers to Fig. \ref{Fig_PS_RT}, top left, whereas, "encircling" refers to Fig. \ref{Fig_PS_RT}, center and bottom left. Consistently with  Eq. \ref{Eq_constr1}, $\underline{\mathcal{S}}^{* RT}_{kjil}(\bar{\rho})$ is, for a fixed value of $\bar{\rho}$, independent  of time and a function of the IC states evaluated at the contingency instants. Taking the partials with respect to $\bar{\nu}$, its smallest and largest values are achieved correspondingly at the perigee and apogee of the eccentric orbit for the not-encircling case, and at the apogee and perigee for the encircling case. Therefore, a lowerbound over the full orbit can be obtained by substituting $\bar{\rho}$ with $1+e$ for the not-encircling case, whereas with $1-e$ for the encircling case. Using Eq. \ref{drRToptecc3}, the sufficient condition of PS in Eq. \ref{Eq_PS_SuffC} can be rewritten as
\begin{equation}\label{RTsepecc}
\begin{aligned}
|c_{2, kjil}| \left\{ \begin{array}{ll} \geq (1 + \bar{\rho}) c_{34, kjil} + \frac{\epsilon_{kjil}}{a} & \textrm{not-encircling} \\
\leq \left\{  \begin{array}{lll} \sqrt{(1 + 2 \bar{\rho})\left(c_{34, kjil}^2 - \frac{(\epsilon_{kjil})^2}{a^2}\right)} & \textrm{if} & \frac{\epsilon_{kjil}}{a} \leq c_{34, kjil} < \frac{1 + \bar{\rho}}{\bar{\rho}} \frac{\epsilon_{kjil}}{a} \\ 
(1 + \bar{\rho}) c_{34, kjil} - \frac{\epsilon_{kjil}}{a} & \textrm{if} & c_{34, kjil} \geq \frac{1 + \bar{\rho}}{\bar{\rho}} \frac{\epsilon_{kjil}}{a} \end{array} \right. & \textrm{encircling}
 \end{array} \right.
\end{aligned}
\end{equation}
$\forall t_i \in [t_0, t_f]$, $\forall t_l \in [t_0, t_f]$, for $k = 1,\dots N-1$, $j = k+1,\dots N$. In IC space, Eq. \ref{RTsepecc}  is represented as in Fig. \ref{Fig_PS_RT}, right, where the most tight and least tight conditions are achieved correspondingly at the perigee and apogee for the not-encircling case, and at the apogee and perigee for the encircling case. In near-circular orbits, using ROE, Eq. \ref{RTsepecc} reduces to the solutions in \cite{koenig_2018_robust}.

\subsection{Eccentric Swarms Designs}\label{Sec_swarm_ecc}
The newly derived analytical PS conditions in RN and RT planes  permit the generalization of the spacecraft swarm designs developed in \cite{koenig_2018_robust} from near-circular orbits to eccentric orbits. The "high-density" swarm design is defined such that all spacecraft are distributed at equal distance in the $c_{34}$-plane and RT PS is achieved by enforcing sufficient separation in along-track direction through Eq. \ref{RTsepecc}. To do so,  design the swarm such that: $c_{1, k} = 0$, and $c_{34, kj} = || \boldsymbol{c}_{34, k} - \boldsymbol{c}_{34, j} ||_2  \geq c_{34, min} \geq \epsilon/a$. Moreover, all spacecraft are distributed at equal distance ($c_{34, min}$) in the $c_{34}$-plane as
\begin{equation}\label{SwarmDesignRTecc}
c_{34, k} = \left(\begin{array}{cc} c_{34, min} \left(W_k \cos(\theta^{hd}_{34}) + X_k \cos(\theta^{hd}_{34} + \frac{\pi}{3})\right) \\  c_{34, min} \left(W_k \sin(\theta^{hd}_{34}) + X_k \sin(\theta^{hd}_{34} + \frac{\pi}{3})\right)  \end{array}\right)
\end{equation}
where, the angular phase $\theta^{hd}_{34}$ is selected to arrange the spacecraft in the plane, and $W_k$ and $X_k$ are not-simultaneously zero integers, that cannot simultaneously satisfy $W_k = W_j$ and $X_k = X_j$, for any $k \neq j$. Given this design, by noticing that $|c_{2, kj}| \leq |c_{2, k}| + |c_{2, j}|$, RT PS is enforced for the entire swarm by using Eq. \ref{RTsepecc} evaluated at the apogee (the most tight condition for the encircling case) as
\begin{equation}\label{SwarmDesignRTecc2}
\begin{aligned}
|c_{2, k}| \leq \left\{  \begin{array}{lll} \sqrt{\frac{3 - 2 e}{4} \left(c_{34, min}^2 - \frac{\epsilon^2}{a^2}\right)} & \textrm{if} & \frac{\epsilon}{a} \leq c_{34, min} < \frac{2-e}{1-e} \frac{\epsilon}{a} \\ 
\frac{2-e}{2} c_{34, min} - \frac{\epsilon}{2 a} & \textrm{if} & c_{34, min} \geq \frac{2-e}{1-e} \frac{\epsilon}{a} \end{array} \right.
\end{aligned}
\end{equation}

The RN separated swarm is designed instead such that all the spacecraft share same phase angles of the $\boldsymbol{c}_{34}/\boldsymbol{c}_{56}$ vectors, while having these vector separated according to Eq. \ref{dediSEPecc}-\ref{dediSEPecc2}. To do so,  design the swarm such that: $c_{1, k} = 0$, $c_{34, kj} = || \boldsymbol{c}_{34, k} - \boldsymbol{c}_{34, j} ||_2 \geq c_{34, min} \geq \epsilon/a$, and $c_{56, kj} = || \boldsymbol{c}_{56, k} - \boldsymbol{c}_{56, j} ||_2 \geq c_{56, min} \geq (\epsilon/a)(1+e)$. Moreover, all the spacecraft must share same phase angles of the $\boldsymbol{c}_{34}/\boldsymbol{c}_{56}$ vectors, $\theta^{sep}_{34}$ and $\theta^{sep}_{56}$, as: $\hat{\boldsymbol{c}}_{34, k} \equiv \big\{\cos(\theta^{sep}_{34}), \sin(\theta^{sep}_{34}) \big\}$, and $\hat{\boldsymbol{c}}_{56, k} \equiv \big\{\cos(\theta^{sep}_{56}), \sin(\theta^{sep}_{56})\big\}$. Given this design, RN PS is enforced for the entire swarm by using Eq. \ref{dediSEPecc}-\ref{dediSEPecc2} evaluated at the perigee (most tight condition) as
\begin{equation}\label{SwarmDesignRN5}
\left|\cos(\theta^{sep}_{34} - \theta^{sep}_{56})\right| \geq \frac{\epsilon (1+e)}{a c_{34, min} c_{56, min}}\sqrt{(c_{34, min})^2 + \left(\frac{c_{56, min}}{1+e}\right)^2 - \frac{\epsilon^2}{a^2}}
\end{equation}

Note that for eccentricity going to zero, these two swarms designs reduce exactly to the ones in \cite{koenig_2018_robust}, with the selected IC state reducing to the quasi-nonsingular ROE through the map in Eq. \ref{MapICtoROE} in Appendix \ref{App:Dyn}. These two passively-safe spacecraft swarm designs are showcased in the next section, together with a reconfiguration between them enforcing PS in $3$-D. This capability is not published in literature and is enabled by Algorithm \ref{ALG_C}.

\section{Numerical Results}\label{Sec_VI}
\subsection{The Virtual Super-resolution Optics with Reconfigurable Swarms (VISORS)}
 Due to launch in 2024, VISORS is a distributed telescope  composed by two $6$U CubeSats that will operate in a near-circular sun-synchronous  low-Earth orbit to image the coronal area of the sun with unprecedented resolution \cite{VISORS}. One spacecraft (OSC) contains the optical payload, the other (DSC) contains the detector. During nominal operations, DSC acts as passive chief, whose center of mass sets the  radial/tangential/normal (RTN) frame's origin, whereas OSC actively controls the relative motion. The test cases considered model the challenging "transfer phase", where a passively-safe control plan that moves the formation from the stand-by mode (initial condition) to the science mode (target condition) has to be computed. Both initial and target modes are designed in closed-form to achieve  radial/normal (RN) passive safety (PS) through relative eccentricity/inclination-separation \cite{damico_proximity_2006} (Eq. \ref{dediSEPecc}-\ref{dediSEPecc2}). Nevertheless, PS has not only to be guaranteed at initial and terminal states, but at every instant of the transfer too. Therefore, Algorithm \ref{ALG_C} has to be used to generalize over pure closed-form analysis and compute the transfer. All the parameters of the reconfiguration scenario are reported in Table \ref{tab:scenario}-\ref{tab:dynamics}.
 Two transfers are considered, $2$ and $5$ orbits long, in which the $L_2$-norm of the control input is minimized, and PS is guaranteed in RN plane for at least two orbits ($T$), of at least $3$ or $5$m ($\epsilon$), at $3$-$\sigma$ confidence ($q$). The uncertainty stems from differential GPS navigation \cite{giralo_phd_2021}, cold-gas impulsive actuation \cite{lightsey_18} and process-noise modeling the discrepancy between the full-force ground-truth \cite{VISORS, giralo_phd_2021} and the on-board dynamics model \cite{koenig_2017_new, guffanti_long-term_2017, guffanti_linear_2019}  used for optimization. In particular, looking at Table \ref{tab:scenario} on the right, the relation between maneuver magnitude and associated actuation execution uncertainty is proportional, which motivates the uncertainty model formalized in Section \ref{Sec_II}. Moreover, looking at Table \ref{tab:dynamics} on the right, the process noise has to account for the unmodeled and possibly stochastic effects of differential drag present in  low-Earth orbit.  Finally, note that the slowly time-varying nature of the orbital elements state used in the numerical integration of the Gauss variational equations (GVE) allows for the large integration time steps considered in Table \ref{tab:dynamics}, while retaining high propagation accuracy \cite{sullivan_phd_2020}.
\begin{table*}[ht!] 
\centering
\caption{VISORS mission's reconfiguration scenario}\label{tab:scenario}
\resizebox{\linewidth}{!}{
\begin{tabular}{ c c c c c c c c c c c c c c c }
 \hline 
 \hline
  \multicolumn{2}{c}{Reference orbit} & & \multicolumn{3}{c}{Relative orbit modes} & & \multicolumn{2}{c}{Transfer \& safety} & & \multicolumn{2}{c}{GPS Navigation \cite{giralo_phd_2021}} & & \multicolumn{2}{c}{Impulsive Actuation \cite{lightsey_18}} \\
  \cline{1-2} \cline{4-6} \cline{8-9} \cline{11-12} \cline{14-15}
   Type & Chief & & & Stand-by ($t_0$) & Science ($t_f$) & & $t_f - t_0$ (orbits) & 2 or 5 & & Uncertainty & $1$-$\sigma$ & & Min. Imp. Bit (Ns) & 1e-3 \\
   \cline{5-6}
   $a$ (km) & $R_E$ + [500 or 600] & & $a c_1$ (m) & 10 or 0 & 0 & & PS type & RN & & & per axis & & Min. Imp. Incr. (Ns) & 6e-5 \\
   \cline{11-12} \cline{14-15}
   $e$ (-) & 0.001 & & $a c_2$ (m) & -100 & 10 or 34.76 & & $T$ (orbits) & 2 & & Abs. Position. (m) & 1 & & Uncertainty & $1$-$\sigma$ \\
   \cline{14-15}
   $i$ (deg) & 98 & & $a c_3$ (m) & 0 & 25 or -34.64 & & $\epsilon$ (m) & 3 or 5 & & Abs. Velocity (cm/s) & 1 & & Maneuver Magn. (\%) & 5 \\
   $\Omega$ (deg) & 252 or 330 & & $a c_4$ (m) & 200 & 0 or 2.96 & & $q$ ($\sigma$) & 3 & & Rel. Position (cm) & 1 & & Maneuver Dir. (arcmin) & 1 \\
   $\omega$ (deg) & 0 or 45 & & $a c_5$ (m) & 0 & 15 or -19.72 & & $\mathcal{J}$ (m/s) & $\sum_{ki} || \tilde{\boldsymbol{u}}_{k}(t_i) ||_{2}$ & & Rel. velocity ($\mu$m/s) & 25 & & & \\
   $\nu(t_0)$ (deg) & -90 or 1 & & $a c_6$ (m) & 200 & 0 or 1.69 & & $n_M$ (-) & 25 or 75 & & & & & & \\
   Period (h) & 1.58 & & & & & & $n$ (-) & 720 or 1800 & & & & & & \\
   Epoch & 1-1-2024 & & & & & & $m$ (-) & 720 or 720 & & & & & & \\
   GPS time & 00$^h$:00$^m$:00$^s$ & & & & & & & & & & & & & \\
  \hline
  \hline
\end{tabular}}
\end{table*}
\begin{table*}[ht!]
\centering
\caption{Spacecraft features and dynamics models}\label{tab:dynamics}
\resizebox{\linewidth}{!}{
\begin{tabular}{ c c c c c c c}
 \hline 
 \hline
  \multicolumn{2}{c}{Spacecraft features} & & \multicolumn{4}{c}{Dynamics models} \\
  \cline{1-2} \cline{4-7}
  Mass (kg) & 12 & & & Ground-truth \cite{giralo_phd_2021, VISORS} & On-board prediction \cite{VISORS} & On-board optimization \cite{koenig_2017_new, guffanti_long-term_2017, guffanti_linear_2019} \\
  \cline{5-7}
  Cross-section area ($m^2$) &  nom: 0.285 & &  Type &  GVE numerical integration &  GVE numerical integration &  Analytical \\
  & avg: 0.206 & & Geopotential & GGM01S ($60 \times 60$) \cite{GRACE} & GGM01S ($20 \times 20$) \cite{GRACE} & $J_2$ secular \& long-period \\
  $C_d$ (-) & 2.3 & & Atmospheric drag & NRLMSISE00 \cite{picone2002nrlmsise} & \begin{tabular}{@{}c@{}} NRLMSISE00 \cite{picone2002nrlmsise} \\ (with constant $1000 \%$ error) \end{tabular} & Secular \& long-period model-free  \\
  $C_r$ (-) & 1.88 & & Solar radiation pressure &  Cannonball model, no eclipses & None & \begin{tabular}{@{}c@{}} Secular \& long-period proportional to \\ differential ballistic coefficients \end{tabular} \\
  & & & Third-body gravity & Analytical lunisolar ephemerides & None & \begin{tabular}{@{}c@{}} Luni-solar secular \& long-period \end{tabular} \\
  & & & Integrator & RK4 & RK4 & Euler \\
  & & & Step size & Fixed: 30s & Fixed: 60s & Fixed: see Table \ref{tab:scenario} \\
  & & & Process-noise on IC (m/orbit) & None & None & 3 \\
  \hline
  \hline
\end{tabular}}
\end{table*}

\subsubsection{Two-orbit Transfer}
Algorithm \ref{ALG_C} is used to solve for the  two-orbit transfer. The solution of the convex two-point-boundary-value-problem (TPBVP) \cite{ecos_13} is used as initial warm-start, $MAX\_ITER$ is set to $3$, and trust-regions varying from $100$m to $5$m along the transfer are implemented. The underlying optimization problem is a  second order cone program (SOCP) \cite{boyd_book}. The optimizer looks for candidate impulsive maneuvers every $30$deg along the orbit, so $n_M = 25$, whereas the discretization at which PS must be enforced is at $1$deg, therefore $n = m = 720$. PS is enforced after each maneuver application to guarantee safety both on the trajectory in between maneuvers and, in case of subsequent maneuver failure, afterwards up to $t_f + T$. 
Given the selected orbit discretization, the proposed efficient $\mathcal{M}^*$ formulation  permits the reduction of the number of constraints required to enforce PS from $> 25 \times 720 = 18000$ to $25$. By explicitly enforcing these constraint within the SOCP, this implies in average a computational saving of a factor $>100$ (from $33$s to $0.27$s\footnotemark[3]) on the solution of an instance of  the optimization problem in Eq. \ref{OCP_IC_multi} using a SOCP solver \cite{ecos_13}. The computational overhead introduced by the precomputation of $t^*_{kjil}$ is minimal, since it is performed jointly with the precomputations required to evaluate $\beta_{kjil}$ as well as other matrices needed in the SOCP, for a total execution time of $1$s. Note that, if all the thousands  of constraints are enforced explicitly and $t^*_{kjil}$ is not computed, precomputations are still needed at each  sequential convex programming (SCP) iteration, but may be lighter depending on the implementation. On the other hand, the computational gain obtained on the solution of the SOCP by enforcing just $25$ constraints is much greater. 
\begin{figure}[ht!]
    \centering
    \subfloat[Using $\mathcal{M}$ formulation of  the passive safety constraint \label{fig:advantage_T2_1}]{\includegraphics[width=0.33\linewidth,trim=2 2 4 4,clip]{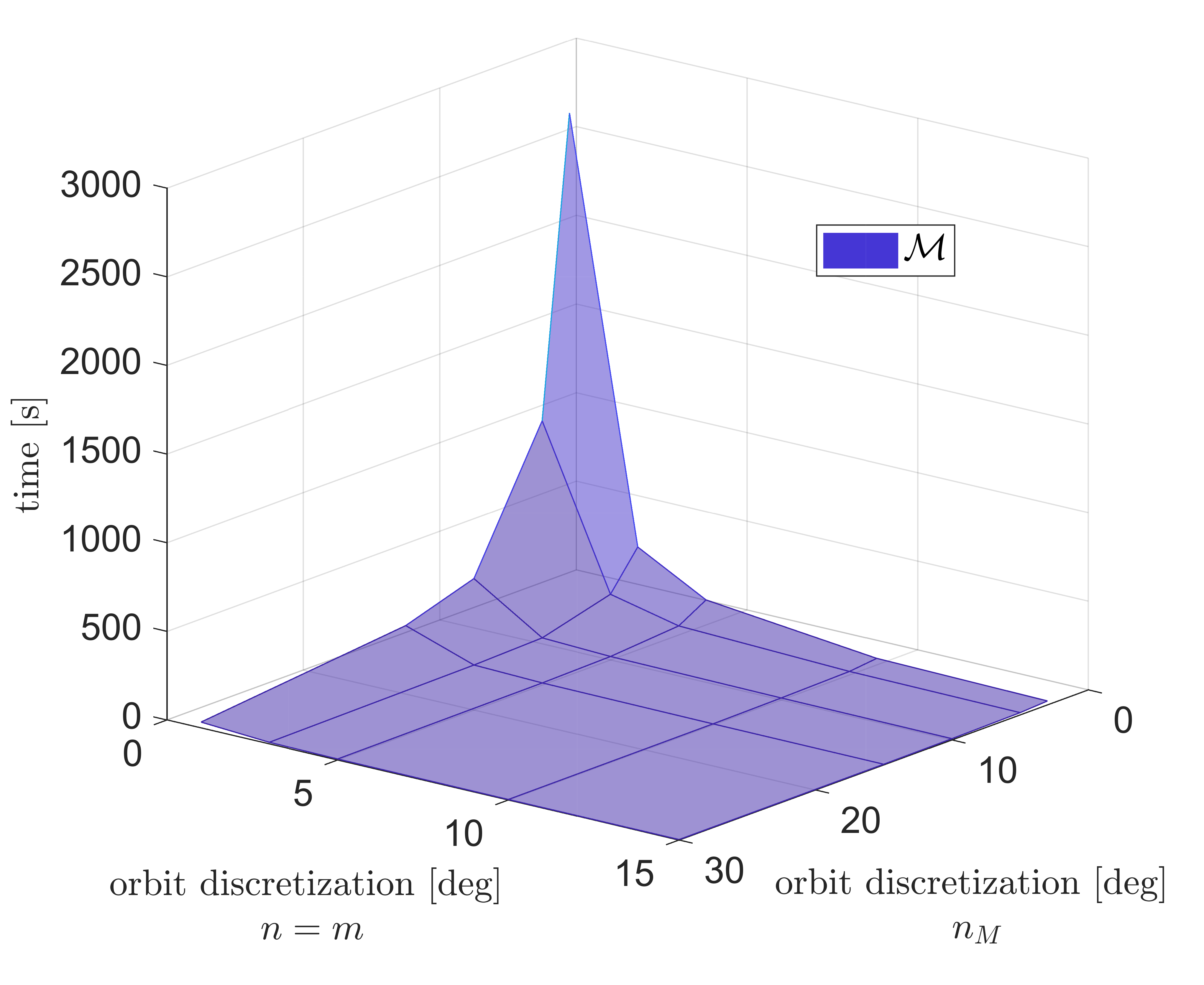}}\hfill \subfloat[Using $\mathcal{M}^*$ formulation of  the passive safety constraint \label{fig:advantage_T2_2}]{\includegraphics[width=0.33\linewidth,trim=2 2 4 4,clip]{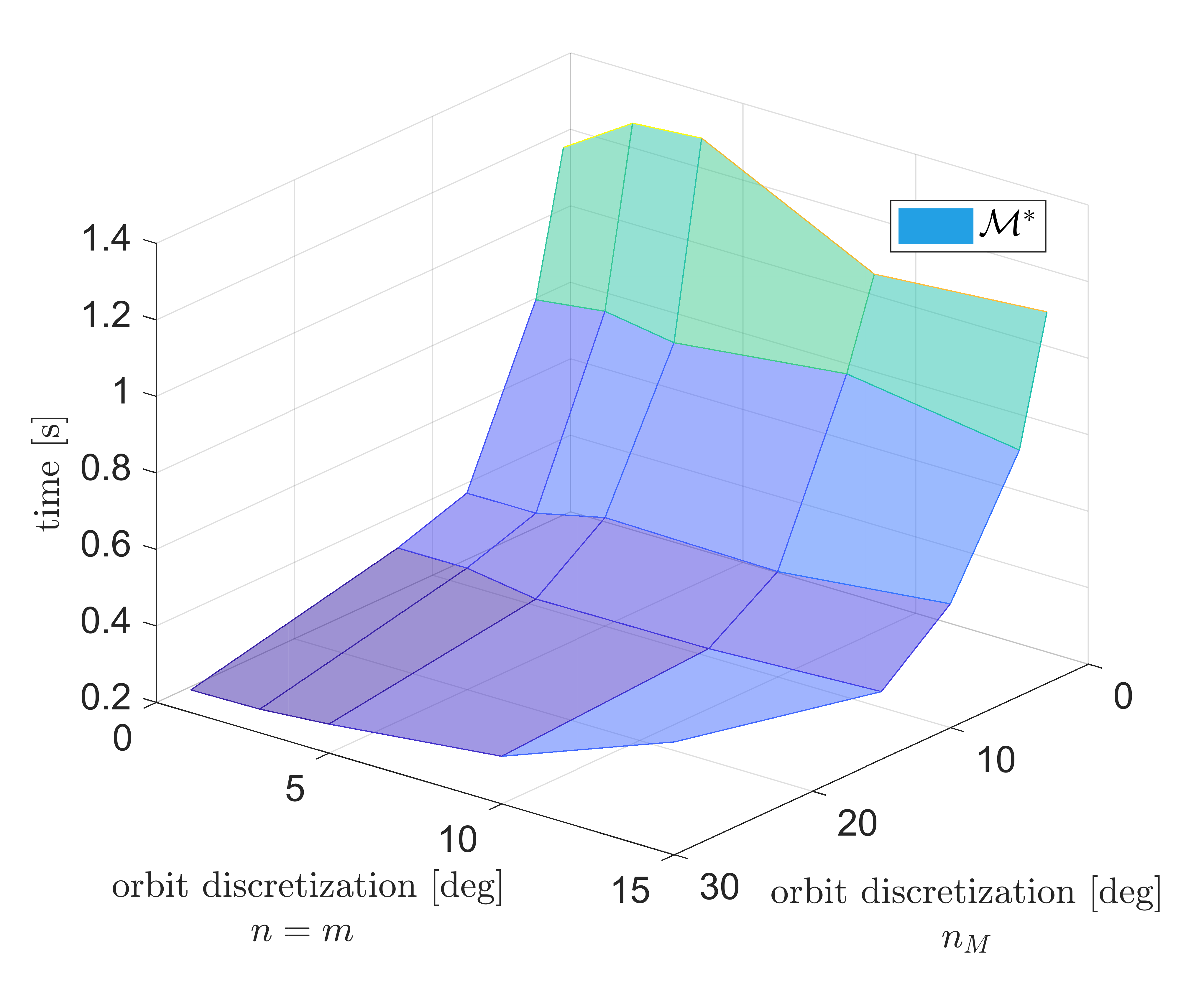}}\hfill
    \subfloat[Corresponding gain factor \label{fig:advantage_T2_3}]{\includegraphics[width=0.33\linewidth,trim=2 2 4 5,clip]{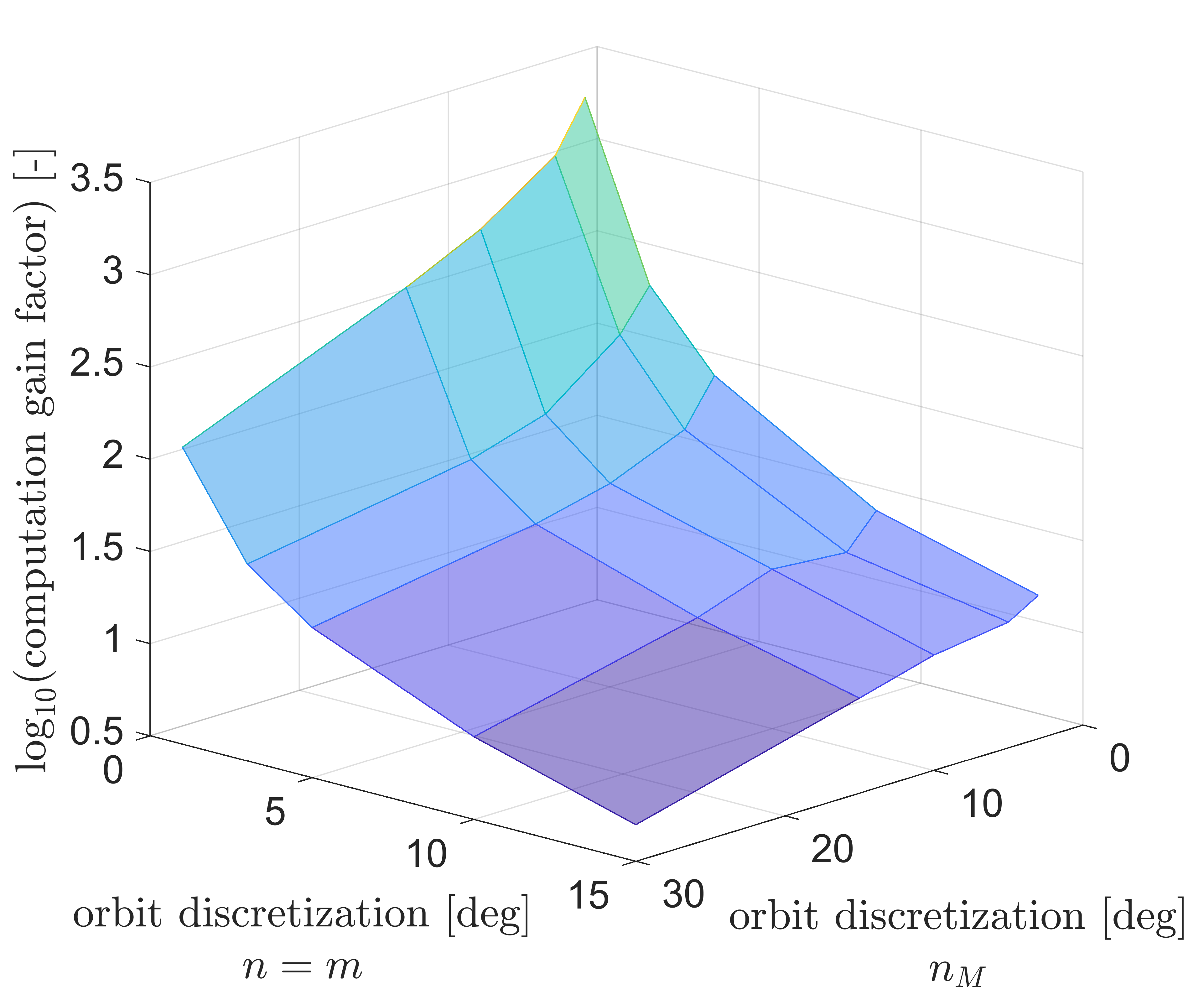}}
    \caption[]{Monte Carlo analysis on the computational effort of solving a  second order cone program} instance \footnotemark[3] \label{fig:advantage_T2}
\end{figure}
\footnotetext[3]{Computation times reported are from a Windows system equipped with a 4GHz processor and 64GB RAM.}
Figure \ref{fig:advantage_T2} enlarges this analysis showing the results of a Monte Carlo on the computational effort  by varying the orbit discretization between maneuvers ($n_M$) and at which safe separation must be enforced ($n = m$). Looking at Figures \ref{fig:advantage_T2_1}-\ref{fig:advantage_T2_2}, the expected correlation between time/orbit discretization, number of constraints, and computation time is clear. In particular, the computational effort trends using the $\mathcal{M}$ and $\mathcal{M}^*$ formulations follow closely the ones of the number of enforced constraints presented in Figure \ref{Fig_constr_num}. Looking at Figure \ref{fig:advantage_T2_3}, the resulting computational saving factor depends on the time/orbit discretization, and it is up to  three orders of magnitude. This validates the benefit of the proposed $\mathcal{M}^*$ formulation of the PS constraint.

\begin{figure}[ht!]
    \centering
    \subfloat[$3$-D RTN space \label{fig:cartresults_T2_1}]{\includegraphics[width=0.58\linewidth,trim=2 2 3 3,clip]{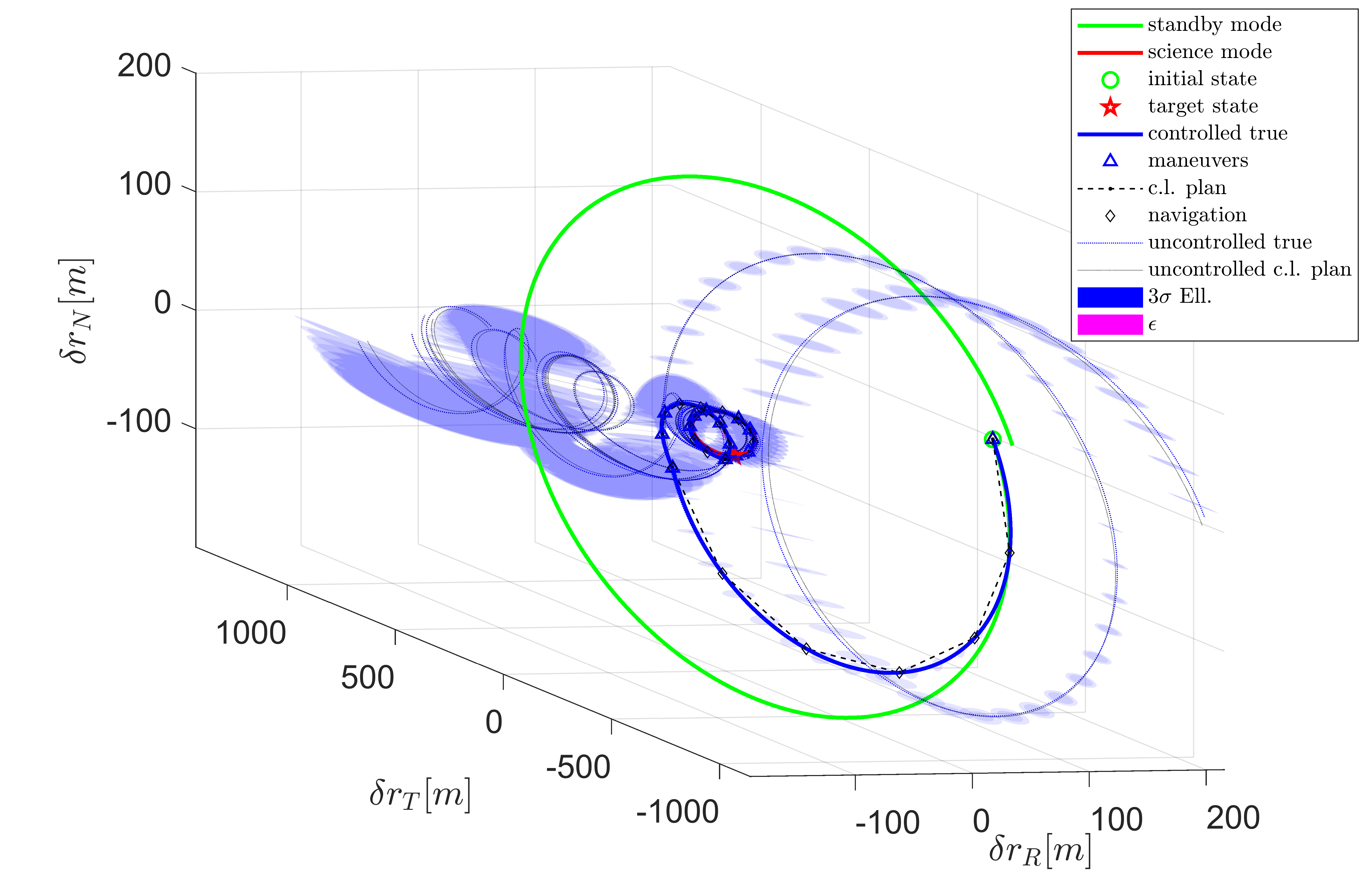}}
    \subfloat[$2$-D RN projection \label{fig:cartresults_T2_2}]{\includegraphics[width=0.38\linewidth,trim=2 2 3 3,clip]{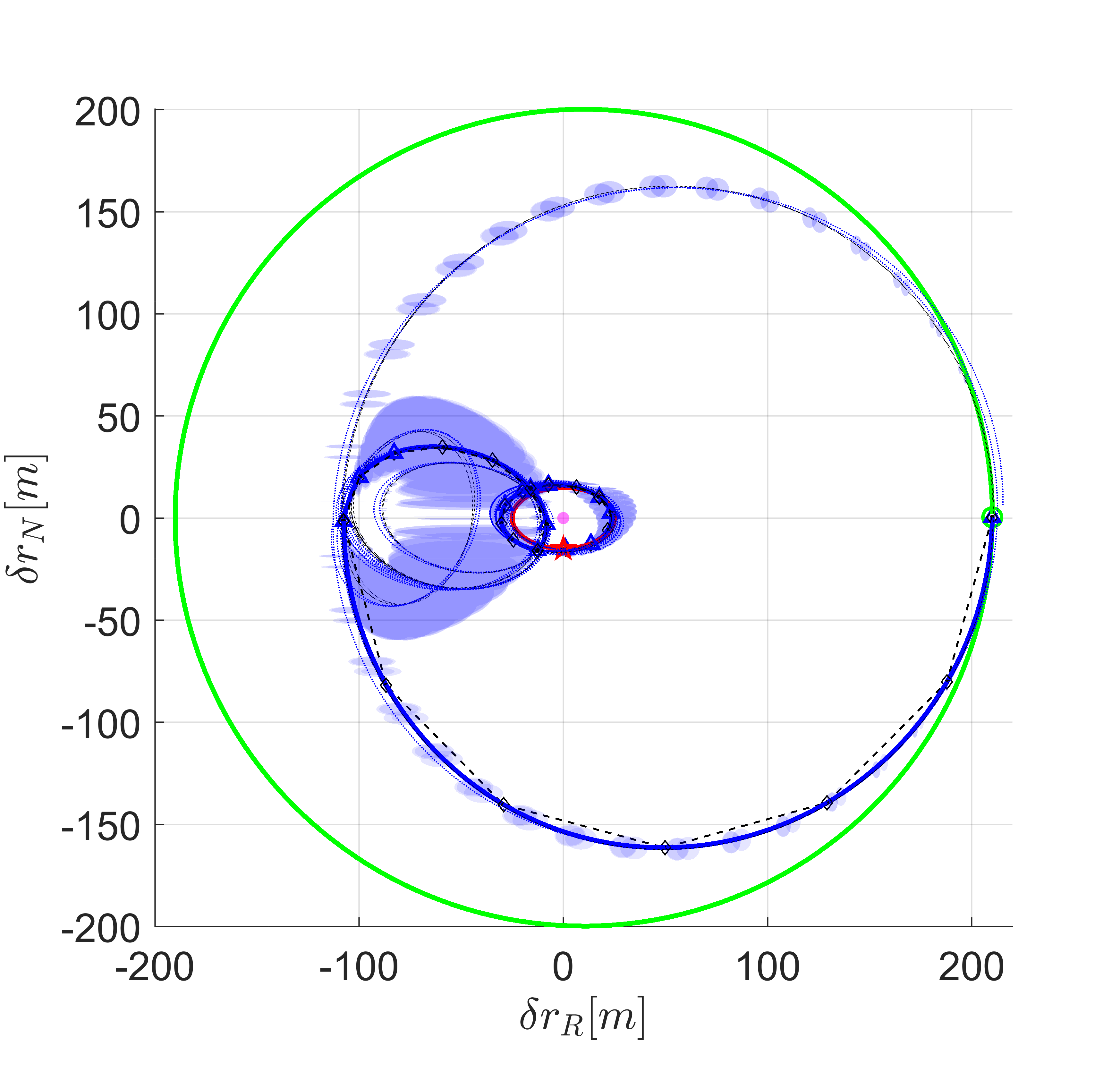}}
     \caption{Two-orbit passively-safe transfer trajectory} \label{fig:cartresults_T2}
\end{figure}

\begin{figure}[ht!]
    \centering
    \subfloat[ Minimum RN separation on uncontrolled trajectories \label{fig:RNsep_T2}]{\includegraphics[width=0.4\linewidth]{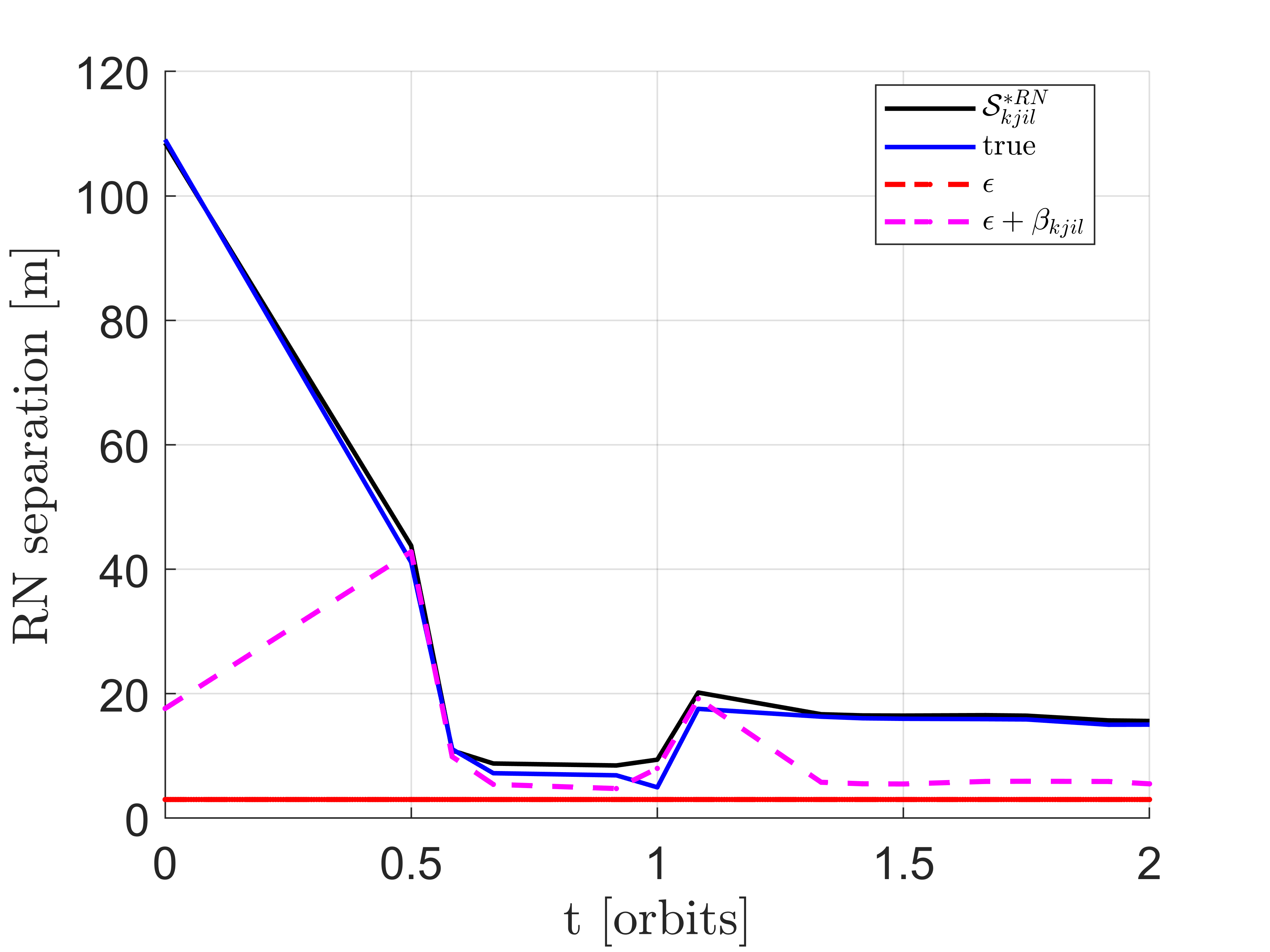}}
    \subfloat[ Minimum $3$-D separation on uncontrolled trajectories \label{fig:3Dsep_T2}]{\includegraphics[width=0.4\linewidth]{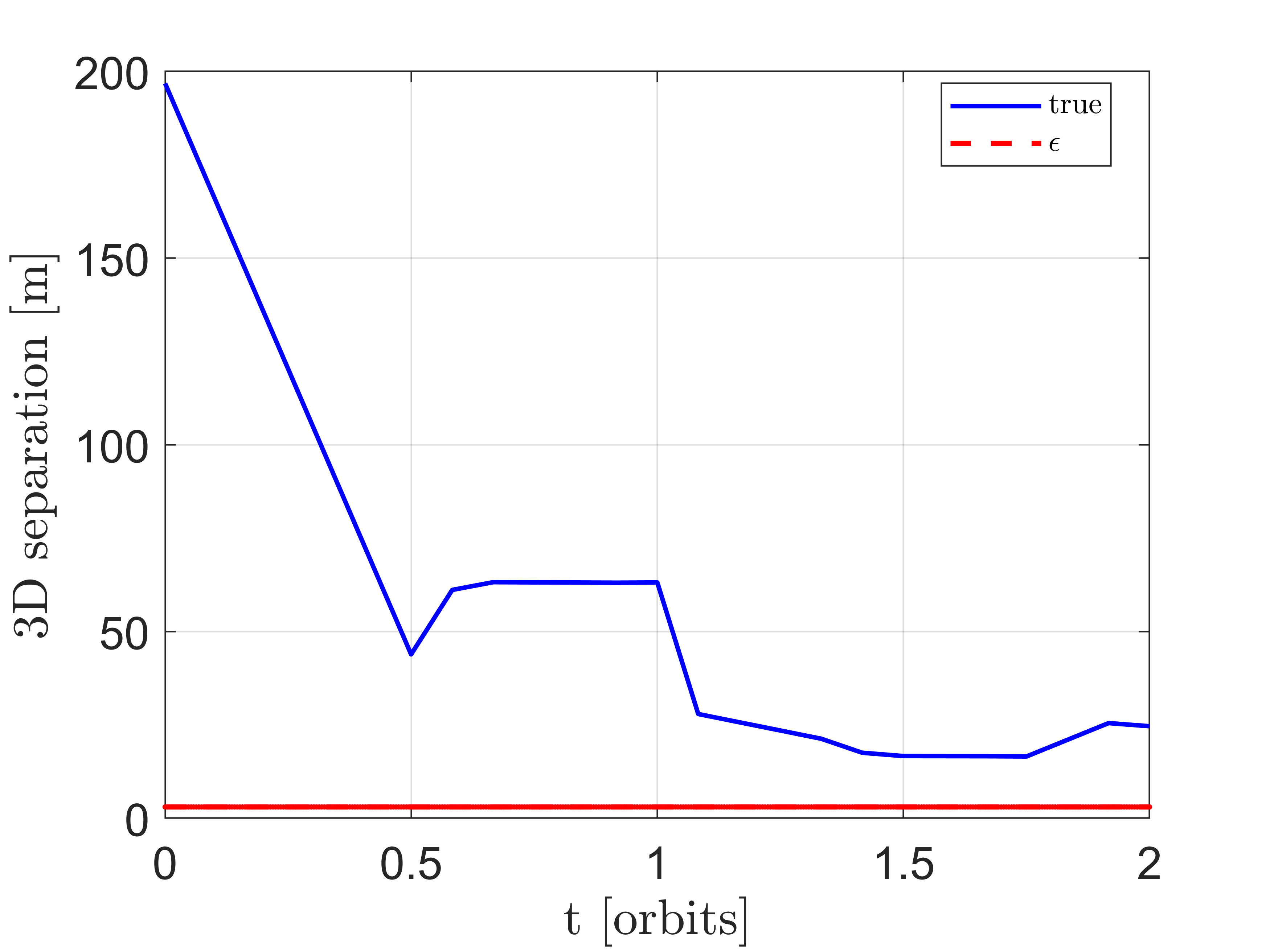}} \\ 
    \centering
    \subfloat[Maneuvers magnitudes \label{fig:ICresults_T2_1}]{\includegraphics[width=0.4\linewidth]{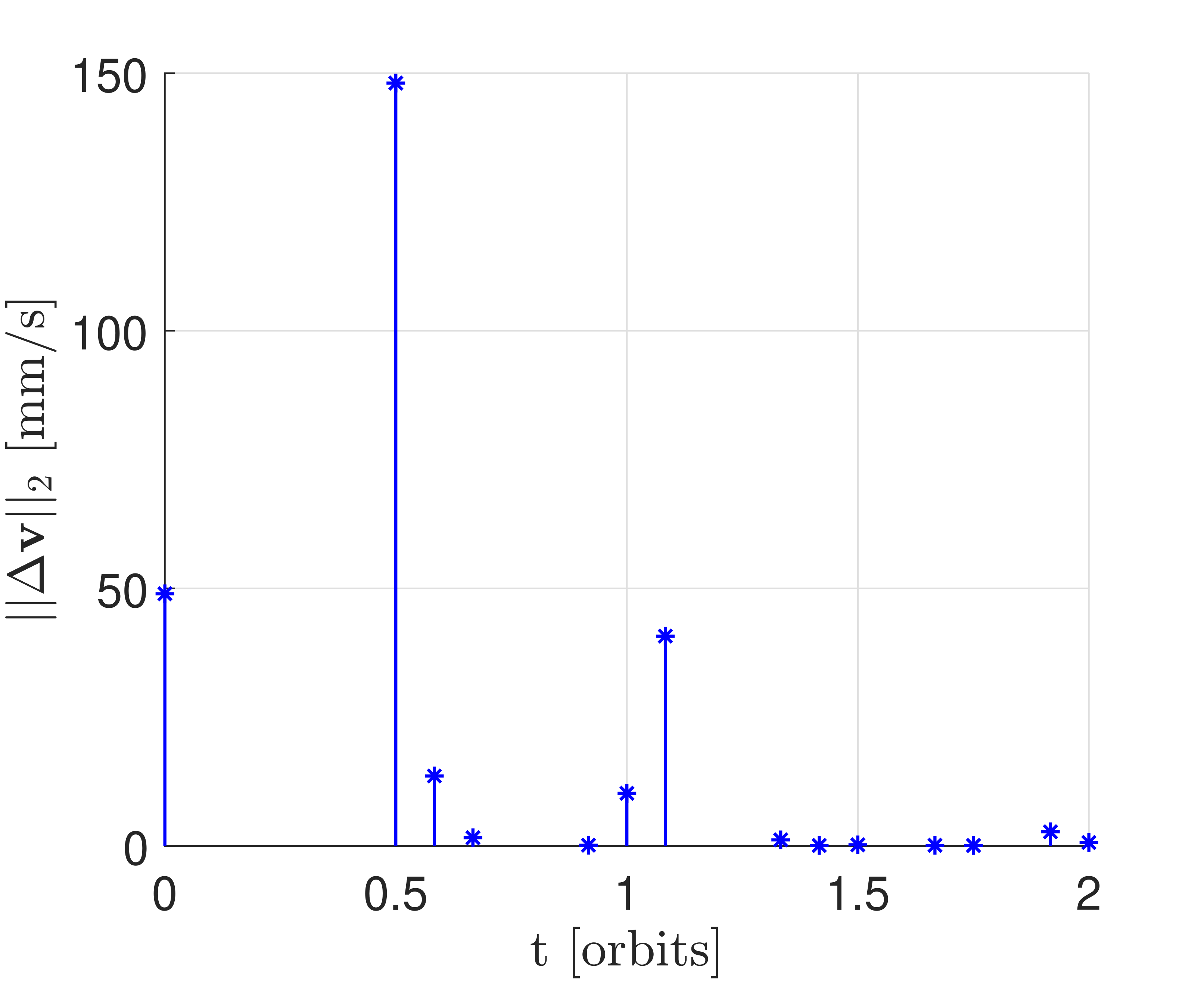}}
    \subfloat[Evolution of $a c_{1}$ \label{fig:ICresults_T2_2}]{\includegraphics[width=0.4\linewidth]{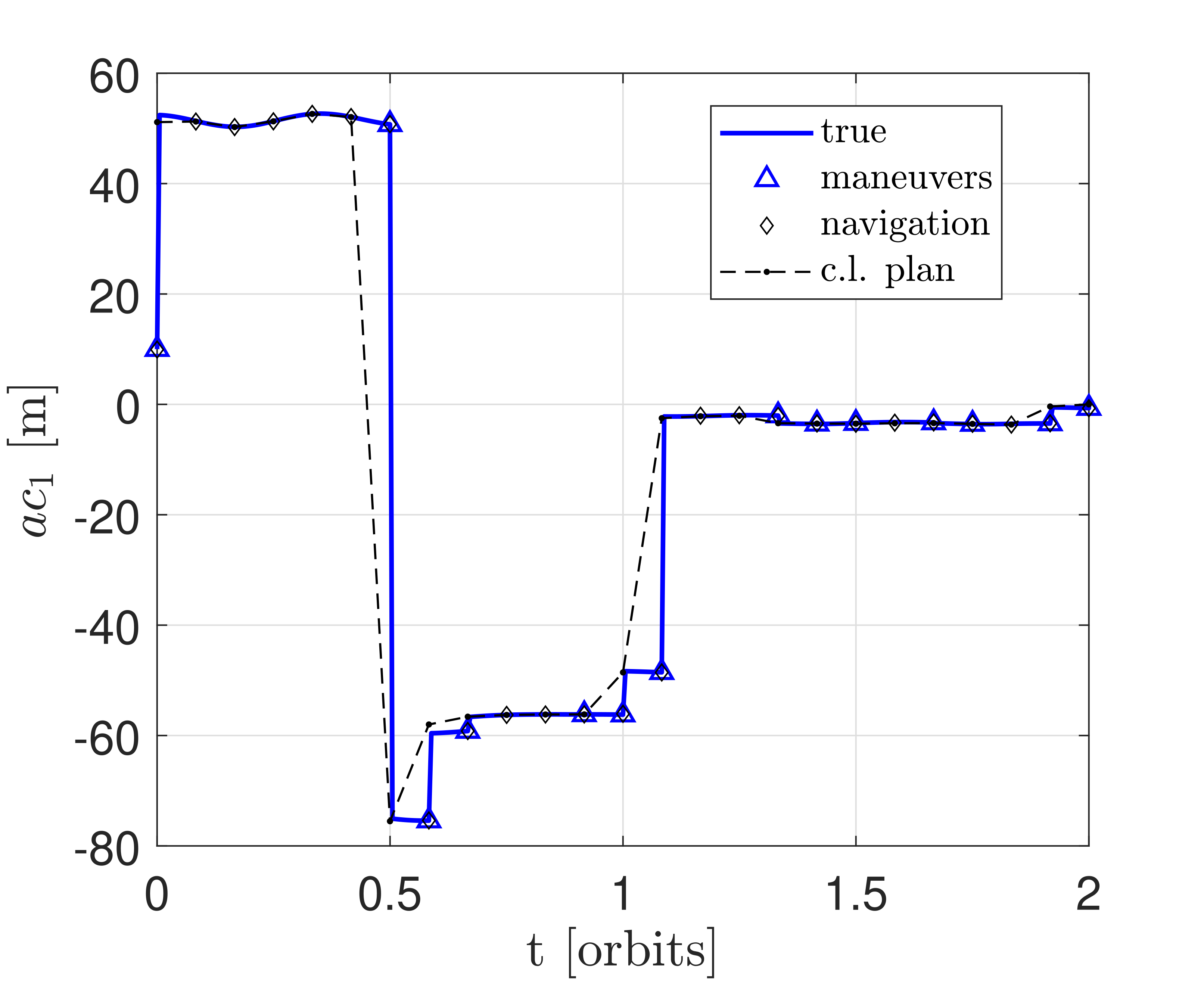}} \\ 
    \subfloat[Evolution of $a \boldsymbol{c}_{34}$ and $a \boldsymbol{c}_{56}$ magnitudes \label{fig:ICresults_T2_3}]{\includegraphics[width=0.4\linewidth]{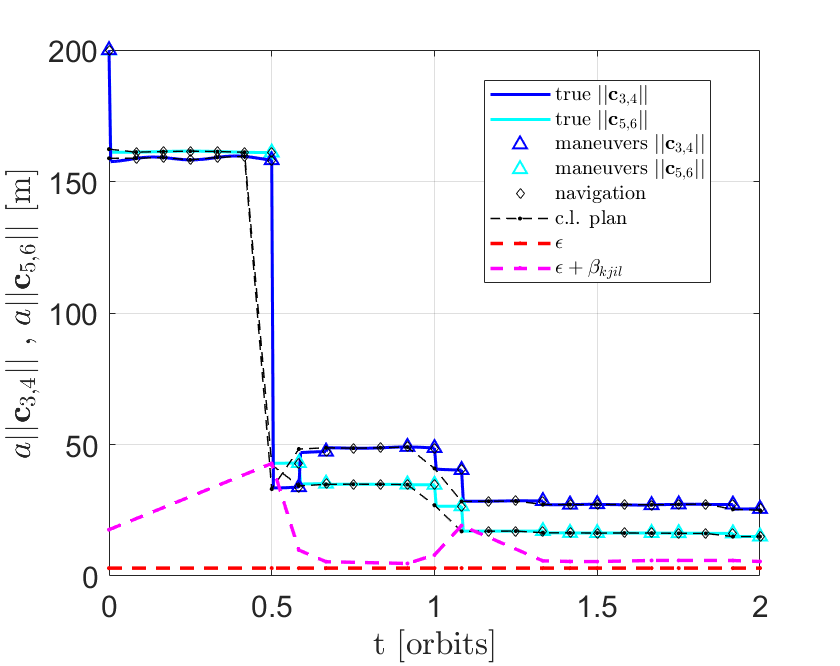}} \subfloat[Evolution of $\boldsymbol{c}_{34}$ and $\boldsymbol{c}_{56}$ angular separation \label{fig:ICresults_T2_4}]{\includegraphics[width=0.4\linewidth]{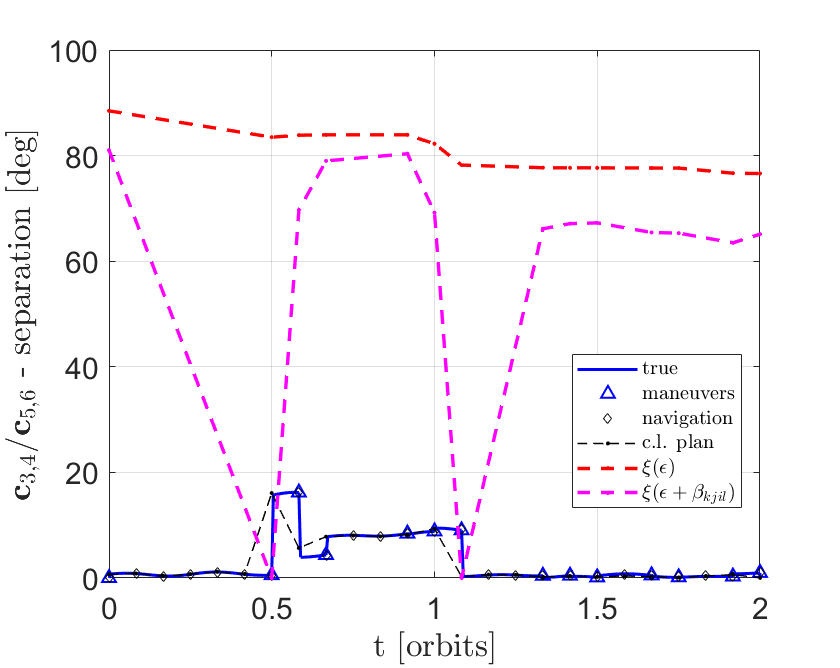}}
    \caption{ Separation, maneuver magnitudes and integration constants profiles over the two-orbit transfer} \label{fig:ICresults_T2}
\end{figure}

In Figure \ref{fig:cartresults_T2}, the resulting optimal transfer trajectory of the OSC around the DSC (at the origin) is presented both in $3$-D RTN space and projected in $2$-D RN plane.  The corresponding separation profiles are shown in Figures \ref{fig:RNsep_T2} and \ref{fig:3Dsep_T2}.
 Looking at Figure \ref{fig:cartresults_T2}, the controlled motion (blue bold line) is always instantaneously tangent to a safe uncontrolled trajectory (dashed lines), which is followed if the OSC loses control at that instant. Spacecraft separation along the uncontrolled trajectories is guaranteed to be $\geq 3$m for at least $2$ orbits, accounting for 3-$\sigma$ uncertainty (blue shaded ellipsoids). This is achieved in closed-loop,  using the model predictive control-like implementation described in Eq. \ref{OCP_IC}, \ref{OCP_IC_multi}, and in Algorithm \ref{ALG_C}. In particular, the black dot-dashed "closed-loop plan" line is the trajectory planned on-board. The blue "controlled true" line is obtained by executing the control plan through actuation and propagating the resulting motion through the ground-truth dynamics. The difference between these two lines lies in the difference between the on-board dynamics model used for optimization and the ground-truth, as well as in the errors introduced by navigation and actuation.  Similarly, the blue dashed "uncontrolled" trajectories are obtained integrating the ground-truth dynamics after each possible contingency. Whereas, the black dashed ones are predicted on-board and are superimposed to the 3-$\sigma$ ellipsoids modeling the uncertainty distributions used to quantify the bound $\beta_{kjil}$. The true uncontrolled trajectories are contained within the 3-$\sigma$ ellipsoids implying that the  bounds rightly capture and compensate for the difference between on-board belief and reality.  This validates the use of the models presented in Eq. \ref{Eq_bound_d}, \ref{Eq_dist_eff}, and \ref{Eq_bound}.
Note how larger uncertainty manifests after application of  larger maneuvers, due to the correlation between maneuver magnitude and actuation uncertainty presented in Table \ref{tab:scenario} on the right.
Comparing the control profile in Figure \ref{fig:ICresults_T2_1} with the uncertainty ellipsoids in Figure \ref{fig:cartresults_T2}, after the largest planned maneuver (second one),  the generated uncontrolled trajectory is characterized by very large uncertainty that has to be kept out of the origin in RN plane to have PS guarantees. This causes the algorithm to allocate a large bound ($\beta_{kjil}$) on the minimum separation to be kept on  the generate uncontrolled trajectory. The varying size of  these quantified bound in relation to the size of the applied maneuvers can be appreciated by comparing the purple ($\epsilon + \beta_{kjil}$) dashed lines in  Figures \ref{fig:RNsep_T2} and \ref{fig:ICresults_T2_3}, with the control profile. Looking at Figure \ref{fig:RNsep_T2}, the black line is $\mathcal{S}^{*RN}_{kjil}$, which assumes integrable dynamics after a contingency, the blue line accounts instead for ground-truth dynamics. $\mathcal{S}^{*RN}_{kjil}$ is enforced being $\geq \epsilon + \beta_{kjil}$, while the bound $\beta_{kjil}$ rightly capture the effects of non-integrable dynamics and uncertainty. This provides guarantee of PS (blue line $> \epsilon$) and validates the bound model presented in Eq. \ref{Eq_bound}.
The total control cost of the transfer is $0.27$m/s. The cost of the TPBVP is $0.253$m/s. A difference of $10\%$, considering the included constraints and uncertainties, represents a very fuel efficient transfer, well within the $0.3$m/s allocated in the mission delta-V budget \cite{VISORS}.

Finally, to compare with the closed-form analysis presented in Section \ref{Sec_Cl}, Figure \ref{fig:ICresults_T2} shows the evolution over the transfer of $a c_{1}$, the magnitude of the vectors $a \boldsymbol{c}_{34}$/$a \boldsymbol{c}_{56}$, and their angular separation. The magnitude of $\boldsymbol{c}_{34}$ and $\boldsymbol{c}_{56}$ sizes the ellipses drawn by uncontrolled trajectories in RN plane (see Figure \ref{Fig_ROE_Cart}, right), their angular separation relates to PS through Eq. \ref{dediSEPecc}, whereas, $c_1$ represents an  offset in radial direction of these ellipses (see Eq. \ref{CARTtoIC}). The trajectory computed using direct optimization achieves PS in RN plane either by having these ellipses centered around the origin with size $\geq \epsilon + \beta_{kjil}$, or, by decentering them in radial direction through a sufficiently large $|c_1|$ if their size is $< \epsilon + \beta_{kjil}$ (see Figure \ref{fig:cartresults_T2_2}). This analysis is validated in Figure \ref{fig:ICresults_T2}. At the second maneuver, the reduction of the $\boldsymbol{c}_{34}$ and $\boldsymbol{c}_{56}$ norms requires $c_1 \approx -80$m to achieve PS in RN plane. Comparing Figures \ref{fig:ICresults_T2_3} and \ref{fig:ICresults_T2_4}, when the norm of $\boldsymbol{c}_{34}$ or $\boldsymbol{c}_{56}$ becomes $< \epsilon + \beta_{kjil}$, $\xi_{kjil}$ (Eq. \ref{dediSEPecc}) shrinks so much that safe separation is not possible by centering the RN ellipse around the origin. Therefore, PS is achieved by properly decentering the ellipse in radial direction. The notable aspect is that  this clear geometrical insight and deep connection with closed-form analysis is obtained using Algorithm \ref{ALG_C} and the underlying convex optimization. This opens the way to much greater generalization and applicability than the specialized closed-form solutions.

\subsubsection{Five-orbit Transfer}
\begin{figure}[ht!]
    \centering
    \subfloat[$3$-D RTN space \label{fig:cartresults_T5_cl_1}]{\includegraphics[width=0.58\linewidth,trim=2 2 4 4,clip]{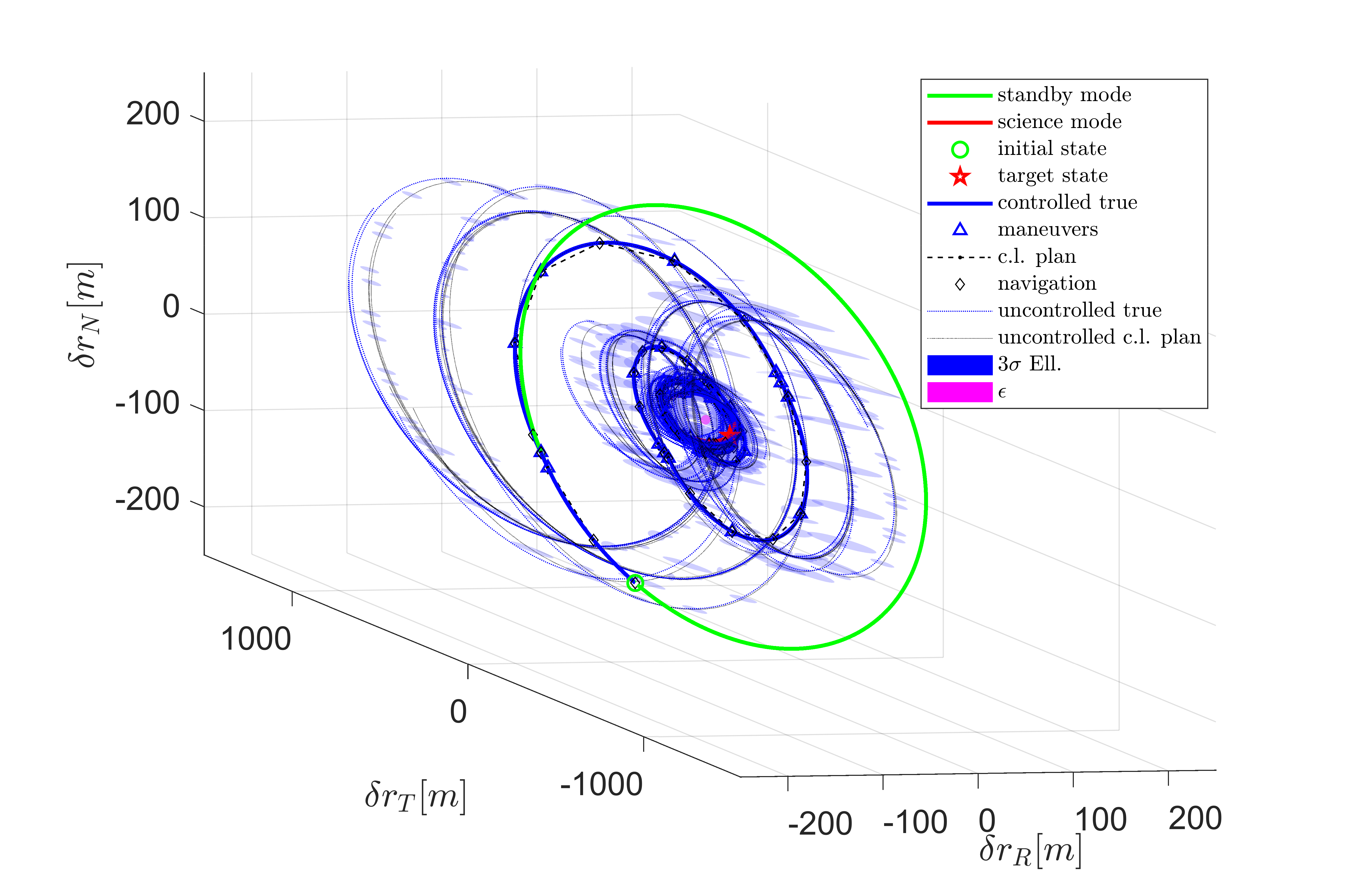}}
    \subfloat[$2$-D RN projection \label{fig:cartresults_T5_cl_2}]{\includegraphics[width=0.38\linewidth,trim=2 2 4 4,clip]{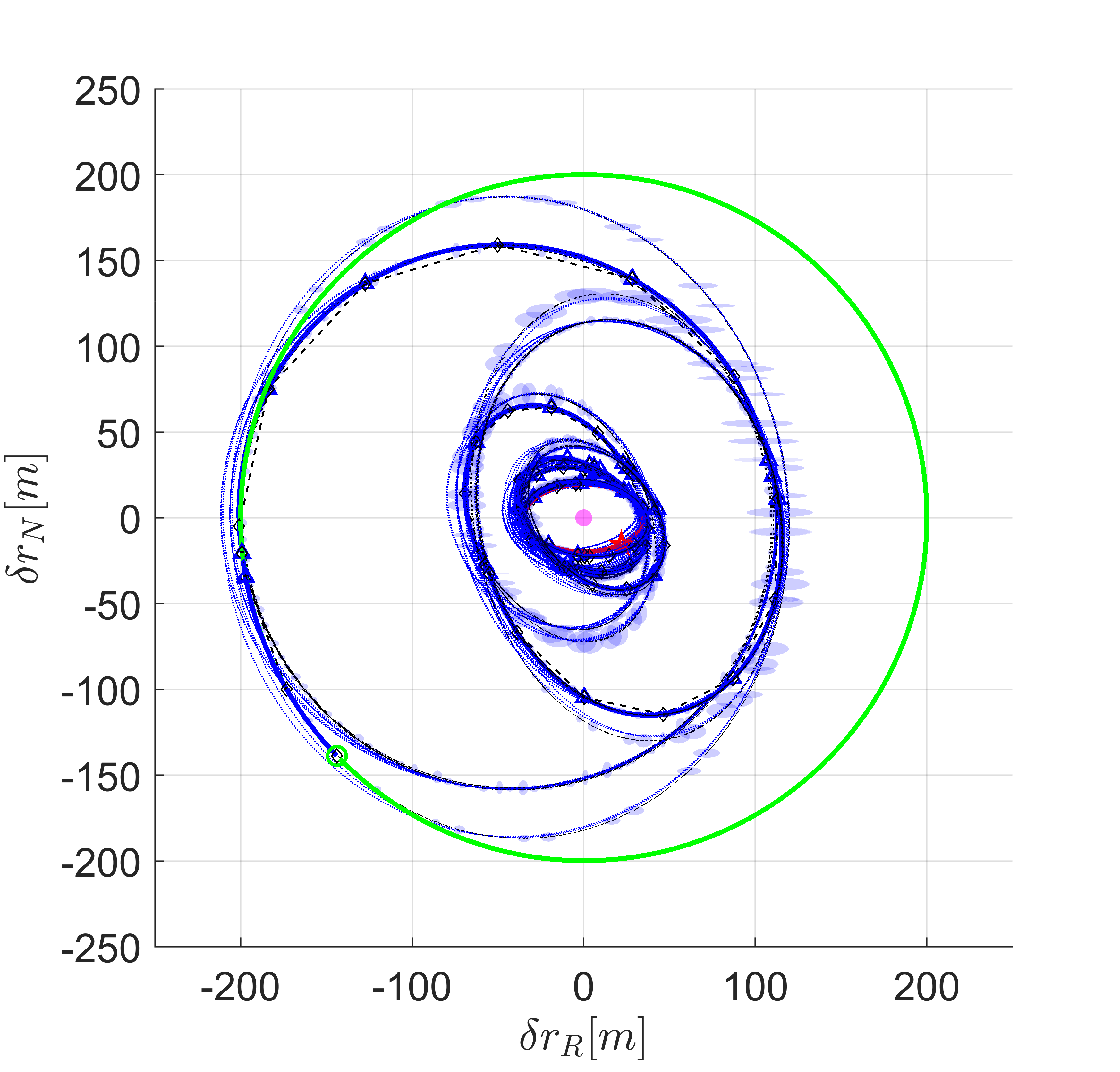}}
    \caption{Five-orbit passively-safe transfer trajectory} \label{fig:cartresults_T5_cl}
\end{figure}

The  five-orbit transfer reflects more closely what it is envisioned for the real mission. Differently from the  two-orbit one, where the TPBVP solution was used as initial warm-start, now an initial open-loop design is made to achieve RN PS throughout the transfer using closed-form analysis under simplifying assumptions on the dynamics and neglecting system uncertainties \cite{VISORS, guffanti_phd_2022}. 
\begin{table}[ht!]
\centering
\caption{Comparison between approaches, Monte Carlo analysis over 25 transfers}
\resizebox{0.98\textwidth}{!}{
\begin{tabu}{c c c c c c c c c c c c c c c}
\hline
\hline
\multicolumn{2}{c}{} & \multicolumn{6}{c}{Terminal state errors} & \multicolumn{1}{c}{} & \multicolumn{6}{c}{Other metrics} \\
\cline{2-8} \cline{10-15}
 \multicolumn{1}{c}{} & & $\epsilon_{ac_1}$ & $\epsilon_{ac_2}$ & $\epsilon_{ac_3}$ & $\epsilon_{ac_4}$ & $\epsilon_{ac_5}$ & $\epsilon_{ac_6}$ & & $\mathcal{J}$ & $min(\mathcal{S}^{*RN})$ & & time\footnotemark[3] & \# calls & PS Guarantees \\
   & & (m) & (m) & (m) & (m) & (m) & (m) & & (m/s) & (m) & & (s) & (-) & \\ \cline{1-15}
  \multirow{2}{*}{\begin{tabular}{@{}c@{}} Open-loop \\ design \cite{VISORS, guffanti_phd_2022} \end{tabular}} & mean & 4.40 & 164.89 & 1.63 & 3.86 & 1.33 & 5.44 &  & 0.324 & 15.91 & & 0.71 & 1 & \multirow{2}{*}{\begin{tabular}{@{}c@{}} Enforced \\ in open-loop \end{tabular}} \\
  & std & 3.01 & 110.25 & 0.88 & 2.99 & 0.73 & 3.91 &  & 0 & 3.42 & & 0.24 & - & \\
  \cline{2-8} \cline{10-15}
  \multirow{2}{*}{\begin{tabular}{@{}c@{}} Algorithm \ref{ALG_C} + \\ Tracking (Sec. \ref{Sec_Alg_Dec} \cite{guffanti_phd_2022}) \end{tabular}} & mean & 0.07 & 1.84 & 0.13 & 0.21 & 0.16 & 0.16 &  & 0.264 & 16.69 & \begin{tabular}{@{}c@{}} Alg. \ref{ALG_C}: \\ Track. : \end{tabular} & \begin{tabular}{@{}c@{}} 3.48 \\ 1.48 \end{tabular} & \begin{tabular}{@{}c@{}} 5 \\ 75 \end{tabular} & \multirow{2}{*}{\begin{tabular}{@{}c@{}} Enforced \\ in closed-loop \end{tabular}} \\
  & std & 0.06 & 1.88 & 0.17 & 0.12 & 0.20 & 0.11 &  & 0.004 & 0.37 & \begin{tabular}{@{}c@{}} Alg. \ref{ALG_C}: \\ Track. : \end{tabular} & \begin{tabular}{@{}c@{}} 0.33 \\ 0.08 \end{tabular} & \begin{tabular}{@{}c@{}} 1 \\ - \end{tabular} & \\
  \cline{2-8} \cline{10-15}
  \multirow{2}{*}{\begin{tabular}{@{}c@{}} Closed-form \\ tracking \cite{gaias_imp_2015} \end{tabular}} & mean & 0.46 & 6.35 & 0.49 & 1.30 & 0.62 & 1.09 &  & 0.333 & 16.70 & & 0.0032 & 3 & \multirow{2}{*}{\begin{tabular}{@{}c@{}} Checked \\ a posteriori \end{tabular}} \\
  & std & 0.32 & 4.73 & 0.55 & 0.95 & 0.63 & 0.90 &  & 0.01 & 0.91 & & 0.0030 & - & \\
  \hline \hline
  \end{tabu}}
  \label{tab:visors_perf_T5}
\end{table}
Table \ref{tab:visors_perf_T5} presents the statistics of a Monte Carlo of 25 transfers realized by implementing the control profiles within the ground-truth dynamics accounting for actuation errors and  in the presence of the various uncertainty sources introduced in Table \ref{tab:scenario}-\ref{tab:dynamics}. Looking at the top, when implemented in the ground-truth, the open-loop control design does not reach the desired terminal state (errors in along-track > $100$m). Therefore, to retrieve closed-loop properties, tracking is needed. This is done in two different ways. The first, Table \ref{tab:visors_perf_T5} on the bottom, uses state-of-the-art closed-form impulsive control solutions \cite{gaias_imp_2015}. The second, Table \ref{tab:visors_perf_T5} on the center, feeds the open-loop plan as initial warm-start to Algorithm \ref{ALG_C} which refines it, and then tracks it using the logic presented in Section \ref{Sec_Alg_Dec}. These two approaches are compared in the following. The transfer trajectory obtained using Algorithm \ref{ALG_C} + tracker is shown in Figure \ref{fig:cartresults_T5_cl}. The controlled motion (blue bold line) is always instantaneously tangent to a safe uncontrolled trajectory (dashed lines), which is followed if the spacecraft loses control at that instant. Spacecraft separation along the uncontrolled trajectories is guaranteed to be $\geq 5$m for at least $2$ orbits, accounting for 3-$\sigma$ uncertainty (blue shaded ellipsoids). The level of uncertainty in the  five-orbit transfer is lower than in the  two-orbit one due to the smaller applied maneuvers, and consequently smaller actuation uncertainty.

The closed-form tracker \cite{gaias_imp_2015} computes a correction to the open-loop control profile every $1.5$ orbits to null the expected terminal state error. The correction is a maneuver sequence  composed of three along-track and two cross-track maneuvers. Looking at Table \ref{tab:visors_perf_T5} on the bottom, good performances in terms of average terminal state error and PS ($min(\mathcal{S}^{*RN})$) are retrieved, at a very small average computational cost ($<0.005$s, for computing the correction over the subsequent $1.5$ orbits). The benefits of this approach are the simplicity, explainability and computational efficiency. The main limit is that formal guarantees of PS are not provided. In particular, the closed-form tracker achieves PS by tracking closely the open-loop trajectory, designed to satisfy RN PS. The closed-loop trajectory  is verified passively-safe a posteriori, but no explicit constrains are enforced during the closed-loop plan computation. In order to employ this algorithmic solution on-board it must be complemented with an additional logic that checks after each maneuver application the PS margin and, if needed, computes an escape maneuver if the safety limit is violated. Moreover, since the closed-loop tracker looks for correction maneuvers on-top of the open-loop plan, it achieves a total control cost which is $\geq$ than the open-loop one.

As presented throughout this paper, Algorithm \ref{ALG_C} instead looks for hard guarantees of PS. Using the open-loop design as a warm-start, Algorithm \ref{ALG_C} looks for candidate maneuvers every $30$deg plus in the locations of the open-loop design ones, for a total of $n_M = 75$ maneuver candidates over $[t_0, t_f]$. As done for the  two-orbit transfer, PS is enforced after each maneuver application to guarantee safety both on the trajectory in between maneuvers and, in case of subsequent maneuver failure, afterwards up to $t_f+T$. The time discretization at which PS must be enforced is at $1$deg, so $1800$ samples on $[t_0, t_f]$ and $720$ samples on $(t_f, t_f+T]$. The proposed efficient $\mathcal{M}^*$ formulation  permits the reduction of the number of constraints required to enforce PS from $> 75 \times 720 = 54000$ to $75$. The benefit from a computational standpoint is analogous to the one presented for the  two-orbit transfer. Looking at Table \ref{tab:visors_perf_T5} in the middle, Algorithm \ref{ALG_C} achieves an average terminal state error $\sim$$80\%$ smaller, and a total delta-V $\sim$$20\%$ smaller, than the closed-form tracker. This is mainly due to the refinement of the initial warm-start over a larger optimization search space of maneuver candidates. In general, the benefits introduced by Algorithm \ref{ALG_C} are the hard safety guarantees, high control accuracy performances, and flexibility of optimizing over a larger design space leading to a more fuel-efficient solution. The main limits are the lower explainability with respect to the closed-form solutions and the increased average computational cost. Looking at the run times presented in Table \ref{tab:visors_perf_T5}, using Algorithm \ref{ALG_C} takes in average $3.5$s, whereas the tracker ( Section \ref{Sec_Alg_Dec} \cite{guffanti_phd_2022}) takes in average $1.5$s. Note that the run time depends on various SCP parameters and tolerances which can be tuned to trade with other performance parameters, such as the fuel optimality reached at convergence. The ratio of Algorithm \ref{ALG_C} and tracking calls over the transfer is dictated by the high-level logic presented at the end of Section \ref{Sec_Alg_Dec}. Specifically, over the $75$ decision points, Algorithm \ref{ALG_C} is used once at the beginning and then  another four times all of which are due to tracking unfeasibility. This shows that the main trigger in recomputing the "guidance" plan is the incapability to remain within the specified tracking bounds, and not having the predicted worst-case minimum RN separation lower than $\epsilon$, which would entail a much greater risk from a safety perspective. Given the computational constraints for on-board implementation on the CubeSats microprocessors and looking at the average computational cost of using Algorithm \ref{ALG_C} with respect to the tracking,  it would be ideal to use Algorithm \ref{ALG_C} just once at the beginning of the transfer and then track. In this sense,  parameters such as the specified size of the tracking bound ($\boldsymbol{\nu}_{ki}^{\mathcal{T}}$) can be properly tuned to prevent unnecessary re-plannings and reduce the overall computational expense.

\subsection{Spacecraft Swarm in Eccentric Orbits}
\begin{figure}[ht!]
    \centering
    \includegraphics[width=0.38\linewidth,trim=2 2 4 4,clip]{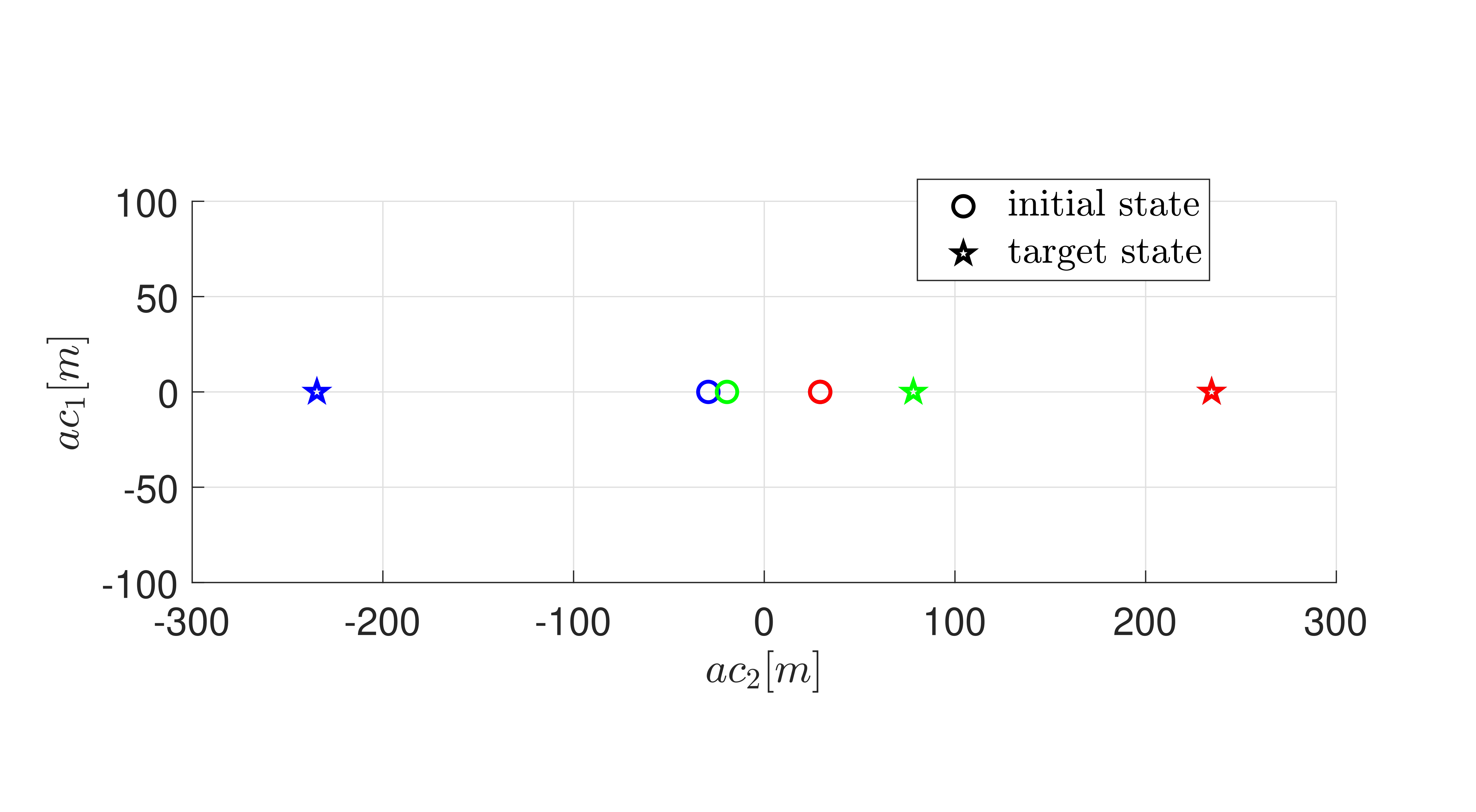}\hfill\includegraphics[width=0.26\linewidth,trim=2 2 4 4,clip]{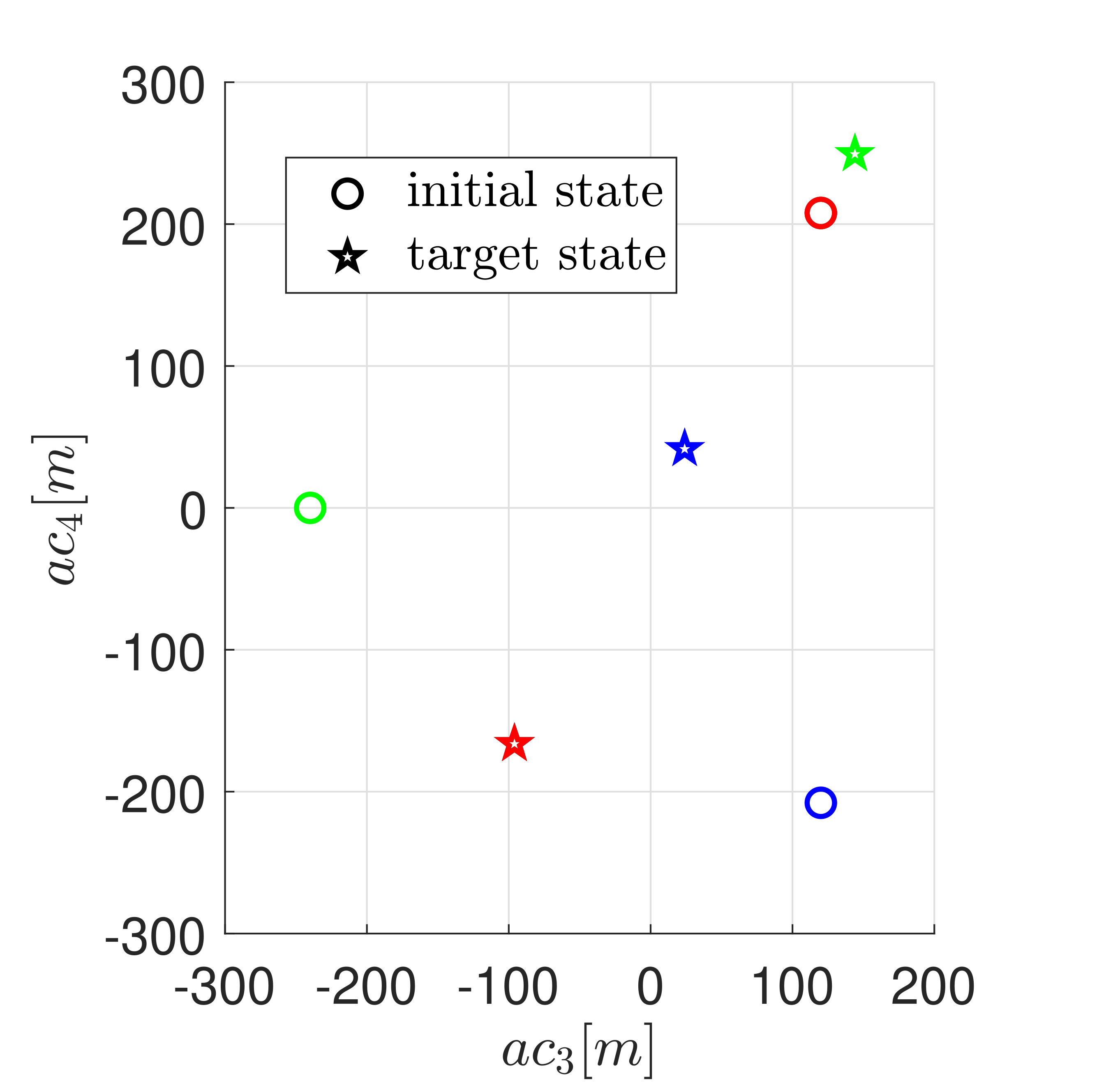}\hfill\includegraphics[width=0.26\linewidth,trim=2 2 4 4,clip]{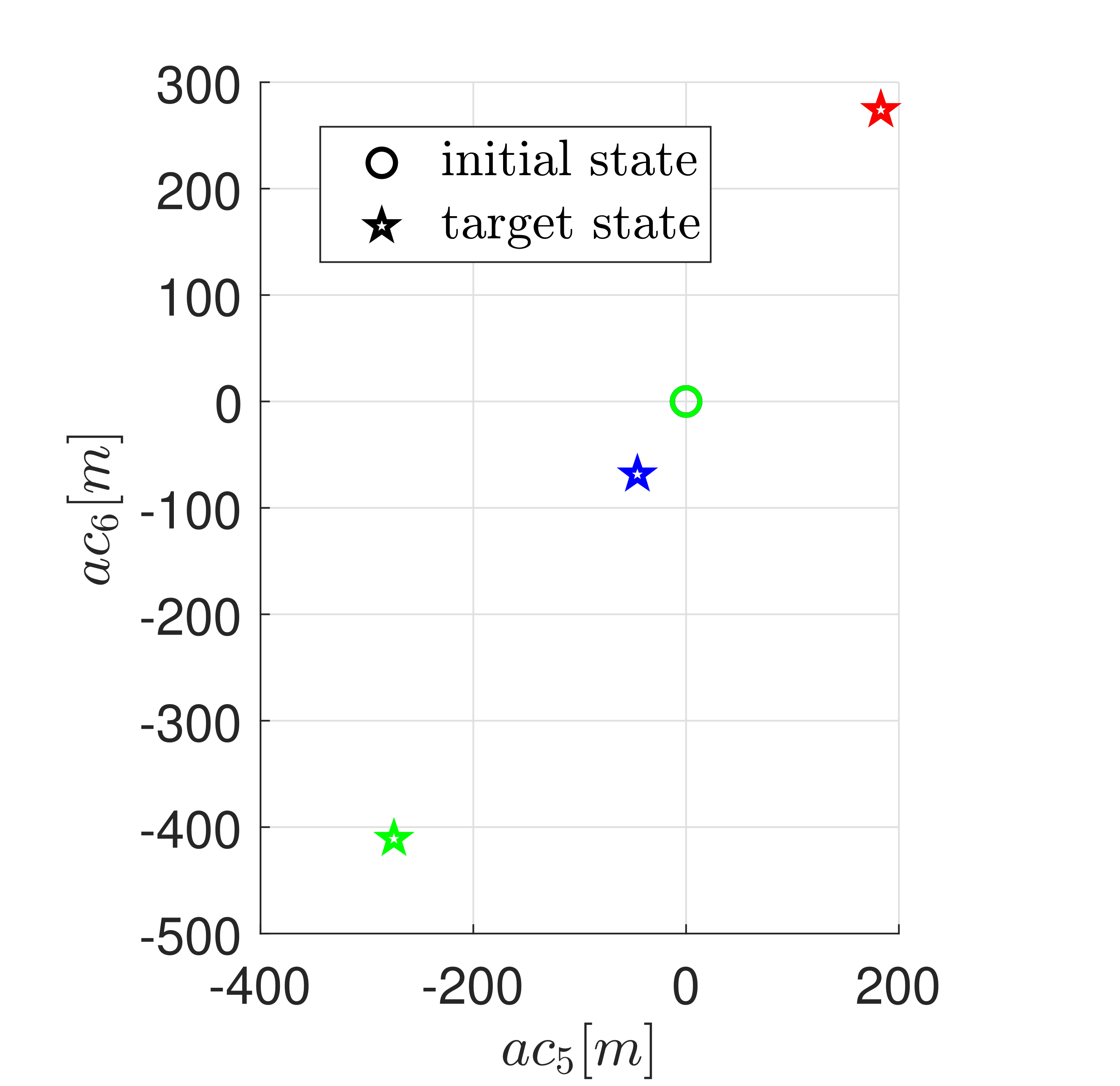}
     \caption{Swarm initial and target configurations in integration constants space}\label{Fig_HEO_0f_IC}
\end{figure}

\begin{figure}[ht!]
    \centering
    \subfloat[Initial RTN relative orbits \label{Fig_HEO_0f_Cart_1}]{\includegraphics[width=0.37\linewidth,trim=2 2 4 4,clip]{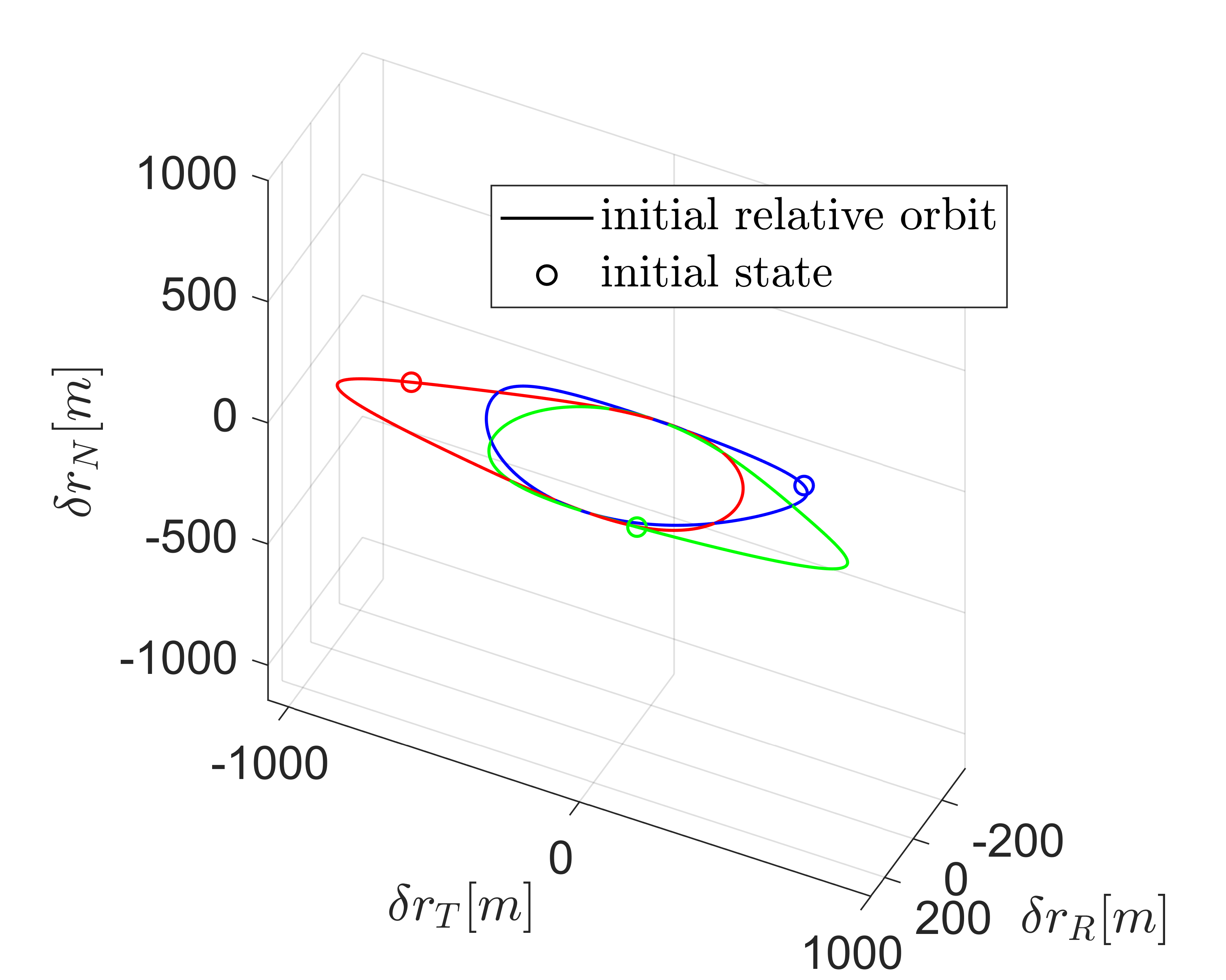}}
    \subfloat[RT (left) and RN (right) spacecraft separations on initial relative orbits \label{Fig_HEO_0f_Cart_2}]{\includegraphics[width=0.31\linewidth,trim=2 2 4 4,clip]{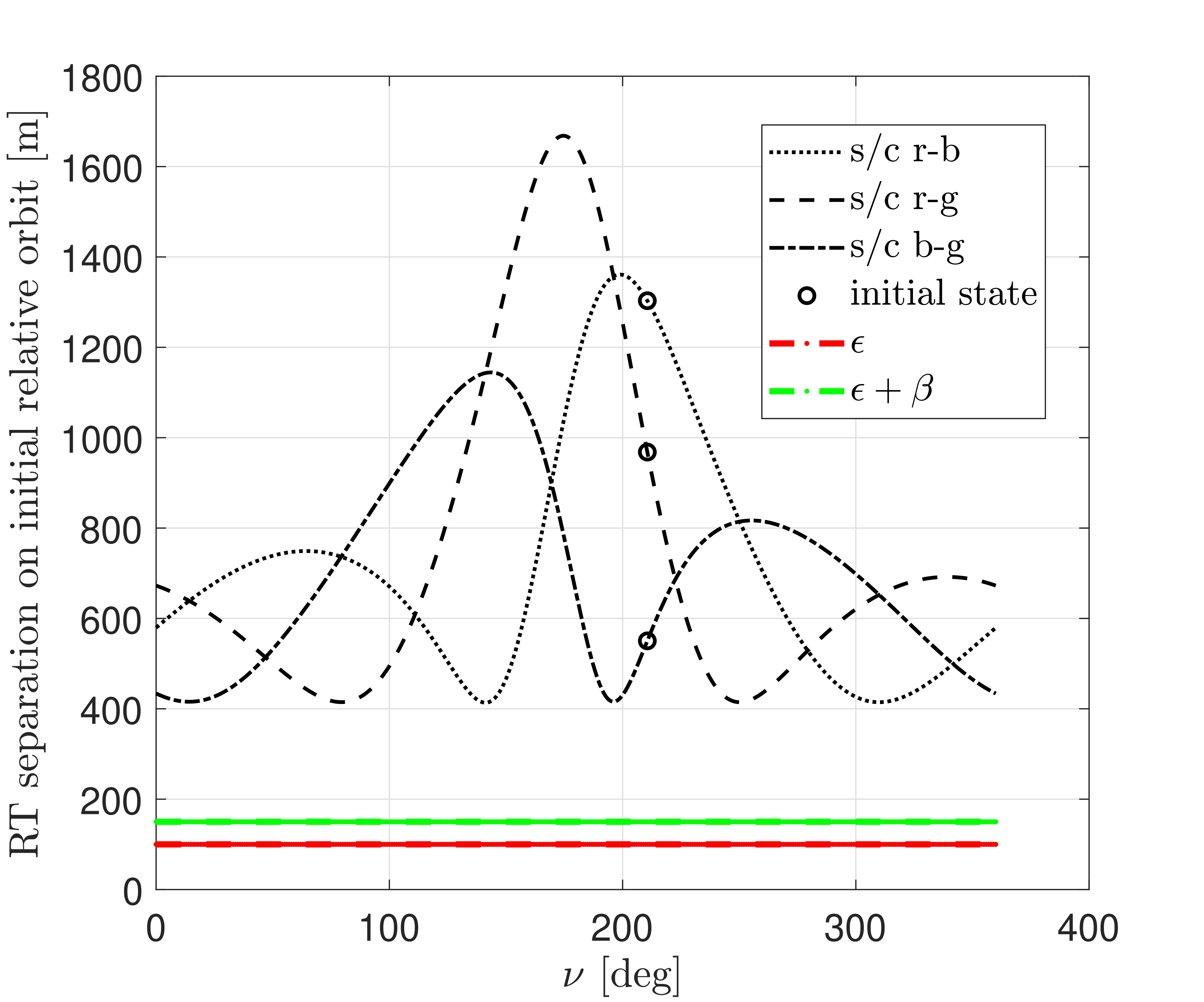} \includegraphics[width=0.31\linewidth,trim=2 2 4 4,clip]{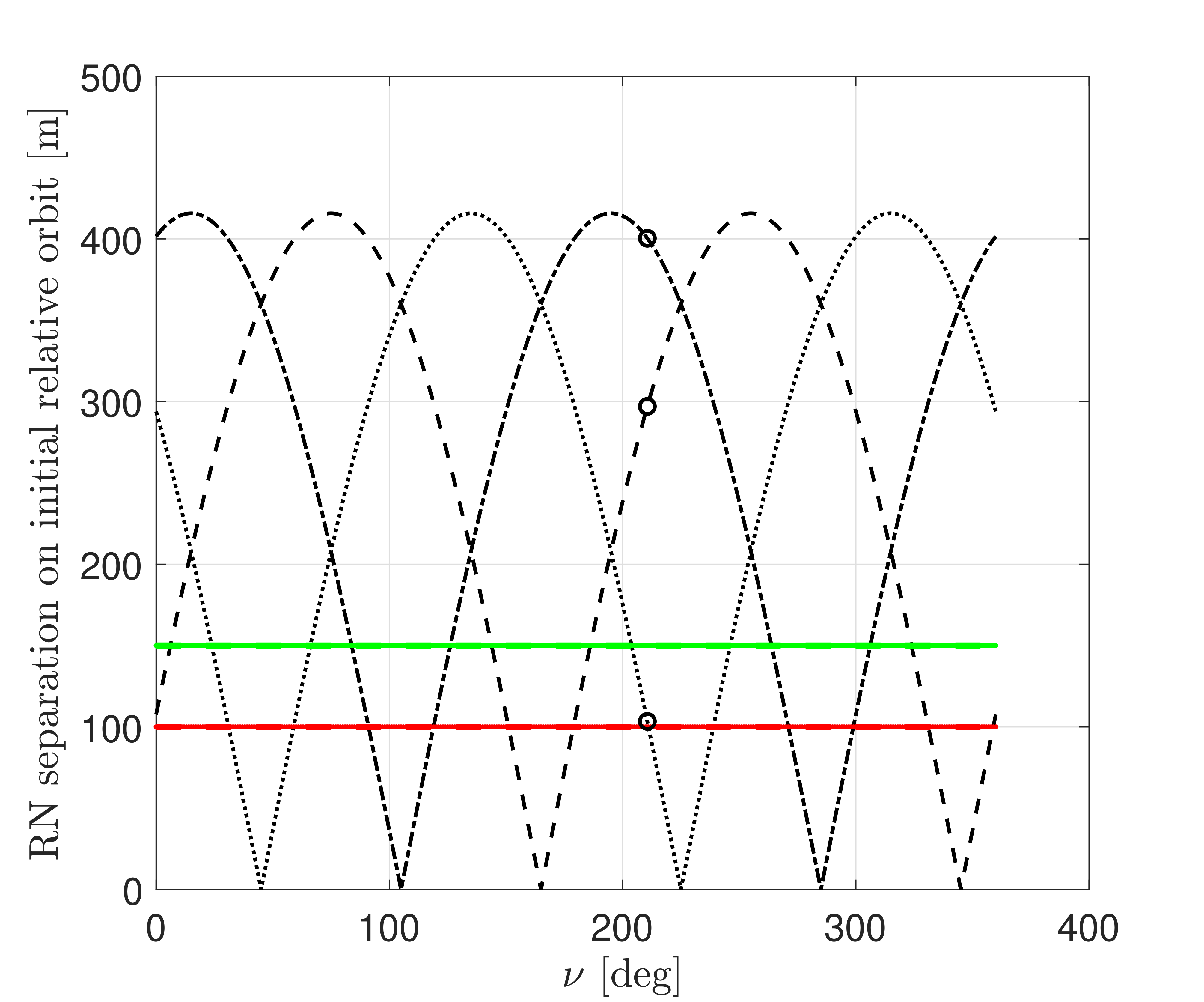}} \\
    \subfloat[Target RTN relative orbits \label{Fig_HEO_0f_Cart_3}]{\includegraphics[width=0.37\linewidth,trim=2 2 4 4,clip]{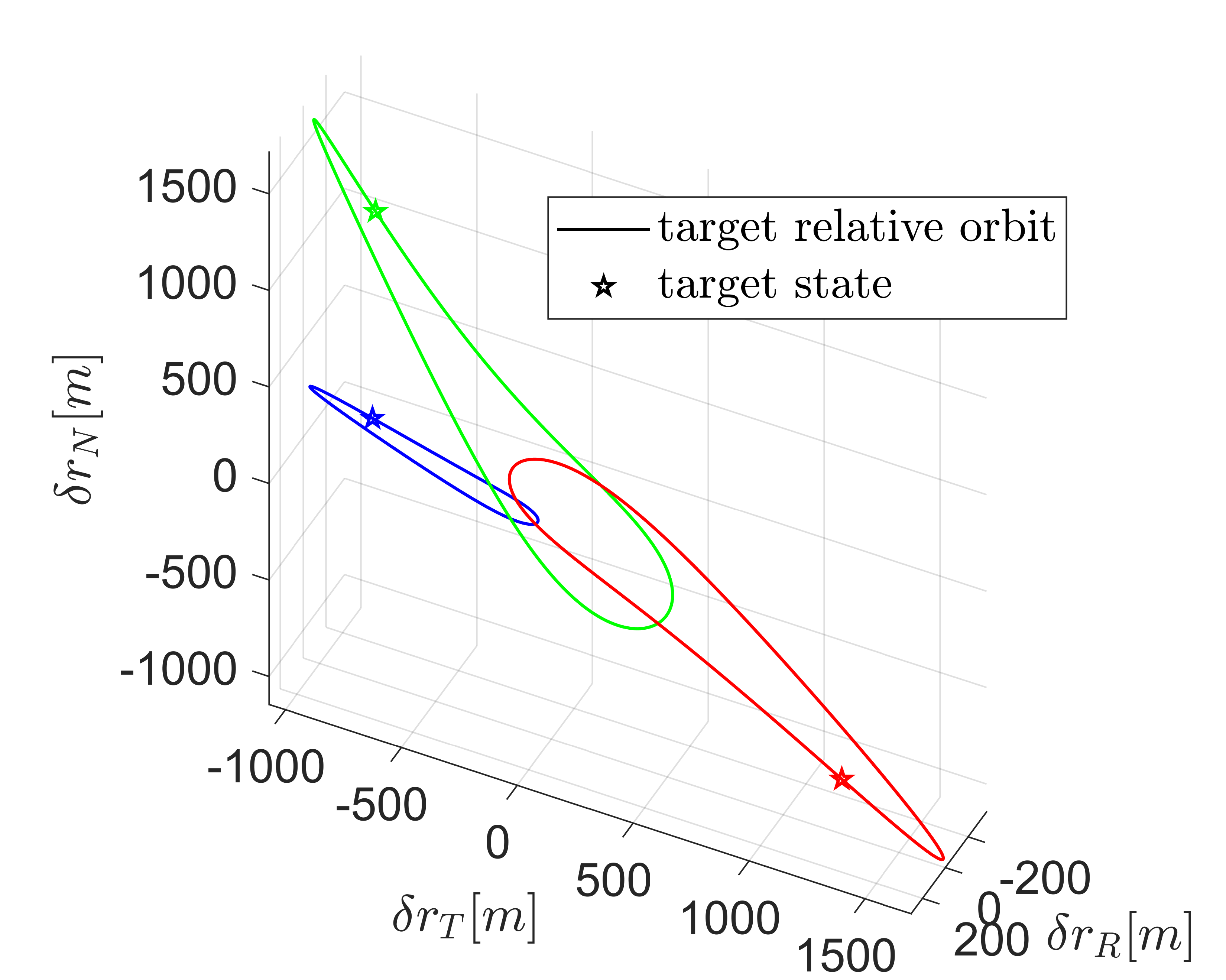}}
    \subfloat[RT (left) and RN (right) spacecraft separations on target relative orbits \label{Fig_HEO_0f_Cart_4}]{\includegraphics[width=0.31\linewidth,trim=2 2 4 4,clip]{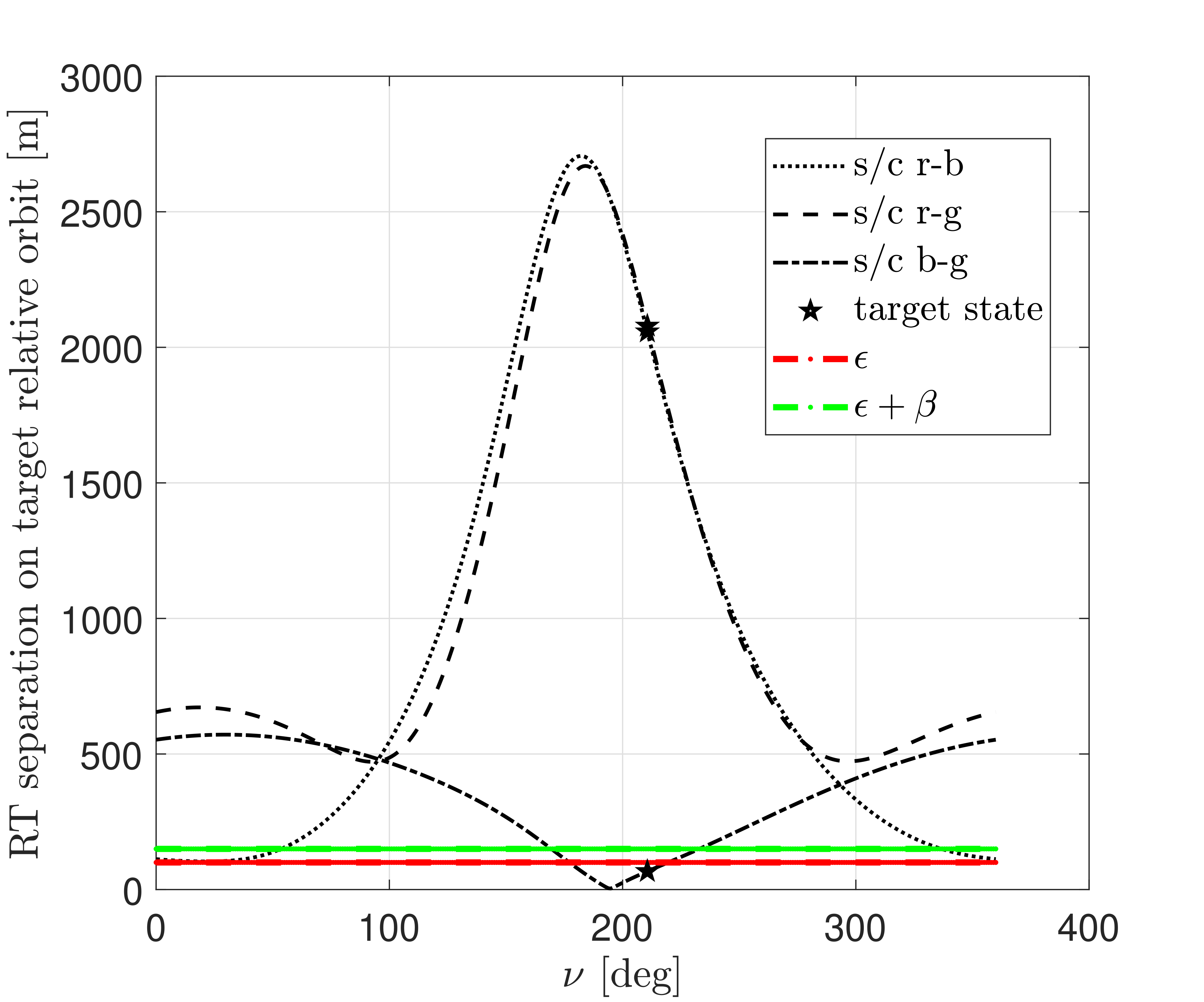} \includegraphics[width=0.31\linewidth,trim=2 2 4 4,clip]{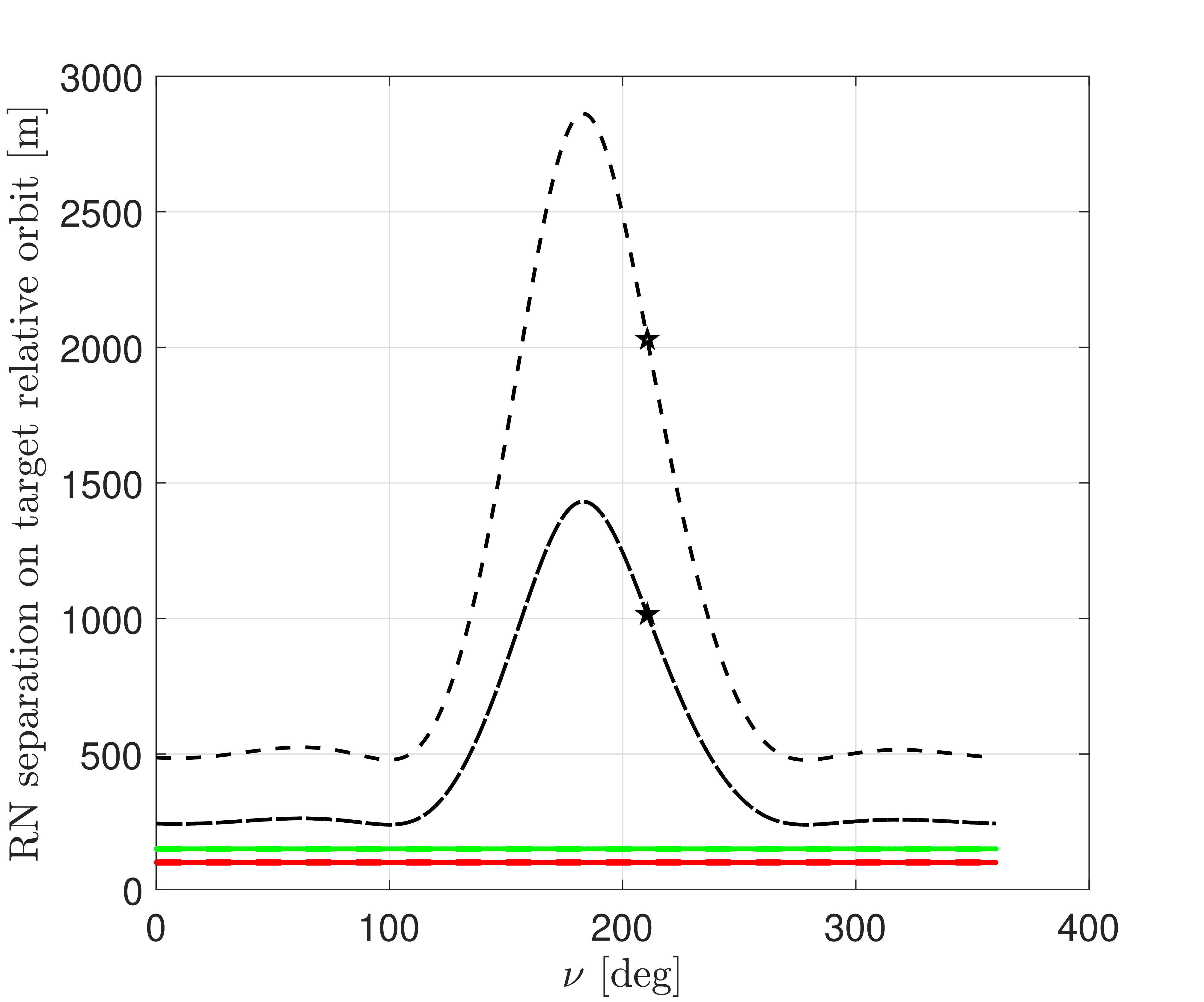}}
     \caption{Swarm initial and target relative orbits, and corresponding separation trends }\label{Fig_HEO_0f_Cart}
\end{figure}

In this section, Algorithm \ref{ALG_C} is assessed on the centralized open-loop trajectory design of a spacecraft swarm reconfiguration test case in highly eccentric orbit (HEO) around Earth. Three spacecraft have to reconfigure with respect to a virtual reference. The initial and target conditions are correspondingly a "high-density" (Eq. \ref{SwarmDesignRTecc2}) and an RN separated (Eq. \ref{SwarmDesignRN5}) swarm designs as defined in Section \ref{Sec_swarm_ecc}. These conditions are depicted graphically in Figure \ref{Fig_HEO_0f_IC} in  integration constants (IC) space, and in Figures \ref{Fig_HEO_0f_Cart_1} and \ref{Fig_HEO_0f_Cart_3} in Cartesian RTN space, and further characterized in Table \ref{tab:HEO_testcase}. Figures \ref{Fig_HEO_0f_Cart_2} and \ref{Fig_HEO_0f_Cart_4} show the spacecraft separations over the initial and target relative orbits. It is clear how the initial swarm configuration is safely separated in  radial/tangential (RT) plane but not in  radial/normal (RN), and the target one, viceversa, is safely separated in RN plane but not in RT.
\begin{table}[ht!] 
\centering
\caption{Highly eccentric orbit swarm test case }\label{tab:HEO_testcase}
\resizebox{\linewidth}{!}{
\begin{tabular}{ c c c c c c c c c c c c }
 \hline 
 \hline
  \multicolumn{2}{c}{Reference orbit} & & \multicolumn{3}{c}{Relative orbit modes} & & \multicolumn{2}{c}{Transfer \& safety} & & \multicolumn{2}{c}{Spacecraft \& features} \\
  \cline{1-2} \cline{4-6} \cline{8-9} \cline{11-12}
   Type & Virtual & & & Initial ($t_0$) & Target ($t_f$) & & $t_f - t_0$ (orbits) & 2 & & Mass (kg) & 80 \\
   \cline{5-6}
   $a$ (km) & 24641 & & $a c_1$ (m) & [0,0,0] & [0,0,0] & & PS type & $3$-D & & Cross Sec. Area (m$^2$) & 0.4 \\
   $e$ (-) & 0.716 & & $a c_2$ (m) & [29.3,-29.3,-19.6] & [234.6,-234.6,78.2] & & $T$ (orbits) & 1 & & $C_r$ (-) & 1.88 \\
   $i$ (deg) & 7 & & $a c_3$ (m) & [120.0,120.0,-240.0] & [-96.0,24.0,144.0] & & $\epsilon$ (m) & 100 & & Thrust (mN) & 50 \\
   $\Omega$ (deg) & -10 & & $a c_4$ (m) & [207.9,-207.9,0] & [-166.3,41.6,249.4] & & $\beta$ (m) & 50 & & $\zeta_i$ (m/s) & $\frac{\textrm{Thrust}}{\textrm{Mass}} \Delta t_i$ \\
   $\omega$ (deg) & 135 & & $a c_5$ (m) & [0,0,0] & [183.1,-45.8,-274.6] & & $n$ (-) & 24 & & &\\
   $\nu(t_0)$ (deg) & 210.6 & & $a c_6$ (m) & [0,0,0] & [273.9,-68.5,-410.9] & & $m$ (-) & 12 & & &\\
   Period (h) & 10.7 & & $a c_{34, min}$ (m) & 240 & 240 & & $N_a$ (-) & 3 & & & \\
   Epoch & 1-1-2021 & & $a c_{56, min}$ (m) & - & 411.8 & & $\mathcal{J}$ (m/s) & $\sum_{ki} || \tilde{\boldsymbol{u}}_{k}(t_i) ||_{2}$ & & & \\
   GPS time & 00$^h$:00$^m$:00$^s$ & & $\theta^{hd}_{34}$ (deg) & 120 & - & & & & & & \\
    & & & $\theta^{sep}_{34}$ (deg) & - & 60 & & & & & & \\
    & & & $\theta^{sep}_{56}$ (deg) & - & 236.2 & & & & & & \\
  \hline
  \hline
\end{tabular}}
\end{table}

A reconfiguration between the two designs while enforcing PS in $3$-D is showcased, demonstrating a capability not published in literature. The reconfiguration horizon is fixed to two reference orbit periods, and PS has to be guaranteed in $3$-D for at least one orbit ($T$), of at least $100$m ($\epsilon$). This open-loop trajectory design is treated as deterministic, therefore no uncertainty sources are included. The spacecraft relative motion dynamics is modeled as Keplerian perturbed by Earth oblateness ($J_2$) and solar radiation pressure (SRP) \cite{guffanti_linear_2019}.  The scenario is considered deterministic, therefore the only contribution to be included in the quantification of the bound $\beta_{kjil}$ in Eq. \ref{Eq_bound} is the integrated effect of $J_2$+SRP over the PS horizon. Moreover, while  in the presence of uncertainties, the quantification of  this bound must be corrected in closed-loop within Algorithm \ref{ALG_C}. When just deterministic effects are included, the worst case expected effect of non-integrable dynamics over the PS horizon can be quantified a priori. This is done for this test case, and a constant PS bound ($\beta = 50$m) is considered. The spacecraft are modeled as continuously controlled by low-thrust actuation. To mimic the realizability of the control profile by the considered actuation a cap $\zeta$ on the control input $L_2$-norm is imposed. In particular, given a constant thrust and the spacecraft mass, using second Newton's law, $\zeta_i = (\textrm{Thrust}/\textrm{Mass}) \Delta t_i$, where $\Delta t_i$ is the maximum interval over which a computed delta-V must be realized by the low-thrust. This implies a convex constraint on the control input. The $L_2$-norm of the control input is minimized. The orbit is discretized at $30$deg, entailing a time-variant time step ($\Delta t_i$) given the spacecraft varying velocity along the eccentric orbits.

\subsubsection{Results}
\begin{figure}[ht!]
    \centering
    \subfloat[Enforcing  passive safety \label{Fig_HEO_Cart_PS_CA_1}]{\includegraphics[width=0.49\linewidth,trim=2 2 2 2,clip]{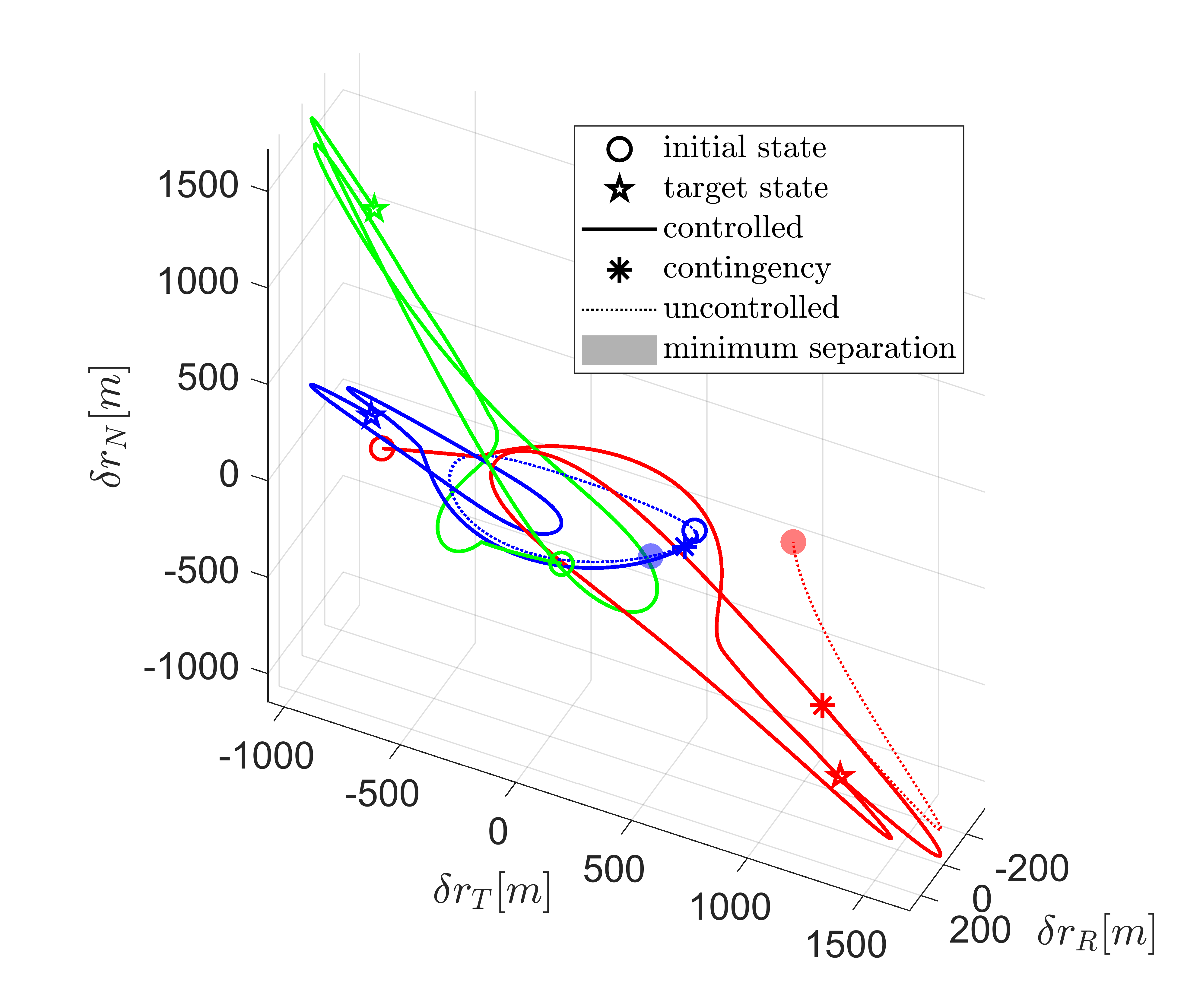}}
    \subfloat[Enforcing  collision avoidance \label{Fig_HEO_Cart_PS_CA_2}]{\includegraphics[width=0.49\linewidth,trim=2 2 2 2,clip]{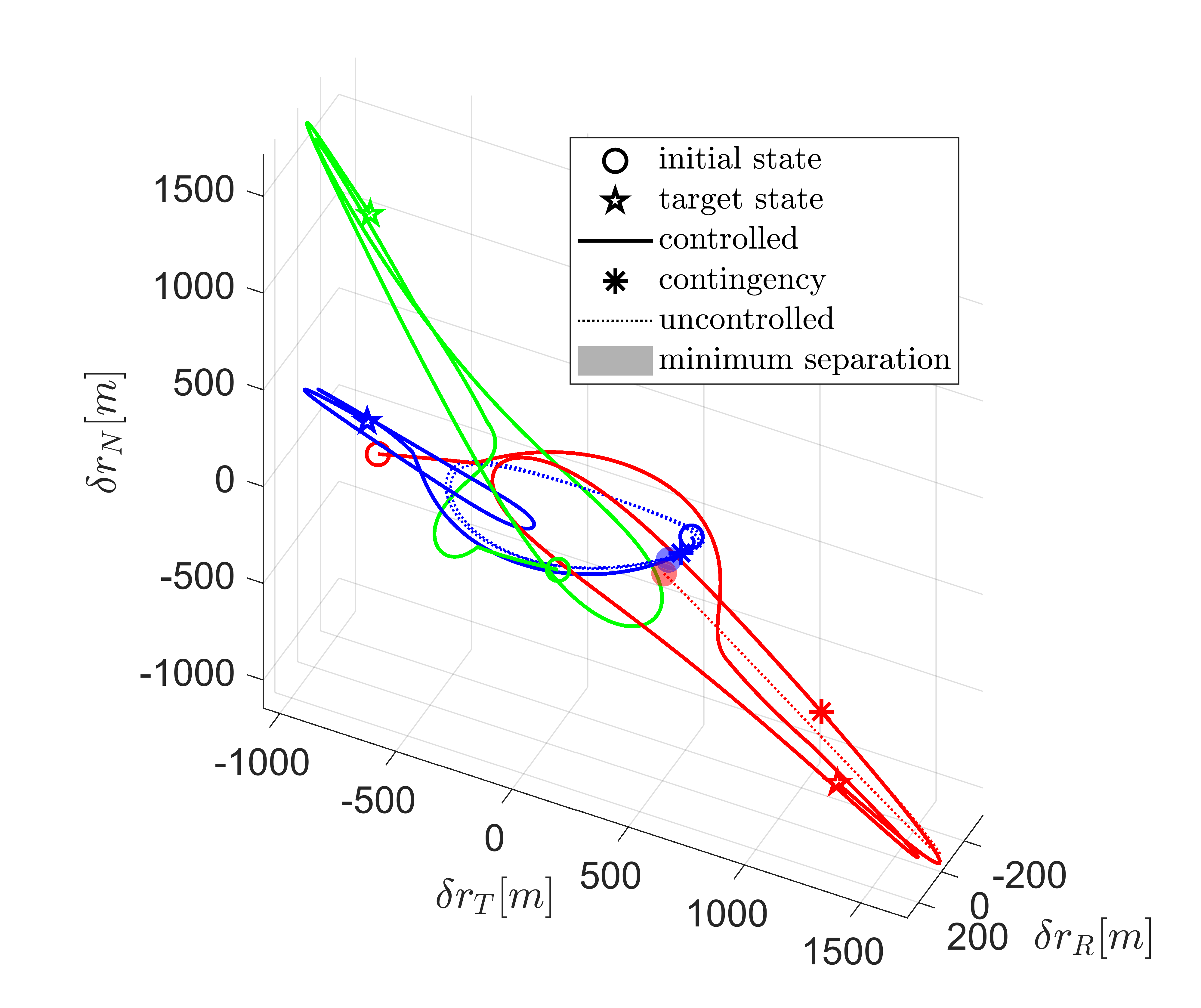}} \caption{ Optimal transfer trajectories.} \label{Fig_HEO_Cart_PS_CA_01}
\end{figure}

\begin{figure}[ht!]
    \centering
    \subfloat[Enforcing  passive safety \label{Fig_HEO_Cart_PS_CA_3}]{\includegraphics[width=0.45\linewidth,trim=2 2 2 2,clip]{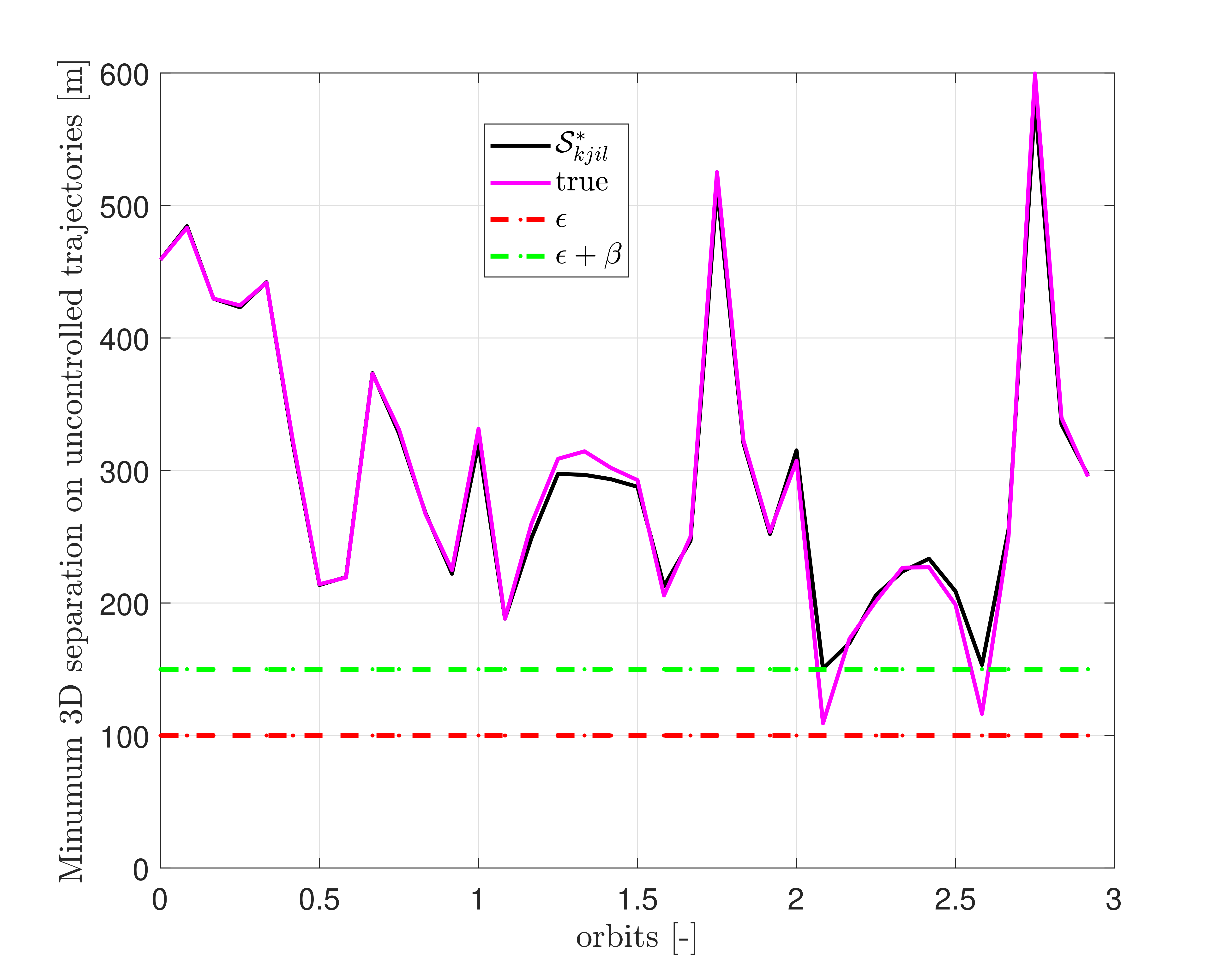}}
    \subfloat[Enforcing  collision avoidance \label{Fig_HEO_Cart_PS_CA_4}]{\includegraphics[width=0.45\linewidth,trim=2 2 2 2,clip]{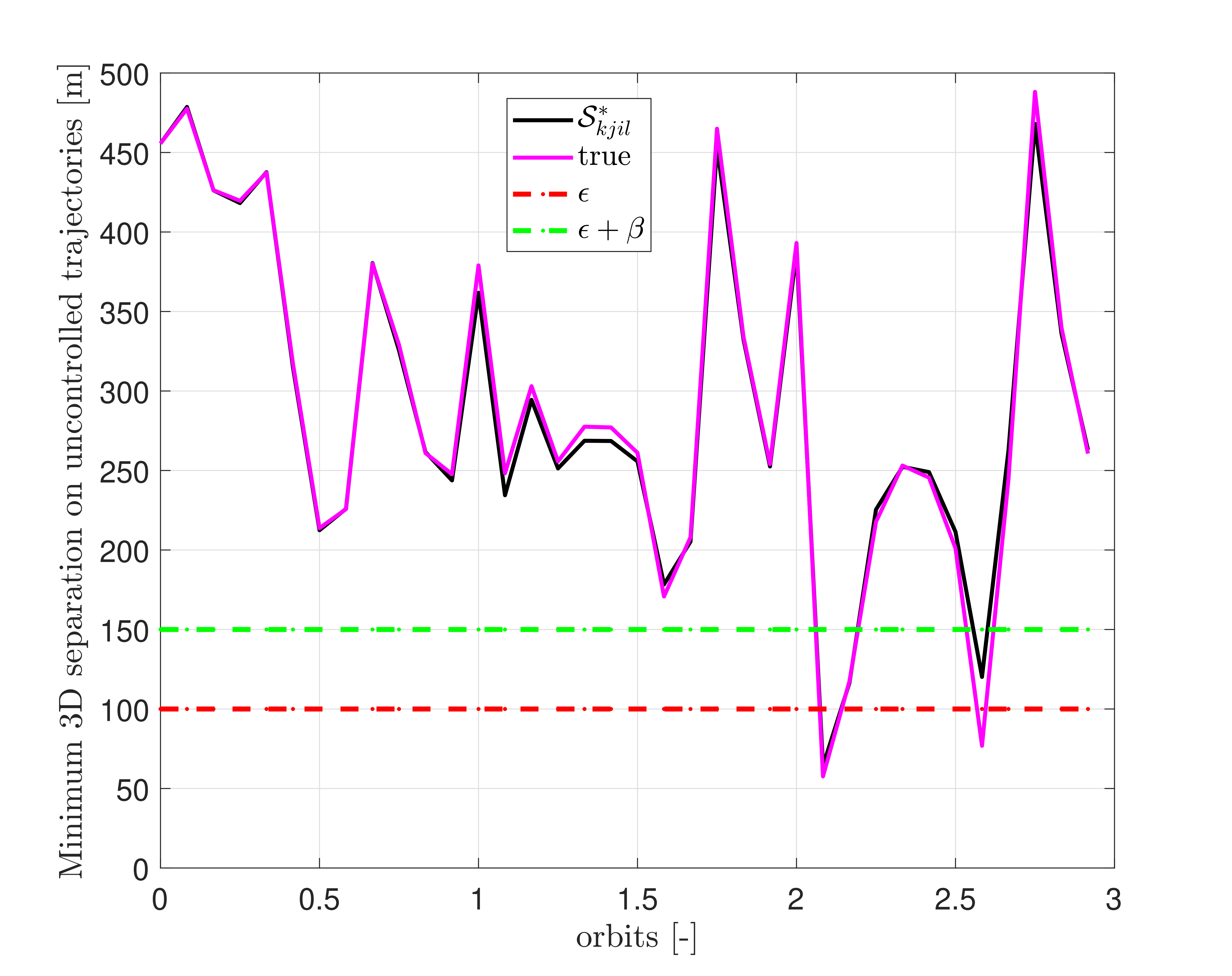}}
     \caption{ Swarm minimum separation accounting for all possible combinations of contingencies}\label{Fig_HEO_Cart_PS_CA_02}
\end{figure}

For this test case, a comparison between PS and standard collision avoidance (CA) is presented. In particular, when enforcing CA, the horizon is $T = 0$ by definition. All the subroutines in Algorithm \ref{ALG_C} related to the computation of $t^*_{kjil}$ disappear, and the constraint is simply linearized in the SCP as done in literature \cite{liu_scp_2014, morgan_2014, malyuta_scp_2022}. To warm-start the CA SCP solver the TPBVP solution \cite{ecos_13} is used. When enforcing PS instead, a safety horizon $T = 1$ orbit is considered, and the solution enforcing CA is used as first attempt of warm-start.  It is relevant to remark that the specific values of $\epsilon$ and $T$ do not affect the scope of the presented comparison. In particular, this comparison aims to show that, given arbitrarily selected values of $\epsilon$ and $T$, the instantaneous collision avoidance constraint does not provide guarantee of safe separation along the uncontrolled trajectories after a contingency, whereas, if feasible, the passive-safety constraint does. There may exist a specific combination of $\epsilon$ and $T$ for which the motion is safe on the uncontrolled trajectories even by enforcing instantaneous collision avoidance, but this does not constitute a guarantee. The objective of this paper is to prove enforcement of fault-tolerant motion safety guarantees, whenever feasible. The results of the comparison are reported in Figures \ref{Fig_HEO_Cart_PS_CA_01}-\ref{Fig_HEO_Cart_PS_CA_02} and in Table \ref{tab:HEO_performance}.  Figure \ref{Fig_HEO_Cart_PS_CA_01} analyzes a specific combination of contingencies. In particular, blue and red spacecraft lose control $0.075$ orbits (blue asterisk) and $1.825$ orbits (red asterisk) after the beginning of the transfer. Looking at Figure \ref{Fig_HEO_Cart_PS_CA_2}, if CA is enforced, the resulting uncontrolled trajectories (blue and red dashed line) reach a minimum separation of $57.6$m $<\epsilon$. Looking at Figure \ref{Fig_HEO_Cart_PS_CA_1}, if PS is enforced, the same contingencies result in uncontrolled trajectories that maintain safe separation. Figure \ref{Fig_HEO_Cart_PS_CA_02} enlarges the analysis to all the possible combinations of contingencies on $[t_0, t_f]$ and corresponding separations up to  $t_f + T = 3$ orbits. The black line is $\mathcal{S}^*_{kjil}$, which assumes integrable dynamics on the uncontrolled arcs, and the purple line accounts instead for $J_2$+SRP effects. In Figure \ref{Fig_HEO_Cart_PS_CA_3}, $\mathcal{S}^*_{kjil}$ is enforced being $\geq \epsilon + \beta$,  while the worst case effect of perturbations is lower than the bound $\beta$. This provides guarantee of PS (purple line $> \epsilon$), and validates the bound model presented in Eq. \ref{Eq_bound_d} and \ref{Eq_bound}. In Figure \ref{Fig_HEO_Cart_PS_CA_4}, both black and purple lines go lower than $\epsilon$ violating the PS constraint. 
This comparison in summarized in Table \ref{tab:HEO_performance} on the right. Both CA and PS guarantee a minimum separation $\geq \epsilon$ when no contingencies happen. Accounting for contingencies instead, PS guarantees the separation on the uncontrolled trajectories being greater than $\epsilon$, whereas CA does not. Looking at the control costs in Table \ref{tab:HEO_performance} on the left, the number reported are totals for the entire swarm. The control cost of the solution enforcing PS is $2$\% larger than the one enforcing CA. This is a fuel efficient solution considering the additional number of constraints included.
\begin{table}[ht!]
\centering
\caption{Comparison of  passive safety (PS) and standard  collision avoidance (CA)}\label{tab:HEO_performance}
\resizebox{0.8\textwidth}{!}{
\begin{tabular}{c c c c c c c c c c c}
 \hline
 \hline
  \multicolumn{1}{c}{} & \multicolumn{1}{c}{Control cost} & \multicolumn{1}{c}{} & \multicolumn{1}{c}{\# Constraints} & \multicolumn{1}{c}{} & \multicolumn{2}{c}{Running time\footnotemark[3]} & & $\epsilon$ & \multicolumn{2}{c}{$\min_{kjil} \mathcal{S}_{kjil}^*$}\\
  \cline{6-7} \cline{9-11}
  & & & & & Warm-start & & & & No contingencies & With contingencies \\
  & (m/s) & & (-) & & (-) & (s) & & (m) & (m) & (m) \\ \hline
  PS & 0.275 & & $n^2 \frac{N_a (N_a-1)}{2}$ & &  \begin{tabular}{@{}c@{}} CA \\ PS \end{tabular} & \begin{tabular}{@{}c@{}} 34.5 \\ 6.3 \end{tabular} & & 100 & 188.2 & 109.2 \\
  \cline{2-11}
  CA & 0.270 & & $n \frac{N_a (N_a-1)}{2}$ & & TPBVP & 2.8 & & 100 & 248.3 & \hspace{0.85cm} 57.6 \hspace{0.01cm} ($< \epsilon$)  \\ 
  \hline
  \hline
\end{tabular}}
\end{table}
Looking at the run times in Table \ref{tab:HEO_performance} on the center, enforcing PS is more computational expensive than enforcing CA. This larger computational expense is related to two factors. Firstly, the number of constraints to be enforced in each underlying SOCP, which is larger for PS ($n^2 N_a (N_a-1)/2$) than for CA ($n N_a (N_a-1)/2$). Secondly, the computations in Algorithm \ref{ALG_C} (lines 9-24) required for  enforcing PS given an initial warm-start that is not necessarily feasible (e.g., that satisfies CA but not PS). The enforcement trend is the one presented in Figure \ref{Fig_tstar_alg_2}. In particular,  the further the initial warm-start is from being feasible, the more inner-loop iterations are needed, the more solutions of SOCP instances \cite{ecos_13} including $n^2 N_a (N_a-1)/2$ constraints are required, the  higher overall computational cost is. Consequently, a bad initial warm-start has a substantial impact. Looking at Table \ref{tab:HEO_performance}, this is further proved by feeding-back as a warm-start a previously found passively-safe solution. Now, Algorithm \ref{ALG_C} is $\sim$$5$ times faster. This test case underlines the importance of properly warm-starting SCP-based algorithms to  reduce their computational expense, in particular for on-board execution.

\section{Conclusion}\label{Sec_VII}
 This paper has addressed the problem of enforcing motion safety guarantees even in the event of sudden loss of control capabilities by any agent within a multi-agent system. The paper shows that, to enforce these so-called passive safety guarantees, point-evaluated constraints on the trajectory following a contingency can be replaced with a single function of the integration constants of the equations of motion at contingency instants. Furthermore, effects of uncertainties on the trajectory can be compensated for through the method of variation of parameters. This permits reducing the number of required constraints by one polynomial degree in the number of discrete time samples, and allows computationally efficient enforcement within a multi-agent optimal control problem, solvable using direct methods. Using this novel methodology, experimental results show that a control profile enforcing passive safety in closed-loop in a multi-agent setting can be computed efficiently through an algorithm based on sequential convex programming and model predictive control. This allows the computation of passively-safe spacecraft reconfigurations in orbital scenarios outside of the assumption limits of the current literature, in both near-circular and eccentric orbits, in the presence of realistic system uncertainties stemming from sensing, actuation, and unmodeled system dynamics. In addition, the application of the proposed approach to the spacecraft relative motion dynamics permits the derivation of novel solutions for multi-spacecraft passive safety in eccentric orbits, with application to relative motion design and guidance. The methodology presented in this paper opens the way to efficient, safe and fault-tolerant real-time motion planning and control on-board, which is the core requirement for autonomous functions in future distributed space systems' missions.

\section{Appendix}
\subsection{Passive Safety Constraint Requirements}\label{App:PS}
In the following, formal proof of the number of constraints required to enforce PS reported in Table \ref{tab_PS} is developed. The proof is developed for the "impulsive" control type, it is left to the reader verifying that for $n_M \rightarrow n$, the number of constraints required in impulsive control tends to the ones required in continuous control. The number of active agents' combinations is $N_a (N_a-1)/2$, where $N_a$ is the number of active agents. $N_p$ is the number of permanently passive agents.\\
 \textit{Proposition.} Assume the interval $[t_0, t_f]$ is discretized in $n$ samples, and the interval $(t_f, t_f+T]$ in $m$ samples.  Assume impulsive control application is planned and possibly failed just in a subset $n_{M} \leq n$ of the $n$ instants, such as $n/n_{M} \in \mathbb{N}_{>0}$. In order to guarantee PS (as in  Definition 1 in Section \ref{Sec_III}) at each of the $n$ and $m$ instants in $[t_0, t_f+T]$, using the $\mathcal{M}$ formulation, the total number of collision avoidance constraints the multi-agent system has to enforce is: $[(n/n_{M})(\sum_{i=1}^{n_{M}}i^2)+(n_{M}+1)^2 m] N_a(N_a-1)/2 = \{(n/n_{M})[n_{M}(n_{M}+1)(2n_{M}+1)/6] + (n_{M}+1)^2 m\} N_a(N_a-1)/2$, with order of magnitude: $\mathcal{O}((n+m) n_{M}^2 N_a^2)$. Whereas, between active and permanently passive agents is: $[(n/n_{M})(\sum_{i=1}^{n_{M}}i)+(n_{M}+1)m] N_a N_p = \{(n/n_{M})[n_{M}(n_{M}+1)/2]+(n_{M}+1)m\} N_a N_p$, with order of magnitude: $\mathcal{O}((n+m) n_{M} N_a N_p)$. Using the $\mathcal{M}^*$ formulation, the total number of IC-dependent functions to be constrained between active agents is: $n_M^2 N_a(N_a-1)/2$. Whereas, between active and permanently passive agents it is: $n_{M} N_a N_p$.
\begin{proof}
Looking at Figure \ref{Fig_CAvsPS_2}, in between each of the instants $t_i \in [t_0, t_f]$, $i = 1, \dots n_{M}$, of  planned impulsive control application, the number of separations' combinations to be constrained is $(n/n_{M}) i^2 N_a(N_a-1)/2$ between active agents, and $(n/n_{M}) i N_a N_p$ between active and permanently passive agents. The factor $(n/n_{M})$ accounts for the need of guaranteeing safe separation in between the impulses, too. Instead on $(t_f, t_f + T]$, at each of the $m$ instants, the number of separation combinations to be constrained is $(n_{M}+1)^2 N_a(N_a-1)/2$ between active agents, and $(n_{M}+1) N_a N_p$ between active and permanently passive ones. The $+1$ in the $n_M+1$ factor accounts for the case in which control is not lost on $[t_0, t_f]$, but the reached target configuration has to be guaranteed safe with respect to the uncontrolled trajectory of another agent. This implies that the total number of separation combinations to be constrained over $[t_0, t_f +T]$ to have guarantee of PS between active agents is: $[(n/n_{M})(\sum_{i=1}^{n_{M}}i^2)+(n_{M}+1)^2 m] N_a(N_a-1)/2$, with order of magnitude: $\mathcal{O}((n+m) n_{M}^2 N_a^2)$. Whereas, between active active and permanently passive agents is: $[(n/n_{M})(\sum_{i=1}^{n_{M}}i)+(n_{M}+1)m] N_a N_p$, with order of magnitude: $\mathcal{O}((n+m) n_{M} N_a N_p)$. 
The total number of IC-dependent functions to be constrained equates the total number of IC states' combinations involved in Eq. \ref{PS_M_star},  which equates the number of contingency instants' combinations. Therefore, to enforce PS between active agents, the total number of IC states' combinations to be constrained is: $n_M^2 N_a(N_a-1)/2$. Whereas, between active active and permanently passive agents is: $n_M N_a N_p$. This concludes the proof.
\end{proof}

\subsection{Relative Motion Dynamics}\label{App:Dyn}
 The map between the integration constants (IC) used in this paper and the IC of the Yamanaka-Ankersen (YA) solution of relative motion in eccentric orbits is
\begin{equation}\label{ICYA}
\begin{aligned}
\boldsymbol{M}_{YA}(\textbf{\oe}) = \left[ \begin{array}{cccccc} 0 & 0 & 0 & -2 & 0 & 0 \\ 1 & 0 & 0 & 0 & 0 & 0 \\ 0 & -\sin(\omega) & \cos(\omega) & 0 & 0 & 0 \\ 0 & \cos(\omega) & \sin(\omega) & 0 & 0 & 0 \\ 0 & 0 & 0 & 0 & -\cos(\omega) & -\sin(\omega) \\ 0 & 0 & 0 & 0 & -\sin(\omega) & \cos(\omega) \\ \end{array} \right]
\end{aligned}
\end{equation}
Moreover, the IC used in this paper are linked to the quasi-nonsingular ROE through the first order map
\begin{equation}\label{MapICtoROE}
\resizebox{0.75\linewidth}{!}{$
\textbf{c}_{k} \approx \boldsymbol{M}_{CW}(\textbf{\oe})  \left[ \begin{array}{cccccc} 1 & 0 & 0 & 0 & 0 & 0 \\ \frac{3}{2} n t & 1 & 0 & 0 & 0 & 0 \\ 0 & 0 & 1 & 0 & 0 & 0 \\ 0 & 0 & 0 & 1 & 0 & 0 \\ 0 & 0 & 0 & 0 & 1 & 0 \\ 0 & 0 & 0 & 0 & 0 & 1 \\ \end{array} \right] \left( \begin{array}{c} (a_{k} - a)/a \\ (\omega_{k} + M_k(t) - \omega - M(t) ) + (\Omega_{k} - \Omega) \cos(i) \\ e_{k} \cos(\omega_k) - e \cos(\omega) \\ e_{k} \sin(\omega_k) - e \sin(\omega) \\ i_k - i \\ (\Omega_k - \Omega) \sin(i) \end{array} \right)$}
\end{equation} 
where
\begin{equation}\label{ICHCW}
\resizebox{0.85\linewidth}{!}{$
\begin{aligned}
\boldsymbol{M}_{CW}(\textbf{\oe}) = \left[ \begin{array}{cccccc} \eta^2 & 0 & 0 & 0 & 0 & 0 \\ 0 & \frac{1}{\eta} & - \Big(\eta^2 - \frac{1}{\eta}\Big) \frac{\sin(\omega)}{e} & \Big(\eta^2 - \frac{1}{\eta}\Big) \frac{\cos(\omega)}{e} & 0 & \Big(\eta^2 - \frac{1}{\eta}\Big) \cot(i) \\ 0 & \frac{e \sin(\omega)}{\eta} & \cos^2(\omega) + \frac{\sin^2(\omega)}{\eta} & \cos(\omega)\sin(\omega) - \frac{\cos(\omega)\sin(\omega)}{\eta} & 0 & - \frac{e \sin(\omega) \cot(i)}{\eta} \\ 0 & - \frac{e \cos(\omega)}{\eta} & \cos(\omega) \sin(\omega) - \frac{\cos(\omega)\sin(\omega)}{\eta} & \sin(\omega)^2 + \frac{\cos(\omega)^2}{\eta} & 0 & \frac{e \cos(\omega) \cot(i)}{\eta} \\ 0 & 0 & 0 & 0 & \eta^2 & 0 \\ 0 & 0 & 0 & 0 & 0 & \eta^2 \\
\end{array} \right]
\end{aligned}$}
\end{equation}
and, in the ROE, $M(t)$ is the mean anomaly at instant $t$. The center matrix in Eq. \ref{MapICtoROE} back-propagates the ROE to initial integration instant ($t = 0$) through the integrable dynamics. For $e \rightarrow 0$, $\boldsymbol{M}_{CW}$ becomes identity (the singularities cancel out), and the IC used in this paper reduce to the quasi-nonsingular ROE and to the  Clohessy-Wiltshire (CW)'s IC \cite{damico_ROEIC_2005}.

\section*{Acknowledgments}
This work was supported by the William R. and Sara Hart Kimball Stanford Graduate Fellowship, and by the VISORS's mission NSF Award \#1936663. The authors are thankful for their support.

\bibliography{main}

\end{document}